\documentclass[12pt]{amsbook}
\usepackage[utf8]{inputenc}
\usepackage{ae,aecompl,aeguill}	
\usepackage[draft]{hyperref} 
\usepackage[all]{hypcap} 

\usepackage{xspace}
\usepackage{ifthen}
\usepackage[francais]{babel}
\usepackage[babel=true,kerning=true]{microtype} 

\usepackage[a4paper,asymmetric,bindingoffset=1.2cm,margin=1.5cm,top=2.5cm,bottom=2.5cm]{geometry}
\makeatletter
\renewcommand{\section}{\@startsection{section}{1}%
  \z@{1.2\linespacing\@plus\linespacing}{\linespacing}%
  {\normalfont\bfseries\centering}}
\def\@makechapterhead#1{\global\topskip 3.5pc\relax
  \begingroup
  \fontsize{\@xivpt}{18}\bfseries\centering
    \ifnum\c@secnumdepth>\m@ne
      \leavevmode \hskip-\leftskip
      \rlap{\vbox to\z@{\vss
          \centerline{\normalsize\mdseries
              \uppercase\@xp{\chaptername}\enspace\thechapter}
          \vskip 3pc}}\hskip\leftskip\fi
     #1\par \endgroup
  \skip@34\p@ \advance\skip@-\normalbaselineskip
  \vskip\skip@ }
\def\@makeschapterhead#1{\global\topskip 3.5pc\relax
  \begingroup
  \fontsize{\@xivpt}{18}\bfseries\centering
  #1\par \endgroup
  \skip@34\p@ \advance\skip@-\normalbaselineskip
  \vskip\skip@ }

\makeatother

\usepackage[usenames,svgnames]{xcolor}
\usepackage{tikz}
\usetikzlibrary{matrix}
\usetikzlibrary{shapes.multipart}
\usepackage{graphicx}
\usepackage[matrix,arrow,curve]{xy}
\graphicspath{{Fig/}}
\long\def\psboxit#1#2{%
\begingroup\setbox0=\hbox{#2}%
\dimen0=\ht0 \advance\dimen0 by \dp0%
    \hbox{%
    \copy0%
    }
\endgroup%
}
\def\Gbox#1{\psboxit{box 0.7 setgray fill}{#1}}
\def\SetTableau#1#2#3#4{%
  \gdef\Tabvrule{\vrule\vrule width-0.4pt}
  \gdef\Tabhrule{\hrule\hrule height-0.4pt}  
  \gdef\Tabstrut{\vrule height#1 depth#2 width0pt\relax}
  \gdef\Tabbox##1{\hbox to #3{\hskip0.4pt\hfill\Tabstrut$#4##1$\hfill}}
} 
\def\PetitTableau{\SetTableau{1.65ex}{0.55ex}{2.2ex}{\scriptstyle}}

\def\Case#1{\vcenter{\Tabhrule%
                   \hbox{\Tabvrule\Tabbox#1\Tabvrule}\Tabhrule}}

\def\CaseGrise#1{\omit\Gbox{$\Case{#1}$}}
\def\GenTab#1{\vcenter{\halign{&$\Case{##}$\cr#1}}\egroup}
\def\Tableau{%
  \bgroup%
  \let\ =\omit%
  \let\\=\cr%
  \offinterlineskip\GenTab}

\usepackage{array}
\usepackage{rotating}
\usepackage{amssymb}
\usepackage{enumerate}
\usepackage{bibunits}
\usepackage{chngcntr}

\counterwithout{figure}{chapter}

\title[Algèbre Combinatoire et Effective]{      
  Algèbre combinatoire et effective:
  des graphes aux algèbres de Kac
  \emph{via} l'exploration informatique}
\author{Nicolas M. Thiéry}

\newcommand{\mupadcombinat}{\href{http://mupad-combinat.sf.net/}{\texttt{MuPAD-Combinat}}\xspace}
\newcommand{\starcombinat}{\href{http://wiki.sagemath.org/combinat/}{{\raisebox{-.1ex}{$\ast$}\text{-Combinat}}}\xspace}
\newcommand{\starcombinatsimple}{\href{http://wiki.sagemath.org/combinat/}{$\ast$\text{-Combinat}}\xspace}
\newcommand{\aldorcombinat}{\href{http://www.risc.uni-linz.ac.at/people/hemmecke/aldor/combinat/}{\texttt{aldor-Combinat}}\xspace}
\newcommand{\sagecombinat}{\href{http://wiki.sagemath.org/combinat/}{\texttt{Sage-Combinat}}\xspace}
\newcommand{\mupad}{\href{http://www.mupad.de/}{\texttt{MuPAD}}\xspace}
\newcommand{\maple}{\href{http://www.maplesoft.com/}{\texttt{Maple}}\xspace}
\newcommand{\magma}{\href{http://magma.maths.usyd.edu.au/magma/}{\texttt{Magma}}\xspace}
\newcommand{\aldor}{\href{http://www.aldor.org/}{\texttt{Aldor}}\xspace}
\newcommand{\axiom}{\href{http://axiom-wiki.newsynthesis.org/}{\texttt{Axiom}}\xspace}
\newcommand{\permuvar}{\href{http://permuvar.sourceforge.net}{\texttt{PerMuVAR}}\xspace}
\newcommand{\gap}{\href{http://www.gap-system.org/}{\texttt{GAP}}\xspace}
\newcommand{\java}{\href{http://www.java.com/}{\texttt{Java}}\xspace}
\newcommand{\perl}{\href{http://www.perl.com/}{\texttt{Perl}}\xspace}
\newcommand{\ruby}{\href{http://www.ruby-lang.org/}{\texttt{Ruby}}\xspace}
\newcommand{\sage}{\href{http://www.sagemath.org/}{\texttt{Sage}}\xspace}
\newcommand{\mathematica}{\href{http://www.wolfram.com/}{\texttt{Mathematica}}\xspace}
\newcommand{\matlab}{\href{http://www.mathworks.com/}{\texttt{Matlab}}\xspace}
\newcommand{\ace}{\href{http://phalanstere.univ-mlv.fr/~ace/}{\texttt{ACE}}\xspace}
\newcommand{\muec}{\href{http://igm.univ-mlv.fr/LabInfo/equipe/combinatoire/MUEC/}{$mu$-\texttt{EC}}\xspace}
\newcommand{\SF}{\href{http://www.math.lsa.umich.edu/~jrs/maple.html}{\texttt{SF}}\xspace}

\newcommand{\End}{{\operatorname{End}}}

\newcommand{\Inv}[1][n]{\mathcal{I}_{#1}} 
\newcommand{\RR}{\mathbb{R}}

\newcommand{\NN}{\mathbb{N}}
\newcommand{\CC}{\mathbb{C}}
\newcommand{\ZZ}{\mathbb{Z}}
\newcommand{\QQ}{\mathbb{Q}}

\newcommand{\age}{{\mathcal A}}
\newcommand{\agealgebra}{{\QQ.\mathcal A}}
\newcommand{\profile}{\varphi}
\newcommand{\id}{{\operatorname{id}}}
\newcommand{\pr}{\mathrm{pr}}
\newcommand{\Des}{{\operatorname{D}_R}}
\newcommand{\Rec}{{\operatorname{D}_L}}

\newcommand{\coweightspace}[1][]{\mathfrak{h}_{#1}}
\newcommand{\coweightlattice}{\coweightspace[\ZZ]}

\newcommand{\q}{-\frac{q_1}{q_2}}

\newcommand{\height}{\operatorname{ht}}

\newcommand{\x}{x}

\newcommand{\coroot}{\alpha^\vee}
\newcommand{\lc}{{\lambda^\vee}}

\newcommand{\rhoc}{{\rho^\vee}}
\newcommand{\Lambdac}{\Lambda^\vee}

\newcommand{\cl}{{\operatorname{cl}}}

\newcommand{\opi}{{\overline{\pi}}}
\newcommand{\W}{W}
\newcommand{\clW}{\mathring{W}}
\newcommand{\Wa}[1][\CC]{{#1[\W]}}

\newcommand{\kW}[1][\CC]{{#1\W}}
\newcommand{\kclW}[1][\CC]{{#1\clW}}
\newcommand{\heckeW         }[2][\W]{{           \operatorname{H} (#1)(#2)}}
\newcommand{\heckeWW        }[1][\W]{{           \operatorname{H}\!#1     }}

\newcommand{\Wmax}{{w_0}}

\newcommand{\sg}[1][n]{{\mathfrak{S}_{#1}}}
\newcommand{\sga}[2][\CC]{{#1[\sg[#2]]}}
\newcommand{\ksg}[2][\CC]{{#1\sg[#2]}}
\newcommand{\hecke}[2][n]           {{\operatorname{H}_{#1}(#2)}}
\newcommand{\heckesg}[1][n]         {{\heckeWW[\mathfrak{S}_{#1}]}}
\newcommand{\affinehecke}[2][n]     {{\widetilde{\operatorname{H}}_{#1}(#2)}}

\newcommand{\ndf}[1]{\operatorname{NDF}_{#1}}
\newcommand{\ndfa}[2][\CC]{{#1[\ndf{#2}]}}
\newcommand{\ndpf}[1]{\operatorname{NDPF}_{#1}}
\newcommand{\ndpfa}[2][\CC]{{#1[\ndpf{#2}]}}

\newcommand{\suchthat}{{\ |\ }}

\newcommand{\KD}{K\!D}

\makeatletter
\newskip\@bigflushglue \@bigflushglue = -100pt plus 1fil

\def\bigcentering{\let\\\@centercr\rightskip\@bigflushglue%
\leftskip\@bigflushglue
\parindent\z@\parfillskip\z@skip}

\makeatother

\newboolean{draft}
\setboolean{draft}{false}
\ifdraft
\newcommand{\TODO}[2][To do: ]{\textcolor{red}{\textbf{#1#2}}}
\else
\newcommand{\TODO}[2][]{}
\fi
\newcommand{\FIXME}{\TODO[Fix-me: ]}

\newcommand{\dynkin}[1]{
  \begin{tikzpicture}[>=latex,join=bevel,baseline=(current bounding box.east)]
    #1
  \end{tikzpicture}}
\newcommand{\dynkinAIIa}{\dynkin{
    \node (N_0) at (21bp,22bp) [draw,draw=none] {$0$};
    \node (N_1) at (6bp,12bp) [draw,draw=none] {$1$};
    \node (N_2) at (36bp,12bp) [draw,draw=none] {$2$};
    \draw [] (N_0) -- (N_1);
    \draw [] (N_1) -- (N_2);
    \draw [] (N_2) -- (N_0);
  }}
\newcommand{\dynkinAVa}{\dynkin{
    \node (N_1) at (46bp,19bp) [draw,draw=none] {$0$};
    \node (N_2) at (86bp,12bp) [draw,draw=none] {$5$};
    \node (N_3) at (6bp,12bp) [draw,draw=none] {$1$};
    \node (N_4) at (26bp,7bp) [draw,draw=none] {$2$};
    \node (N_5) at (46bp,6bp) [draw,draw=none] {$3$};
    \node (N_6) at (66bp,8bp) [draw,draw=none] {$4$};
    \draw [] (N_1) -- (N_2);
    \draw [] (N_3) -- (N_1);
    \draw [] (N_3) -- (N_4);
    \draw [] (N_4) -- (N_5);
    \draw [] (N_5) -- (N_6);
    \draw [] (N_6) -- (N_2);
  }}
\newcommand{\dynkinCIIa}{\dynkin{
    \node (N_0) at (6bp,6bp) [draw,draw=none] {$0$};
    \node (N_1) at (36bp,6bp) [draw,draw=none] {$1$};
    \node (N_2) at (66bp,6bp) [draw,draw=none] {$2$};
    \draw [->,double] (N_0) -- (N_1);
    \draw [<-,double] (N_1) -- (N_2);
  }}
\newcommand{\dynkinGIIa}{\dynkin{
    \node (N_1) at (10bp,6bp) [draw,draw=none] {$0$};
    \node (N_2) at (36bp,6bp) [draw,draw=none] {$1$};
    \node (N_3) at (66bp,6bp) [draw,draw=none] {$2$};
    \draw [] (N_1) -- (N_2);
    \draw [->] (N_2) -- node[above] {$3$} (N_3);
  }}

\newtheorem{theorem}{Théorème}[section]
\newtheorem{lemma}[theorem]{Lemme}
\newtheorem{proposition}[theorem]{Proposition} 
\newtheorem{corollary}[theorem]{Corollaire} 

\newtheorem{definition}[theorem]{Définition}
\newtheorem{problem}[theorem]{Problème}

\newtheorem{conjecture}[theorem]{Conjecture}
\theoremstyle{remark}

\newtheorem{question}[theorem]{Question}

\usepackage{verbatim}
\usepackage{moreverb}
\makeatletter
\g@addto@macro\verbatim{\microtypesetup{kerning=false}}
\def\Mexin@processline{>{}> \the\verbatim@line\par}
\newenvironment{Mexin} {\vspace{-.5ex}\verbatim\small\addtolength\parskip{-.5ex}\let\verbatim@processline=\Mexin@processline}{\endverbatim}
\newenvironment{Mexout}{\vspace{-.5ex}\verbatim\tiny\addtolength\parskip{-.9ex}}{\endverbatim}
\makeatother
\newcommand{\Mup}[1]{\texttt{\microtypesetup{kerning=false}#1}}

\begin{document}
\enlargethispage{10cm}
\thispagestyle{empty}

\begin{center}
  \large
%
  \includegraphics[height=2cm]{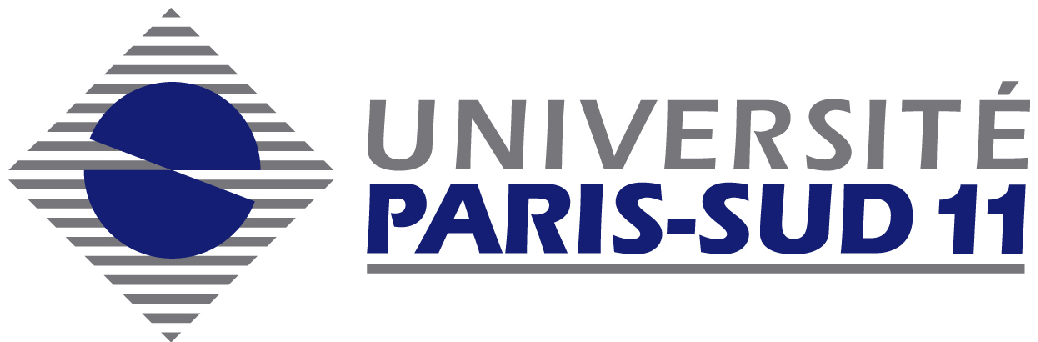}
  \hfill \includegraphics[height=2cm]{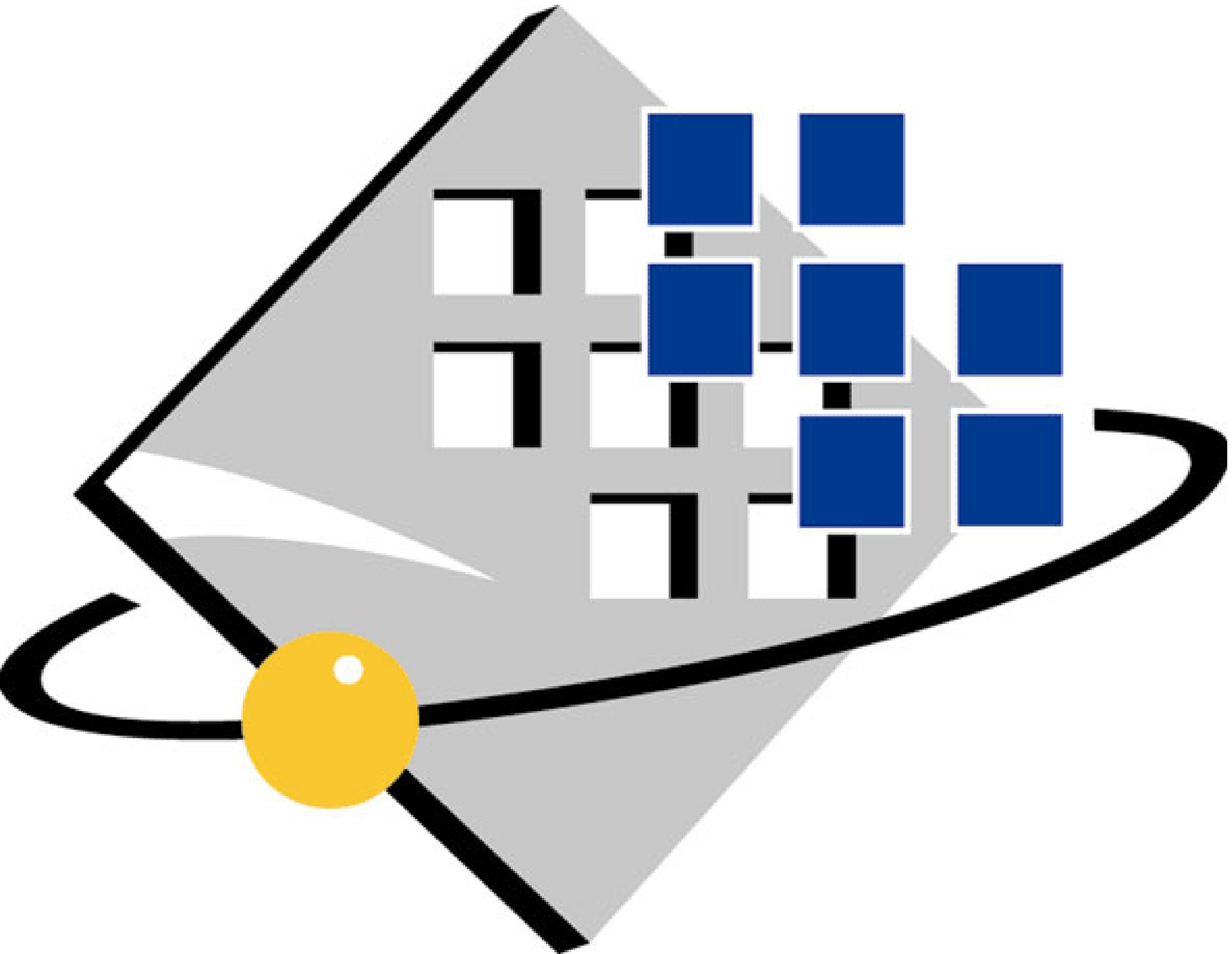}
  \bigskip\\
  \framebox[14cm]{
    \vbox{
      \vbox to .5 ex{}
      {
      \LARGE
      UNIVERSITÉ PARIS-SUD\\
      FACULTÉ DES SCIENCES D’ORSAY}\\
      \vbox to .5ex{}
    }
  }
  \ \\
  \bigskip
  {\Large MÉMOIRE}\\
  \bigskip
  Présenté pour obtenir\\
  \bigskip
  {\Large
    LE DIPLÔME\\
    D'HABILITATION À DIRIGER DES RECHERCHES\\
    DE L'UNIVERSITÉ PARIS XI\\}
  \bigskip
  Spécialité: Mathématiques\\
  \bigskip
  par\\
  \bigskip
  \textbf{{\LARGE Nicolas M. THIÉRY}}\\
  \bigskip
  \bigskip
  \hrule
  \bigskip
  {\textbf{
      \huge
      Algèbre combinatoire et effective:\\
      des graphes aux algèbres de Kac\\
      \emph{via} l'exploration informatique\\}}
  \bigskip
  \hrule
  \bigskip
  \bigskip
  Soutenu le 10 décembre 2008\\
  devant la commission d'examen:\\
  \bigskip
  \bigskip
  \begin{tabular}{ll}
    François Bergeron		& Rapporteur\\
    Peter Cameron		& Rapporteur\\
    Bernard Leclerc		& Rapporteur\\
    \\
    Jean-Benoît Bost		& Examinateur\\
    Mireille Bousquet-Mélou	& Examinatrice\\
    Alain Lascoux		& Examinateur\\
    Jean-Yves Thibon		& Examinateur\\
    Leonid Vainerman		& Examinateur\\
    Paul Zimmermann		& Examinateur\\
  \end{tabular}
\end{center}
  %
  %

\clearpage

\TODO{dimension finie  /  snob  /  }


\setcounter{page}{0}
\setcounter{tocdepth}{1}
\tableofcontents
\listoffigures


\makeatletter
\newenvironment{publist}[2]{%
  \begingroup
  \renewcommand{\bibliofont}{\normalsize}
  \makeatletter
  \def\subsection{\@startsection{subsection}{1}%
    \z@{.7\linespacing\@plus\linespacing}{.5\linespacing}%
    {\normalfont\bfseries}}
  \makeatother
  \subsection*{#1 }%
  \endgroup
  \begin{bibunit}[unsrt]%
    \renewcommand{\@bibunitname}{publist-#2}
    \renewcommand{\section}[2]{}
    \renewcommand{\chapter}[2]{}
  }{
    \putbib[isil]%
  \end{bibunit}}
\makeatother

\chapter*{Liste de publications}


\begin{publist}{Articles dans des revues d'audience
    internationale avec comité de rédaction}{journeaux}
  \nocite{Hivert_Schilling_Thiery.HeckeGroupAffine.2008}
  \nocite{Hivert_Schilling_Thiery.HeckeGroupAffine.2007}
  \nocite{Hivert_Thiery.HeckeGroup.2007}
  \nocite{Gaudry_Schost_Thiery.2004}
  \nocite{Hivert_Thiery.MuPAD-Combinat.2004}
  \nocite{Novelli_Thibon_Thiery.2004}
  \nocite{Thiery_Thomasse.SAGBI.2002}
  \nocite{Hivert_Thiery.SA.2002}
  \nocite{Pouzet_Thiery.IAGR.2001}
  \nocite{Thiery.AIG.2000}
\end{publist}
\begin{publist}{Articles dans des actes de conférences internationales
    avec comité de rédaction}{confs}
  \nocite{Hivert_Thiery.HeckeSg.2006}
  \nocite{Pouzet_Thiery.AgeAlgebra.2005}
  \nocite{Thiery.CMGS.2001}
\end{publist}
\clearpage
\begin{publist}{Articles soumis}{soumises}
  \nocite{Bandlow_Schilling_Thiery.2008.Promotion}
  \nocite{David_Thiery.2008.Kac}
\end{publist}
\begin{publist}{Articles en préparation}{preparation}
  \nocite{Pouzet_Thiery.IAGR}
  \nocite{Pouzet_Thiery.AgeAlgebra1}
  \nocite{Pouzet_Thiery.AgeAlgebra2}
\end{publist}
\begin{publist}{Thèse de doctorat et autres communications}{autres}
  \nocite{Martinez_Molinero_Thiery.2006}
  \nocite{Thiery.PDemo.2000}
  \nocite{Thiery.AIG.1999}
  \nocite{Thiery.IAGR}
\end{publist}
\newpage
\chapter*{Remerciements}
\enlargethispage{10cm}
\thispagestyle{empty}

Je suis très honoré que François Bergeron, Bernard Leclerc et Peter
Cameron aient accepté d'être rapporteurs de ce mémoire. Les nombreuses
discussions que j'ai eues avec eux, lors de diverses conférences, ont
toujours été captivantes et éclairantes.

L'influence de Jean-Yves Thibon sur mes travaux est évidente dans ce
mémoire; les problèmes qu'il a soulevés sont à l'origine d'un chapitre
entier. Celle d'Alain Lascoux est moins visible. J'aimerais d'autant
plus la faire apparaître à sa juste valeur: certes nous, les «gamins»
du phalanstère, n'appliquons pas souvent à la lettre ses nombreux
conseils avisés; c'est sans doute par besoin de les mettre à notre
sauce, à notre niveau, pour mieux les intégrer. Mais la démarche est
là, profondément ancrée, et c'est grâce à elle que nous avons abouti à
bien des résultats de ce mémoire. Après m'avoir soutenu, et avoir
épaulé sans réserve le projet \starcombinat depuis sa création,
Jean-Yves et Alain ont naturellement accepté d'être dans mon
jury. Pour tout cela, et pour le phalanstère qu'ils font vivre, merci.

Je voudrais exprimer ma gratitude envers Marie-Claude David et Léonid
Vainerman pour m'avoir guidé dans la découverte des facteurs et des
algèbres de Kac et envers Léonid pour sa participation à mon
jury. Je tiens aussi à remercier Paul Zimmermann, pour son soutien et
nos multiples échanges et actions autour du calcul formel libre, en
particulier pour la combinatoire. Enfin, c'est un plaisir et un
honneur de compter Mireille Bousquet-Mélou dans mon jury.

Cette habilitation doit beaucoup à ces deux dernières années
consacrées à la recherche; merci à la NSF, au MSRI et au CNRS pour le
financement, et surtout à tous mes collègues de l'IUT d'Orsay; je sais
l'effort qu'ils ont dû consentir pour cela. Plus généralement, je
voudrais remercier toutes les équipes qui m'ont toujours accueilli à
bras ouverts, et auprès desquelles je me suis formé: le Laboratoire de
Probabilités, Combinatoire et Statistiques de Lyon I (Bernard Roux,
Stéphan Thomassé, etc.), les laboratoires GAGE et LIX de l'École
Polytechnique (Marc Giusti, Éric Schost, Pierrick Gaudry), le
phalanstère de combinatoire, à Rouen et Marne-la-Vallée
(Jean-Christophe Novelli, Florent Hivert, Teresa Gomez-Diaz, etc.),
l'équipe \mupad à Paderborn (Christopher Creutzig, Benno Fuchssteiner,
Ralf Hillebrand, Walter Oevel, Stefan Wehmeier, etc.), les
départements de mathématiques de l'Université de Californie à San
Diego (Adriano Garsia, Nolan Wallach) et à Davis (Anne Schilling,
Brant Jones, Jason Bandlow, Monica Varizani, Jesus de Loera, etc.), le
projet NSF Affine Schubert Calculus (Grant No 0652641, Anne Schilling,
Luc Lapointe, Jennifer Morse, Mark Schimozono, Mike Zabrocki, Thomas
Lam). Un remerciement tout particulier va à l'équipe d'arithmétique et
géométrie algébrique d'Orsay, qui m'offre chaleureusement depuis
quatre ans un espace de liberté où seule la production scientifique
compte, quels que soient ma thématique et mes collaborateurs.  Merci
spécifiquement à Jean-Benoît Bost pour sa participation à mon jury et
à David Harari pour l'accompagnement de cette habilitation.


Chers codéveloppeurs de \starcombinat, chers coauteurs, Anne, Conrado,
Éric, Florent, Jason, Jean-Yves, Jean-Christophe, Marie-Claude,
Pierrick, Stéphan, Xavier, grand merci pour ce bout de chemin
ensemble, hier, aujourd'hui et demain. Au fond, c'est ce chemin, plus
que la destination, qui donne tout son sens humain à notre travail.

Cher Maurice, ce mémoire vous est dédié. Vous m'avez pris sous votre
aile et appris à être libre et à faire miens les vers de Cyrano:
\begin{quotation}
  «Ne pas monter bien haut, peut-être, mais tout seul!»
\end{quotation}
De tout cœur, merci.

\newpage

Et puis ... Ce mémoire ne serait pas ce qu'il est sans le courage de
ses relecteurs attentifs et dévoués: Albane, Corinne, Émilia, Florent,
Jean, Sandrine. Il n'existerait tout simplement pas sans la patience
d'Adèle et d'Élise pendant que papa écrivait «son petit bouquin» et
sans le soutien indéfectible et réconfortant de leur maman.

\chapter*{Prélude}


Voilà venu le temps de l'habilitation, où je suis censé démontrer que,
non content d'avancer vaille que vaille dans ma propre recherche, je
peux prendre quelqu'un sous mon aile. Ce présent mémoire a donc
vocation à répondre à deux questions: «qui suis-je?» et «où vais-je?».

«Mais sur quoi travaillez-vous au juste?», me demande-t-on. S'il faut
vraiment me mettre dans une case, je choisis la combinatoire
algébrique; ou peut-être plutôt l'algèbre combinatoire.  «Mais
encore?» Avec un néophyte, je peux répondre que c'est l'art de compter
en utilisant les miraculeuses propriétés de l'addition et de la
multiplication; «c'est bien, mais compter quoi? Et pour quoi faire?»
J'y reviendrai. En revanche, avec un collègue qui aimerait me ranger
dans une petite case, je suis plus perplexe. La
figure~\ref{figure.topics} traduit au mieux l'image que j'ai en
tête. Théorie des graphes?  Des invariants? Des représentations?
Combinatoire bordelaise? Fonctions symétriques? Algèbres de Hopf,
voire de Kac?  Calcul formel? Un peu de tout cela; mais je n'ose
répondre oui à aucune de ces questions, de peur d'y passer pour un
béotien. Je ne pourrais pas, à l'instar de nombreux collègues que
j'envie pour l'occasion, écrire un grand mémoire de synthèse résumant
leurs vastes connaissances sur \emph{leur} sujet et ouvrant de grandes
portes sur un avenir radieux pour toute une communauté derrière eux.

Suis-je un rêveur éclectique ne sachant rien sur tout? Peut-être.
J'espère que ce mémoire, présentant mes contributions à ces différents
sujets de recherche, mettra en valeur ce qui les unit: outils,
méthodologies, points de vue. 

Que suis-je donc. Avant tout un explorateur. Foin de l'image
romantique de l'archéologue devinant la huitième merveille du monde à
partir de quelques tessons épars et établissant un délicat plan de
fouilles sur vingt ans pour la mettre à jour millimètre par millimètre
au pinceau. Non, moi c'est plutôt la dynamite et la tronçonneuse, le
GPS et les drones télécommandés. En clair l'exploration informatique,
lorsqu'elle s'y prête. Pragmatisme et efficacité. Je suis peut-être
bête, mais mon marteau-piqueur est plus gros que le tien, et je sais
le manier. Heureusement que le champ des idées est renouvelable à
l'infini, sinon j'aurais changé de métier depuis belle lurette (déjà
que j'ai mauvaise conscience pour le CO${}_2$ relâché à l'occasion de
mes voyages professionnels ou pour faire tourner mes calculs).

Ma stratégie favorite est simple. Choisir une belle montagne perdue
dans la brume et chercher un guide pour faire équipe. Ou l'inverse. Au
spécialiste du domaine de me transmettre sa science, de m'expliquer
les subtilités locales du climat, de m'emmener à la frontière du connu
(un problème précis et l'état de l'art environnant). À moi de sortir
la tronçonneuse pour défricher (déchiffrer?) la forêt vierge,
d'inventer une nouvelle machine à chasser les nuages. Où sont les
obstacles? Les abîmes? Les sommets les plus abrupts? N'y aurait-il pas
une vallée suspendue pour traverser le massif en douceur, un pont
provisoire pour enjamber le ravin? Nommer, aplanir grossièrement et
cartographier. Définir, remarquer et conjecturer. Puis passer à autre
chose. Aux alpinistes des théorèmes de vaincre les plus hauts
sommets. Aux bâtisseurs de théories de faire un jardin japonais du
champ de bataille jonché de faits incongrus que je laisse derrière
moi.

J'apporte un savoir faire et une caisse à outils. En échange, mon
équipier apporte un savoir. Un savoir que j'intégrerai d'autant mieux
qu'il aura fallu que je l'implante. Pas de flou artistique
possible. Ici, la collaboration interdisciplinaire fonctionne, car
elle est emmenée par une succession de questions concrètes, évitant le
piège des considérations générales et oiseuses. Et à la fin, grâce au
logiciel libre, nous repartons tous les deux avec la caisse; caisse
qui à l'occasion s'est enrichie de nouveaux outils, peut-être même
d'un nouveau forgeron.
Et puis, ensemble, nous publions les résultats. Car il ne s'agit pas
de se faire esclave programmeur.
Je suis avant tout chercheur. Le travail d'implantation est
intéressant dans la mesure où un nouveau problème nécessite une large
part de conception et d'algorithmique nouvelle. Il y a aussi un
savoir-explorer. Quelle est la bonne question à poser? Où faut-il
passer en force, où en ratissant au peigne fin? Où tenter sa chance et
lancer un hameçon au hasard?  Quelle confiance accorder à des premiers
signes?  Et puis celui qui tient la tronçonneuse et avance devant a
toutes les chances d'être le premier à entrevoir de nouveaux
phénomènes. D'autant qu'avec le temps l'œil se forme.

L'exploration informatique n'est pas une idée nouvelle en combinatoire
algébrique. Schützenberger en fut un pionnier dès les années
1950. Maintenant quasiment tous les chercheurs y ont recours à un
moment ou un autre, que ce soit en tapant quelques commandes \maple,
ou en développant, trop souvent de manière isolée, des bibliothèques
de plusieurs dizaines voire centaines de milliers de lignes de
code. Mon rêve: mutualiser tout ces efforts de développement pour qu'à
la fin chacun ait à sa disposition les meilleurs outils, tout en
perdant moins de temps à faire de la technique. C'est tout le sens du
projet logiciel \starcombinat que j'ai lancé en 2000 avec Florent
Hivert et que je décrirai en détail dans le
chapitre~\ref{chapter.combinat}. La stratégie est de cristalliser une
communauté transversale autour de \starcombinat, en tissant peu à peu
des liens à l'échelle internationale.  Partis à deux, nous sommes
maintenant plus d'une vingtaine, avec 130\,000 lignes de code.  Cela
nécessite de puiser dans le savoir-faire des informaticiens: d'une
part, les outils et modèles de développement collaboratifs (par
ex. logiciel libre) et, d'autre part, les techniques de conception
(par ex. programmation orientée objet) adaptés à notre situation. En
bref, industrialiser le processus pour maîtriser le changement
d'échelle. Et au final permettre des calculs d'un niveau de complexité
supérieur, intégrant simultanément plusieurs techniques algorithmiques
(comme de l'algèbre linéaire creuse avancée, de l'élimination type
base de Gröbner et des calculs combinatoires sur des objets à
isomorphie près; voir section~\ref{section.invariants}), ou combinant
plusieurs constructions conceptuelles (dualité, tenseurs, changements
de base, torsion de (co)produits, sous-algèbres et quotients; voir
section~\ref{section.kac}).
Cet effort de mutualisation existe depuis longtemps dans d'autres
domaines (par exemple avec GAP pour la théorie des groupes). Fait
nouveau, il se met en place à l'échelle des mathématiques (par exemple
avec Sage). Ma modeste contribution est de faire avancer la situation
dans mon domaine.

On l'a dit, la combinatoire algébrique se prête en général bien à
l'exploration informatique. Mais pas toujours. Une bonne partie de mes
recherches (algèbres d'âge, théorie autour des invariants de groupes
de permutations) a été faite au tableau noir, avec une bonne
vieille craie. Cependant l'approche est restée la même: explorer des
exemples concrets, voir et comprendre ce qui se passe, puis abstraire
autant que faire se peut; bref tenter d'appliquer la maxime:
%
\emph{The art of doing mathematics consists in finding that special
  case which contains all the germs of generality.} --David Hilbert
Quoted in N Rose Mathematical Maxims and Minims (Raleigh N C 1988).

L'informatique permet d'abord l'étude d'exemples plus conséquents, ce
qui peut être essentiel lorsque les premiers exemples non triviaux ne
sont déjà plus traitables à la main. Mais transparaissent aussi en
filigrane des questions qui me tiennent à cœur: qu'est-ce qui est
calculable, en pratique? Cet objet mathématique, puis-je le modéliser
sur ma machine pour pouvoir ensuite lui poser des questions
intéressantes?  Jusqu'où peut-on aller avec l'exploration
informatique? J'ai ces questions en tête dès que j'aborde un nouveau
sujet. Cela offre un point de vue, certes forcément réducteur, mais
qui donne un angle d'attaque, un fil conducteur et une succession de
prises pour rentrer dans le sujet.  Cela sans jugement de valeurs ni
prétention à l'universalité: c'est ce qui fonctionne, pour moi et dans
une certaine gamme de problèmes.


\chapter*{Introduction}

Ce mémoire fait la synthèse de presque quinze années de recherche,
afin d'en dégager les perspectives. Ces années ont été pour moi une
période de grande liberté, pendant laquelle j'ai pris le temps de me
forger une voie et une démarche personnelle, à mi-chemin entre
l'informatique et les mathématiques.

Ce qui m'a attiré vers la combinatoire algébrique, c'est l'ouverture
vers d'autres disciplines, en mathématiques, en informatique, ou en
physique théorique.  Ma démarche est en effet de me construire, petit
à petit, une boîte à outils, en élargissant progressivement mon champ
de recherches et en tissant un réseau de collaborateurs dans des
communautés variées: théorie des invariants, des graphes, des groupes,
combinatoire algébrique ou non, calcul formel, etc.

Lorsque, au fil des rencontres scientifiques et des séjours,
j'envisage d'aborder un nouveau sujet de recherche, je me pose deux
questions: «Ma boîte à outils actuelle me donne-t-elle un point de vue
original, m'offrant une chance de voir ce que d'autres n'ont pas
encore vu?»; et «Quels outils et concepts vais-je apprendre, qui
seraient susceptibles de déclencher des progrès sur des sujets en
suspens?».

Pour aborder un nouveau sujet, je travaille systématiquement en
collaboration. Mon ou mes partenaires sont les garants de l'intérêt et
de l'originalité dans un domaine où je n'ai pas forcément encore de
recul. En retour, j'apporte des outils et une expertise. Mon fonds de
commerce est l'exploration de domaines relativement vierges:
construire et étudier des exemples, repérer des conjectures. Et
surtout, chercher le bon point de vue où les énoncés et, idéalement,
les démonstrations s'expriment simplement.

Le langage forme souvent une barrière de communication entre domaines
éloignés. C'est pourquoi, dans ma démarche, la combinatoire joue un rôle
essentiel pour modéliser simplement, et souvent de manière effective,
des problèmes en les abstrayant de leur contexte. Cela permet de nouer
de nouvelles collaborations sur des problèmes concrets et précis. Pour
la même raison, ma question favorite est: «Comment cela se
calcule?». Si je suis capable de retranscrire le problème dans un
ordinateur, c'est qu'aucune subtilité ne m'a échappé; les deux parties
sont bien sur la même longueur d'onde. La compréhension du contexte et
des motivations, le plus souvent essentielle pour parvenir à une
solution, vient ensuite naturellement au fur et à mesure de l'échange
qui se met en place.

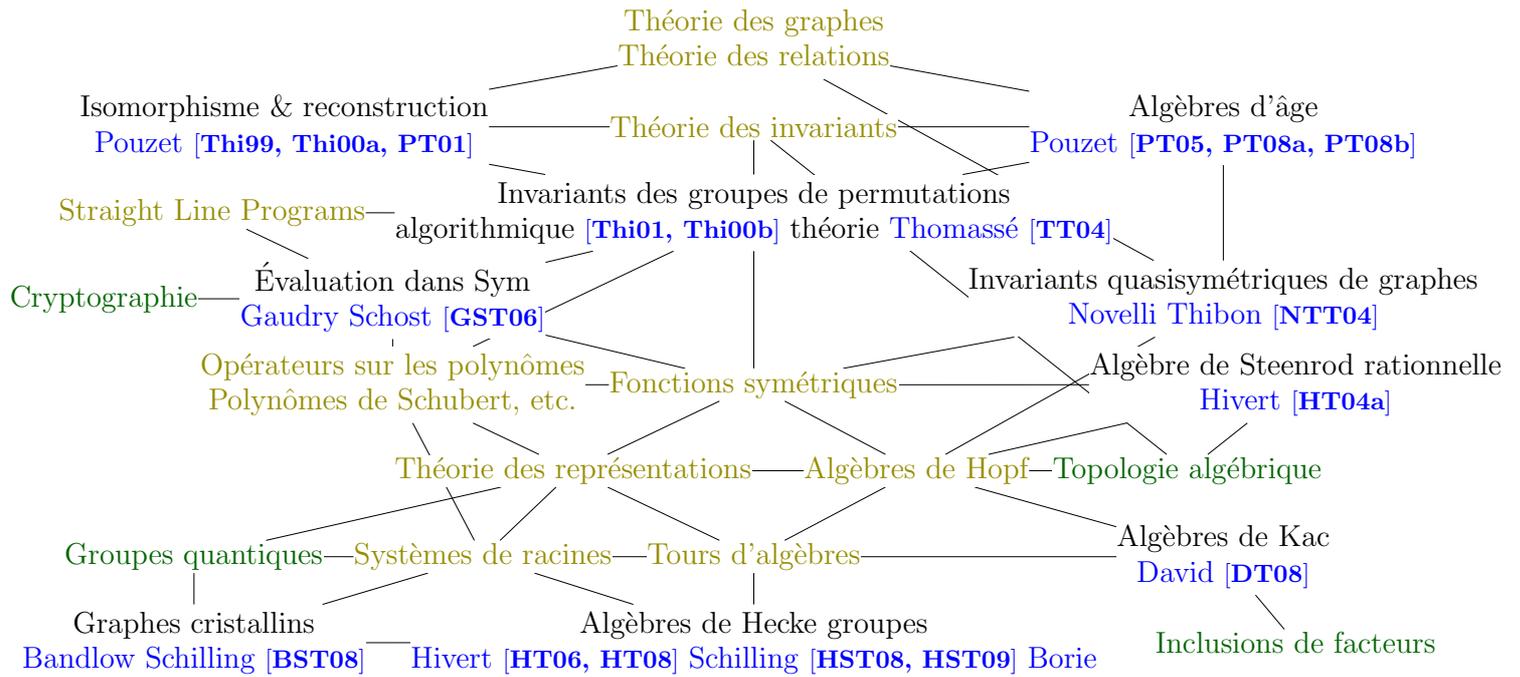
\begin{figure}[h]
  \begin{bigcenter}
    \scalebox{.95}{
            \pgfdeclarelayer{nodes}
      \pgfsetlayers{main,nodes}
      \newcommand{\fr}[1]{#1}
      \newcommand{\en}[1]{}
      \begin{tikzpicture}[>=latex,join=bevel,xscale=.5,yscale=.6]
        \newcommand{\mycite}[1]{{\Small{\cite{#1}}}}
        \tikzstyle{every node}=[fill=white,inner ysep=1pt, inner xsep=0pt]
        \tikzstyle{area}=[text=olive]
        \tikzstyle{topic}=[text centered]
        \tikzstyle{far area}=[text=DarkGreen,text centered]
        \newcommand{\mynode}[1]{\begin{tabular}{@{}c@{}}#1\end{tabular}}
        \newcommand{\coauthor}[1]{{\color{blue}{#1}}}
        \begin{pgfonlayer}{nodes}
          \node[area]     (graph)    at ( 0,20) {\mynode{
              \fr{Théorie des graphes\\Théorie des relations}
              \en{Graph Theory\\Theory of Relations}}};
          \node[area]     (inv)      at ( 0,18)         {
            \fr{Théorie des invariants}
            \en{Invariant Theory}};
          \node[area]     (sym)      at (  0,12)        {
            \fr{Fonctions symétriques}
            \en{Symmetric Functions}};
          \node[area]     (schubert) at (-10,12)        {\mynode{
              \fr{Opérateurs sur les polynômes\\Polynômes de Schubert, etc.}
              \en{Operators on Polynomials\\Schubert Polynomials, etc.}}};
          \node[topic]    (iso rec)  at (-13,18)         {\mynode{
              \fr{Isomorphisme \& reconstruction}
              \en{Isomorphism  \& Reconstruction}\\
              \coauthor{Pouzet~\mycite{Thiery.IAGR,Thiery.AIG.2000,Pouzet_Thiery.IAGR.2001}}}};
          \node[topic]    (age)      at ( 13,18)         {\mynode{
              \fr{Algèbres d'âge}
              \en{Age Algebras}\\
              \coauthor{Pouzet~\mycite{Pouzet_Thiery.AgeAlgebra.2005,Pouzet_Thiery.AgeAlgebra1,Pouzet_Thiery.AgeAlgebra2}}}};
          \node[topic]    (inv perm) at ( 0,16)         {\mynode{
              \fr{Invariants des groupes de permutations}
              \en{Invariants of Permutation Groups}\\
              \fr{algorithmique}\en{Algorithmic}
              \coauthor{\mycite{Thiery.CMGS.2001,Thiery.PDemo.2000}}
              \fr{théorie}\en{Theory} 
              \coauthor{Thomassé~\mycite{Thiery_Thomasse.SAGBI.2002}}}};
          \node[topic]    (sym eval) at (-10,14)         {\mynode{
              \fr{Évaluation dans Sym}
              \en{Evaluation properties of SF}\\
              \coauthor{Gaudry Schost~\mycite{Gaudry_Schost_Thiery.2004}}}};
          \node[far area] (crypto)   at (-18,14)        {
            \fr{Cryptographie}
            \en{Cryptography}};
          \node[area]     (slp)      at (-15,16)        {Straight Line Programs};
          \node[area]     (hopf)     at (4.5,10)        {\fr{Algèbres de Hopf}\en{Hopf algebras}};
          \node[area]     (RT)       at ( -5,10)        {\fr{Théorie des représentations}\en{Representation Theory}};
          \node[area]     (fd alg)   at (  0,8)         {\fr{Tours d'algèbres}\en{Fin. Dim. Algebras}};
          \node[area]     (roots)    at ( -7.5,8)         {\fr{Systèmes de racines}\en{Root Systems}};
          \node[far area] (lie)      at (-15.5,8)       {\fr{Groupes quantiques}\en{Quantum groups}};
          \node[topic]    (crystals) at (-15.5, 6)      {\mynode{
              \fr{Graphes cristallins}
              \en{Crystals}\\
              \coauthor{Bandlow Schilling~\mycite{Bandlow_Schilling_Thiery.2008.Promotion}}}};
          \node[topic]    (hecke)    at (  0, 6)        {\mynode{
              \fr{Algèbres de Hecke groupes}
              \en{Hecke Group Algebras}\\
              \coauthor{Hivert~\mycite{Hivert_Thiery.HeckeSg.2006,Hivert_Thiery.HeckeGroup.2007}
                Schilling~\mycite{Hivert_Schilling_Thiery.HeckeGroupAffine.2007,Hivert_Schilling_Thiery.HeckeGroupAffine.2008} 
                Borie}}};
          \node[topic]    (qgraph)   at ( 13,14)        {\mynode{
              \fr{Invariants quasisymétriques de graphes} 
              \en{Quasisymmetric Invariants of Graphs}\\
              \coauthor{Novelli Thibon~\mycite{Novelli_Thibon_Thiery.2004}}}};
          \node[topic]    (steen)    at ( 15,12)        {\mynode{
              \fr{Algèbre de Steenrod rationnelle}
              \en{Rational Steenrod Algebra}\\
              \coauthor{Hivert~\mycite{Hivert_Thiery.SA.2002}}}};
          \node[far area] (top alg)  at (12,10)        {\fr{Topologie algébrique}\en{Algebraic Topology}};
          \node[far area] (factors)  at (15, 6)        {\fr{Inclusions de facteurs}\en{Inclusions of Factors}};
          \node[topic]    (kac)      at (13, 8)        {\mynode{
              \fr{Algèbres de Kac}
              \en{Kac Algebras}\\
              \coauthor{David~\mycite{David_Thiery.2008.Kac}}}};
        \end{pgfonlayer}
        \draw(schubert) to (sym);

        \draw(iso rec)  to (graph);
        \draw(iso rec)  to (inv);

        \draw(age)      to (graph);
        \draw(age)      to (inv);

        \draw (inv perm)to (iso rec);
        \draw (inv perm)to (inv);
        \draw (inv perm)to (sym);
        \draw (schubert)to (inv perm);
        \draw (age)     to (inv perm);

        \draw (sym eval)to (sym);
        \draw (sym eval)to (schubert);
        \draw (sym eval)to (inv perm);
        \draw (sym eval)to (crypto);
        \draw (sym eval)to (slp);
        \draw (slp)     to (inv perm);

        \draw(sym)      to (RT);
        \draw(sym)      to (hopf);
        \draw(RT)       to (hopf);
        \draw(RT)       to (schubert);

        \draw(fd alg)   to (RT);
        \draw(fd alg)   to (hopf);

        \draw(roots)    to (fd alg);
        \draw(roots)    to (RT);
        \draw(roots)    to (schubert);

        \draw(lie)      to (RT);
        \draw(lie)      to (roots);

        \draw(crystals) to (lie);
        \draw(crystals) to (roots);

        \draw(hecke)    to (fd alg);
        \draw(hecke)    to (roots);
        \draw(hecke)    to (crystals);

        \draw(kac)      to (fd alg);
        \draw(kac)      to (hopf);

        \draw(kac)      to (factors);

        \draw(qgraph)   to (hopf);
        \draw(qgraph)   to (sym);
        \draw(qgraph)   to (graph);
        \draw(qgraph)   to (age);

        \draw (steen)   to (sym);
        \draw (steen)   to (hopf);
        \draw (top alg) to (hopf);
        \draw (top alg) to (inv);
        \draw (steen)   to (top alg);
      \end{tikzpicture}

    }
\end{bigcenter}
\caption{Mes sujets de recherche et thèmes avoisinants}
\label{figure.topics}
\end{figure}

La figure~\ref{figure.topics} résume les sujets de recherche que j'ai
abordés et leurs interconnexions. Deux thèmes principaux se
dégagent. Au cœur du premier on trouve les problèmes d'isomorphisme en
combinatoire et leur algébrisation. Dans le second, les modèles
combinatoires deviennent un outil pour étudier des représentations
d'algèbre. Un troisième thème essentiel de mon travail, sous-tendant
les deux autres, est le développement d'outils pour l'exploration
informatique, en particulier dans le cadre du projet logiciel
international \starcombinat que j'ai fondé en 2000.

\section{Algèbres commutatives et isomorphisme en combinatoire}

Le premier volet de mes recherches, présenté au
chapitre~\ref{chapter.algebresCombinatoires}, concerne les problèmes
d'isomorphisme en combinatoire. Ces problèmes sont notoirement
difficiles, l'isomorphisme de graphe étant en quelque sorte l'instance
phare. Ainsi, la fameuse conjecture de reconstruction de graphes de
Ulam n'est toujours pas résolue malgré un demi-siècle de recherches
intenses. Le fil directeur de ce volet est: est-ce que l'encodage
algébrique des problèmes d'isomorphisme peut aider à mieux les
comprendre?

\subsection{Invariants algébriques de graphes et reconstruction}

Le point de départ est mon travail de thèse sous la direction de
Maurice Pouzet. L'objet central en était une algèbre de polynômes
invariants pour une certaine action par permutation du groupe
symétrique qui encode l'isomorphisme de graphe. Il s'agissait
d'évaluer ce que l'étude de cette algèbre, à l'aide de la théorie des
invariants et d'une utilisation intensive du calcul
formel~\cite{Thiery.AIG.2000}, pouvait apporter à la conjecture de
reconstruction de
Ulam~\cite{Thiery.IAGR,Pouzet_Thiery.IAGR.2001,Pouzet_Thiery.IAGR}. Cette
problématique et mes résultats sont notamment repris dans~\cite[5.5
Graph Theory]{Kemper_Derksen.CIT.2002}.  En marge de cette étude, j'ai
introduit avec Jean-Christophe Novelli et Jean-Yves Thibon de
multiples variantes quasi-symétriques des invariants de
graphes~\cite{Novelli_Thibon_Thiery.2004} à la structure simple et
riche.

\subsection{Théorie des invariants effective}

L'algèbre des invariants de graphes s'est révélée être un objet
complexe. En caricaturant, les théorèmes et algorithmes de la théorie
des invariants des groupes finis sont trop généraux pour donner des
résultats fins sur cet exemple. Ceci m'a amené à développer des outils
(bibliothèque \permuvar{}~\cite{Thiery.PDemo.2000} pour \mupad) pour
étudier les invariants de groupes de permutations et à m'intéresser
par la suite aux aspects effectifs et aux applications de la théorie
des invariants. J'ai par exemple mis au point un nouvel algorithme de
calcul de systèmes générateurs de ces
invariants~\cite{Thiery.CMGS.2001}, basé sur des techniques
d'élimination respectant les symétries (bases SAGBI-Gröbner). D'un
autre côté, j'ai obtenu avec Stéphan Thomassé un résultat structurel
sur le comportement de ces invariants vis-à-vis de
l'élimination~\cite{Thiery_Thomasse.SAGBI.2002}.

L'algorithmique buttant sur les limites intrinsèques des techniques
d'élimination, je me suis intéressé aux approches par
évaluation. Après des premiers résultats, dans le cas des fonctions
symétriques et au moyen du modèle SLP (Straight Line
Program)~\cite{Gaudry_Schost_Thiery.2004}, je viens de charger Nicolas
Borie, qui entame une thèse sous ma direction, de l'étude d'une
nouvelle approche des calculs d'invariants de groupes de permutations
par transformée de Fourier.

\subsection{Algèbres d'âge}

En parallèle, j'ai élargi mes recherches, de nouveau avec Maurice
Pouzet, aux algèbres d'âges des structures relationnelles. Ici, les
objets combinatoires, introduits par Rolland Fraïssé, sont les
restrictions finies d'une structure relationnelle infinie $R$ (par
exemple les sous-graphes finis d'un graphe infini), considérés à
isomorphie près. La collection de ces objets est appelée
\emph{âge}. La fonction qui les compte par taille est le \emph{profil}
$\phi_R(n)$. En dépit de la simplicité et de la grande généralité du
cadre, le comportement du profil semble très contraint.
\begin{conjecture}[Pouzet]
  \label{conjecture.profil.rationel}
  Sous des hypothèses faibles, la série génératrice du profil
  $\profile_R(n)$ est une fraction rationnelle dès lors que la croissance
  de $\profile_R(n)$ est sous-exponentielle.
\end{conjecture}

L'encodage algébrique est donné par l'algèbre d'âge de Peter
Cameron. Cet encodage permet d'exploiter la richesse des âges comme
\emph{modèles combinatoires}. Nous avons montré que l'on peut réaliser,
comme algèbres d'âge, outre les invariants de groupes de permutations,
de nombreuses algèbres combinatoires commutatives au cœur de travaux
récents: en premier plan les \emph{polynômes quasi-symétriques} et de
nombreuses variantes.
%

Notre objectif est d'obtenir des informations sur le profil en
utilisant l'algèbre d'âge (dont il donne la série de Hilbert).  Ainsi,
nous démontrons la conjecture~\ref{conjecture.profil.rationel} sous
certaines conditions incluant tous les exemples précités. Plus
généralement, nous cherchons à établir un dictionnaire entre
propriétés combinatoires de l'âge et propriétés de l'algèbre
(engendrement fini, Cohen-Macaulay, etc.). Pour cela, nous tentons de
généraliser les théorèmes et outils que j'avais utilisés en théorie
des invariants. Les résultats ont été annoncés au fur et à mesure à
FPSAC'05~\cite{Pouzet_Thiery.AgeAlgebra.2005},
CGCS'07\footnote{\href{http://cgcs2007.lri.fr/}{International
    Combinatorics, Geometry and Computer Science Conference}},
et font l'objet de deux publications en fin de
préparation~\cite{Pouzet_Thiery.AgeAlgebra1,Pouzet_Thiery.AgeAlgebra2}.

\section{Combinatoire pour la théorie des représentations}

Ce premier volet de ma recherche relève principalement de la
combinatoire algébrique au sens strict: l'objectif est d'algébriser
des objets et problèmes combinatoires pour mieux les comprendre. Dans
le deuxième volet, présenté au
chapitre~\ref{chapter.algebresCombinatoires}, la tendance
s'inverse. Le leitmotiv est la recherche de modèles combinatoires
simples (mais cependant riches!)  pour décrire des structures
algébriques et leur représentations.  En ce sens, il s'agit plutôt
d'\emph{algèbre combinatoire}.
On fait le pari que beaucoup de problèmes d'algèbre ne sont difficiles
qu'en apparence; la clef est alors de trouver le bon point de vue, le
bon modèle dans lequel la démonstration devient courte et élémentaire.
%
%
L'exploration informatique joue donc un rôle inestimable pour essayer
rapidement de nombreux points de vue. 
\FIXME{liant}
En filigrane apparaissent les algèbres de Hopf, les tours d'algèbres
non commutatives, les groupes quantiques, les graphes de
représentations, les systèmes de racines (affines) et les algèbres de
Hecke associées.

\subsection{Algèbres de Hecke groupes}

Dans cette thématique, mon sujet principal est l'étude d'un nouvel
objet, l'algèbre de Hecke groupe associée à un groupe de Coxeter
(section~\ref{section.heckegroupe}). Pour comprendre son intérêt, il
faut d'abord en situer le contexte.

Un thème récurrent du Phalanstère de combinatoire de Marne-la-Vallée
est l'interprétation des algèbres de Hopf combinatoires comme groupes
de Grothendieck des tours d'algèbres de dimension
finie~\cite{Krob_Thibon.NCSF4.1997,Bergeron_Hivert_Thibon.2004,HNTAriki}.
L'exemple originel, dû à Frobenius, est l'algèbre de Hopf des
fonctions symétriques (cf.~\cite{Macdonald.SF.1995,Zelevinsky.1981}),
les fonctions de Schur étant les caractères des représentations
irréductibles du groupe symétrique.

\begin{figure}[h]
  \centering
  \scalebox{.5}{\input{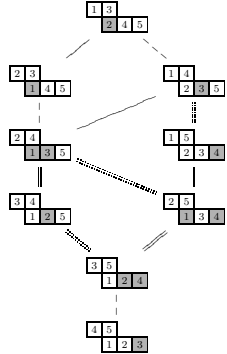}}
  \caption{Un module combinatoire pour l'algèbre de Hecke dégénérée
    $\hecke[5]{0}$}
  \label{figure.Hn0module}
\end{figure}

Dans le cas général de tours d'algèbres non-semi-simples, il faut
distinguer entre représentations simples et projectives. Cela donne
une paire d'algèbres en dualité.  Ainsi, le rôle central joué par la
paire d'algèbres duales Fonctions Symétriques Non Commutatives /
Fonctions Quasi-symétriques vient en particulier du fait qu'elles
encodent la théorie des représentations des $0$-algèbres de
Iwahori-Hecke $H_n(0)$~\cite{Krob_Thibon.NCSF4.1997}, l'algèbre de
Hecke $H_n(q)$ étant une déformation de l'algèbre du groupe symétrique
avec $H_n(1) = \CC[\sg]$. La combinatoire sous-jacente est celle des
rubans et des classes de descentes dans le groupe symétrique (voir
figure~\ref{figure.Hn0module}).

Dans le même temps, l'algèbre de Hecke affine apparaît en filigrane
dans de nombreux travaux du Phalanstère, en particulier comme algèbre
d'opérateurs sur les polynômes. C'est par exemple un outil fondamental
pour l'étude des polynômes de Macdonald. Ses modules irréductibles de
dimension finie ont été classifiés par
Zelevinsky~\cite{Zelevinsky.1980} au moyen de la combinatoire des
multisegments. Dans le cas de la spécialisation centrale principale,
cette combinatoire devient à nouveau celle des classes de
descentes. Depuis quelques années, un objectif du groupe était donc de
résoudre le problème suivant.
\begin{problem}[Jean-Yves Thibon]
  \label{problem.0hecke.affine}
  Expliquer pourquoi les représentations de la $0$-algèbre de Hecke et
  celles de la spécialisation centrale principale de l'algèbre de
  Hecke affine font intervenir la même combinatoire.
\end{problem}
Sa résolution complète a été le fil conducteur notre étude.

À l'occasion d'un groupe de travail en 2003 où nous regardions un
problème de physique faisant intervenir l'algèbre de Hecke affine
(modèle de Frahm-Polychronakos), nous avons construit avec Florent
Hivert, Jean-Christophe Novelli et Jean-Yves Thibon une algèbre
d'opérateurs $\heckesg$ en recollant la $0$-algèbre de Hecke
$\hecke{0}$ et l'algèbre du groupe symétrique $\sg$ \emph{via} leur
représentation régulière à droite. À notre surprise, le calcul sur
ordinateur de petits exemples a révélé une structure riche faisant
intervenir les classes de descentes. Nous nous sommes immédiatement
attelés à en découvrir l'origine. Après de nombreux tâtonnements, j'ai
enfin mis la main sur la bonne description intrinsèque de $\heckesg$
comme algèbre d'opérateurs préservant certaines symétries (ou
certaines anti-symétries). De là, nous avons déroulé le fil avec
Florent Hivert: dimension, base, théorie des représentations, anneau
de Grothendieck des caractères.
Comme effet de bord, cela a donné notre premier exemple d'anneau de
caractères n'étant pas une algèbre de Hopf, la structure d'algèbres et
de cogèbres étant incompatibles; ce fait est maintenant trivial, la
dimension de $\heckesg$ n'étant pas de la forme
$r^nn!$~\cite{Bergeron_Lam_Li.2007}.

Pour cette étude, nous disposions depuis peu d'un outil mis au point
par Florent Hivert et calculant automatiquement la théorie des
représentations des premiers étages d'une tour d'algèbres. Pour
expérimenter, nous avons alors considéré quelques tours d'algèbres
jouets comme l'algèbre du monoïde des fonctions (de parking)
croissantes. À notre grande surprise, celles-ci se sont naturellement
intégrées dans le schéma, nous permettant de définir des
représentations de $(\heckesg)_n$ et $(H_n(q))_n$ sur les puissances
extérieures de la représentation naturelle, et de retrouver comme cas
particulier la tour d'algèbres de Temperley-Lieb.  Ces résultats sont
présentés dans~\cite{Hivert_Thiery.HeckeSg.2006}.  Nous avons depuis
généralisé la construction de $\heckesg$ à n'importe quel groupe de
Coxeter fini. La théorie des représentations reste essentiellement
inchangée: dans tous les cas, elle est Morita-équivalente à celle de
l'algèbre du treillis booléen~\cite{Hivert_Thiery.HeckeGroup.2007}.

Nous avons finalement résolu le problème~\ref{problem.0hecke.affine}
avec Anne Schilling.
\begin{theorem}[Hivert, Schilling, T.~\cite{Hivert_Schilling_Thiery.HeckeGroupAffine.2007,Hivert_Schilling_Thiery.HeckeGroupAffine.2008}]
  Pour tout groupe de Weyl fini, et sauf pour quelques racines de
  l'unité, l'algèbre de Hecke groupe est le quotient naturel de la
  $q$-algèbre de Hecke affine \emph{via} son action de niveau zéro.

  Ce quotient est de plus compatible avec la spécialisation centrale
  principale, et les modules simples associés de l'algèbre de Hecke
  affine donnent, par restriction, les modules projectifs de la
  $0$-algèbre de Hecke.
\end{theorem}



Il reste maintenant plusieurs pistes à explorer ou en cours
d'exploration: le comportement du quotient de l'algèbre de Hecke
affine par son action de niveau zéro lorsque $q$ est une racine de
l'unité, les liens avec les polynômes de Macdonald non symétriques, la
généralisation à tout type du lien avec les fonctions de parking
croissantes, etc. Surtout, il reste à répondre à la question: la
richesse de la structure des algèbres de Hecke groupes est-elle
intrinsèque, ou simplement une ombre portée des algèbres de Hecke
affines? 

\subsection{Opérateurs de promotion sur les graphes cristallins affines}

En marge du sujet que je viens de décrire, j'ai participé à trois
projets de recherche. Le premier, présenté en
section~\ref{section.cristaux}, s'y rattache directement \emph{via} les
outils utilisés (groupes de Weyl affines, actions de niveau zéro,
graphes combinatoires). En effet, ceux-ci jouent un rôle essentiel
dans l'étude des graphes cristallins provenant des représentations de
dimension finie des groupes quantiques affines. Une problématique
importante, faisant l'objet d'une conjecture de Masaki Kashiwara, est
la caractérisation de ces derniers comme produits tensoriels de
graphes cristallins de Kirillov-Reshetikin. Avec Anne Schilling et
Jason Bandlow nous étudions le type $A_n^{(1)}$. La combinatoire
sous-jacente est celle des tableaux. Nous pensons que le cœur du
problème est de montrer que, sur les produits tensoriels de $k$
tableaux, le seul opérateur de promotion est induit par celui défini
sur les tableaux par Schützenberger au moyen du jeu de taquin.  La
démonstration pour $k=2$ fait l'objet d'un article de 31
pages~\cite{Bandlow_Schilling_Thiery.2008.Promotion}.

\subsection{Algèbres de Kac}

L'étude de tours d'algèbres et d'algèbres de Hopf m'a naturellement
amené au deuxième projet de recherche (section~\ref{section.kac}), en
collaboration avec Marie-Claude David, autour des algèbres de Kac de
dimension finie. Cette catégorie d'algèbres de Hopf contient
simultanément les algèbres de groupe et leurs duales, et le point de
vue est proche de celui de la théorie des groupes et de la théorie de
Galois. Les questions centrales sont, par exemple, la détermination du
groupe d'automorphismes et surtout du treillis des sous-structures.
Ce dernier point est principalement motivé par l'existence d'une
correspondance de Galois entre ce treillis et celui des facteurs
intermédiaires de certaines inclusions de facteurs de type $II_1$.
L'étude de deux familles infinies d'exemples fait l'objet d'une
publication de 80 pages~\cite{David_Thiery.2008.Kac}.

\subsection{Polynômes harmoniques pour les opérateurs de Steenrod}

Le dernier projet de recherche que je présente dans ce mémoire
(section~\ref{section.steenrod}) est à l'intersection de mes deux
volets de recherche. Il concerne une conjecture de Reg Wood venant de
topologie algébrique et faisant intervenir l'algèbre de Steenrod. On
peut la formuler comme suit:
\begin{conjecture}[Reg Wood~\cite{Wood.DOSA.1997}, Hivert,
  T.~\cite{Hivert_Thiery.SA.2002}]
    Le sous-espace des polynômes $p$ de $\QQ[x_1,\dots,x_n]$ satisfaisant
  pour tout $k$ l'équation aux dérivées partielles linéaire:
  \begin{displaymath}
    \left(
      \left(1 + x_1\frac{\partial}{\partial x_1}\right) \frac{\partial}{\partial x_1}^k
      + \dots + 
      \left(1 + x_n\frac{\partial}{\partial x_n}\right) \frac{\partial}{\partial x_n}^k
    \right)
    p = 0
  \end{displaymath}
  est isomorphe à la représentation régulière graduée du groupe
  symétrique. En particulier, il est de dimension $n!$\,.
\end{conjecture}
Cette conjecture est un analogue exact d'un résultat très classique
sur les coinvariants du groupe symétrique. Avec Florent Hivert, nous
avons donné une formulation de cette conjecture comme analogue
quantique, en construisant l'algèbre de Steenrod comme déformation non
commutative de l'algèbre de Hopf des fonctions symétriques. Cela nous
a permis d'en déduire des résultats
partiels~\cite{Hivert_Thiery.SA.2002}. Cependant, malgré les efforts
de plusieurs chercheurs, et non des moindres, la conjecture de Reg
Wood résiste toujours.

\section{Exploration informatique et \starcombinatsimple}

Les projets de recherche présentés dans ce mémoire ont en commun
l'exploration, et en particulier l'exploration informatique.  Elle
sert de guide, suggérant des conjectures, occasionnellement donnant
des preuves, ou, au contraire, produisant des contre-exemples.  La
combinatoire algébrique s'y prête bien, car les modèles combinatoires
utilisés donnent des représentations concrètes et effectives des
objets mathématiques à l'étude. 

De plus, on s'intéresse le plus souvent à des familles $(A_n)_{n\in
  \NN}$ d'objets présentant de fortes régularités: typiquement, $A_0$,
$A_1$ sont triviaux, mais les propriétés intéressantes apparaissent
dès $n=3,4,5$ et, si c'est le cas, ont toutes les chances de se
prolonger. En ce sens, nous sommes très loin des expérimentations en
arithmétique où, du fait de la combinatoire des nombres premiers, les
contre-exemples apparaissent souvent très loin. En échange, nous avons
le plus souvent à faire face à une explosion combinatoire: les
exemples triviaux sont les seuls traitables à la main, et $A_5$ sera
par exemple déjà à la limite de ce que les algorithmes classiques
peuvent traiter.

Le défi est de contrôler l'explosion combinatoire, par la
modélisation et l'algorithmique, pour gagner un ou deux crans
supplémentaires. Cela se fait souvent par approximations
successives. La découverte d'un nouveau modèle combinatoire ou d'une
nouvelle propriété permet de mieux comprendre les objets; en retour,
cela permet de calculer plus loin et d'en découvrir de nouvelles
propriétés.

Bien entendu, mener à bien de tels calculs sous-entend un important
travail de programmation, et requiert une large panoplie de techniques
(calcul formel, algèbre linéaire creuse, groupes et représentations,
fonctions symétriques, manipulations de classes combinatoires, séries
génératrices, solveurs divers, etc.). Lors de ma thèse, j'ai regretté
l'absence d'une plate-forme bien établie pour la recherche en
combinatoire algébrique donnant un accès aisé à tous ces outils.

Cela m'a amené à fonder en décembre 2000 le projet \starcombinat, avec
l'aide de Florent Hivert puis, progressivement, de toute une
équipe. Sa mission est de fournir une boîte à outils extensible pour
l'exploration informatique en combinatoire algébrique, avec comme
objectif affiché de fédérer les efforts de développement logiciel
dans la communauté de la combinatoire
algébrique~\cite{Hivert_Thiery.MuPAD-Combinat.2004}.  L'important
investissement initial que m'a demandé ce projet est en train de
porter ses fruits, avec une communauté à l'échelle internationale et
plus d'une quarantaine de publications afférentes (voir
section~\ref{section.combinat.publications} du
chapitre~\ref{chapter.combinat}).

Je présenterai l'apport de l'exploration informatique à chacun de mes
projets de recherche au fil des
chapitres~\ref{chapter.algebresCombinatoires}
et~\ref{chapter.theorieDesRepresentations}. Le
chapitre~\ref{chapter.combinat} est de toute autre nature. J'y
décrirai plus en profondeur le projet \starcombinat. Je détaillerai
notamment les défis particuliers rencontrés lors de son développement,
et les solutions originales que ceux-ci m'ont amené à mettre au point,
tant du point de vue de l'algorithmique que de la conception ou du
choix du modèle de développement.

\bigskip
J'espère montrer, à travers ce mémoire, comment le travail de
recherche et celui d'ingénierie informatique se complètent et se
renforcent mutuellement, le second apportant non seulement des
solutions pratiques au premier, mais aussi une source de questions et
d'inspiration.

\chapter[Algèbres commutatives et isomorphisme en combinatoire]
{Algèbres commutatives graduées\\et problèmes d'isomorphisme en combinatoire}
\label{chapter.algebresCombinatoires}

Le fil directeur de ce chapitre est l'encodage de familles d'objets
combinatoires munies d'une relation d'isomorphisme par des algèbres
commutatives graduées. Dans un premier temps
(section~\ref{section.invariants}), les objets combinatoires sont les
(multi)graphes étiquetés, et l'algèbre est une algèbre de polynômes
invariants pour une certaine action par permutation du groupe
symétrique. Par la suite, le cadre est généralisé à un groupe de
permutation fini quelconque
(section~\ref{section.invariantsPermutation}), puis aux âges des
structures relationnelles (section~\ref{section.age}).

La construction est toujours le même: les objets à un isomorphisme
près forment la base de l'algèbre, la graduation étant donnée par la
taille des objets. Le produit traduit alors toutes les manières de
combiner deux objets pour en construire un plus gros; dans certains
cas, un coproduit traduit réciproquement comment un objet peut se
décomposer en objets plus petits.

L'objectif premier est d'appliquer des outils algébriques à cette
construction pour obtenir des informations sur les problèmes
d'isomorphisme sous-jacent. Mais en retour cette construction donne
des modèles combinatoires riches sur lesquels certaines propriétés
algébriques peuvent être lues. Enfin la question de calcul efficace
dans ces algèbres, et donc sur ces modèles combinatoires, est
centrale, en particulier pour l'exploration informatique.

\section{Invariants algébriques de graphes et reconstruction}
\label{section.invariants}

\subsection{Conjecture de reconstruction de graphes de Ulam}

Le point de départ de ma thèse est la fameuse conjecture de
reconstruction de graphes de Ulam. Elle peut être expliquée en quelques
minutes à un non mathématicien. À cet effet, j'ai eu pendant des années
en permanence dans ma poche le jeu de cartes (transparentes!) présenté
dans la figure~\ref{figure.cartes}.
\newcommand{\graph}[1]{\includegraphics[width=2.5cm]{#1}}
\begin{figure}
  \label{figure.cartes}
  \graph{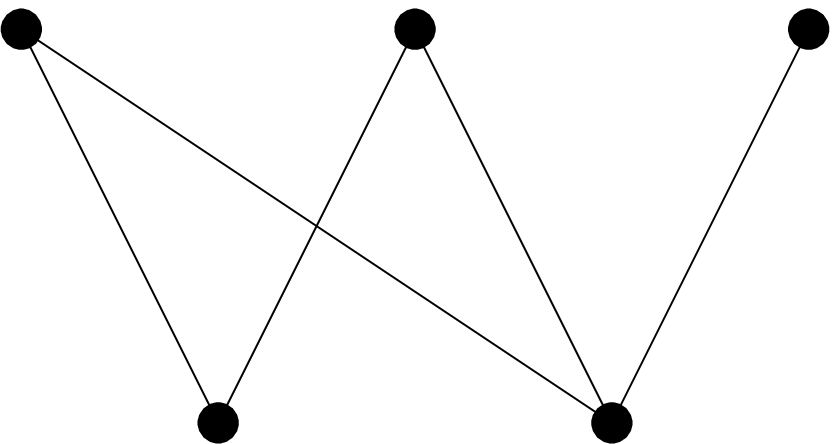}\\
  \bigskip
  \bigskip
  \graph{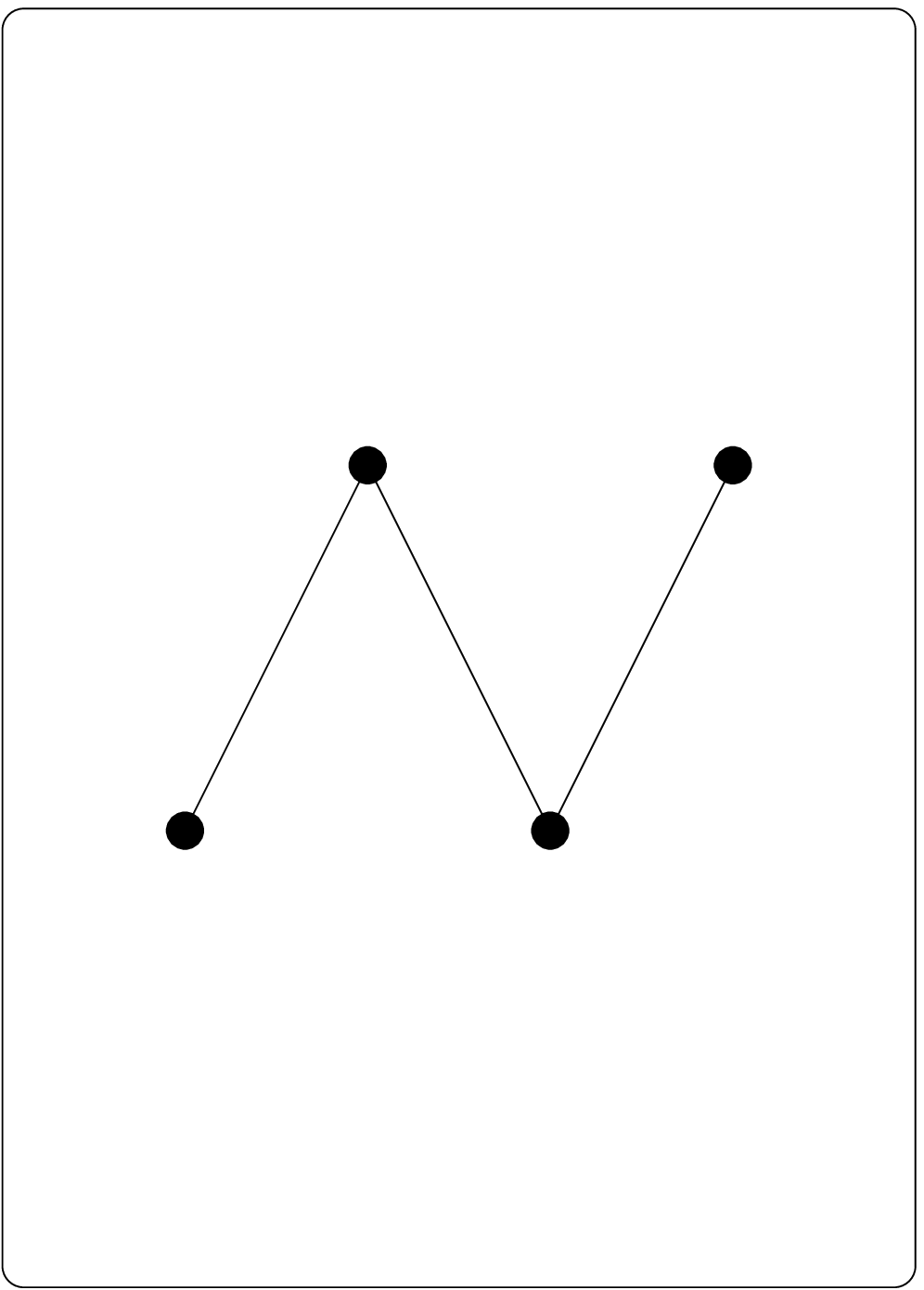}\graph{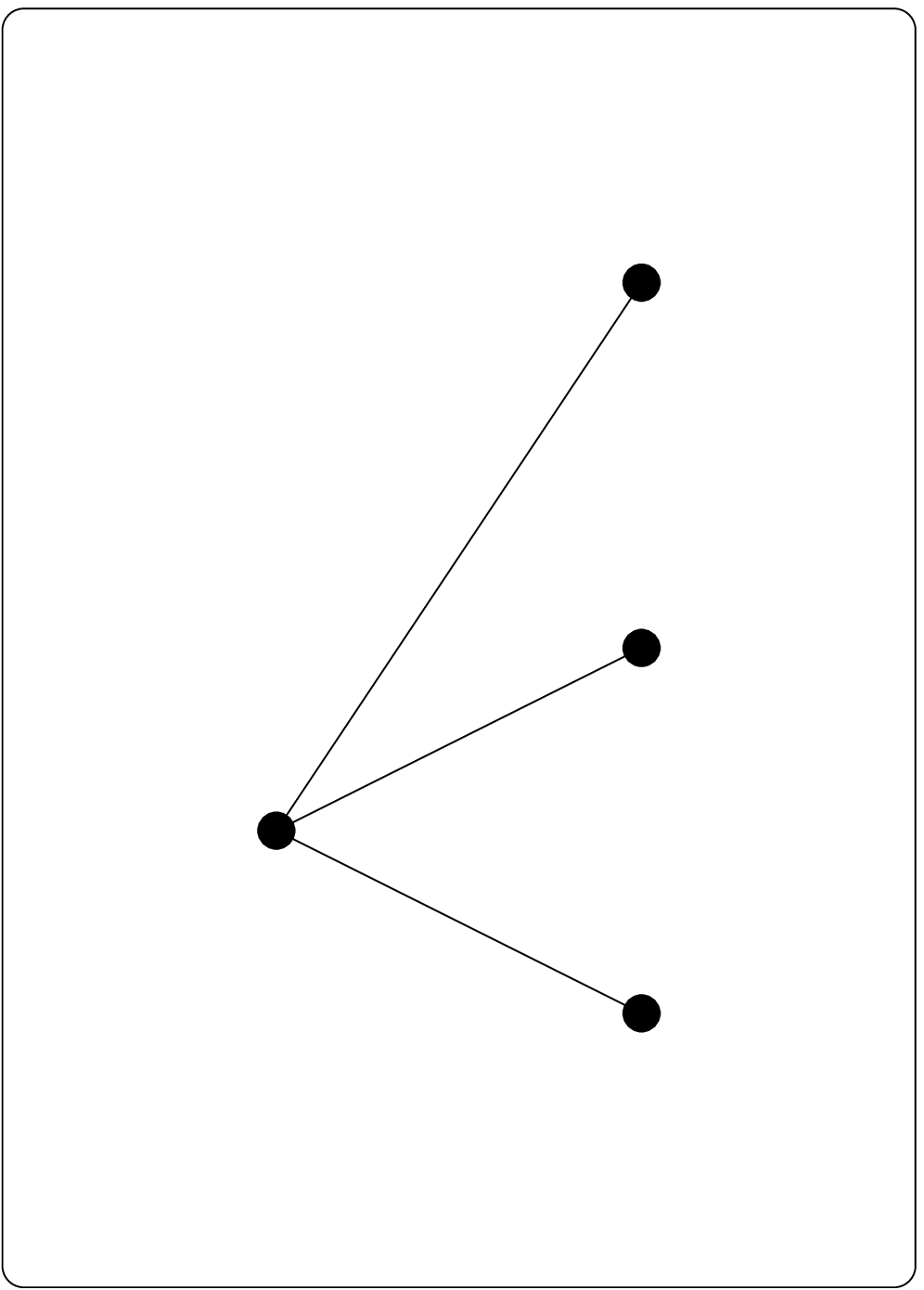}\graph{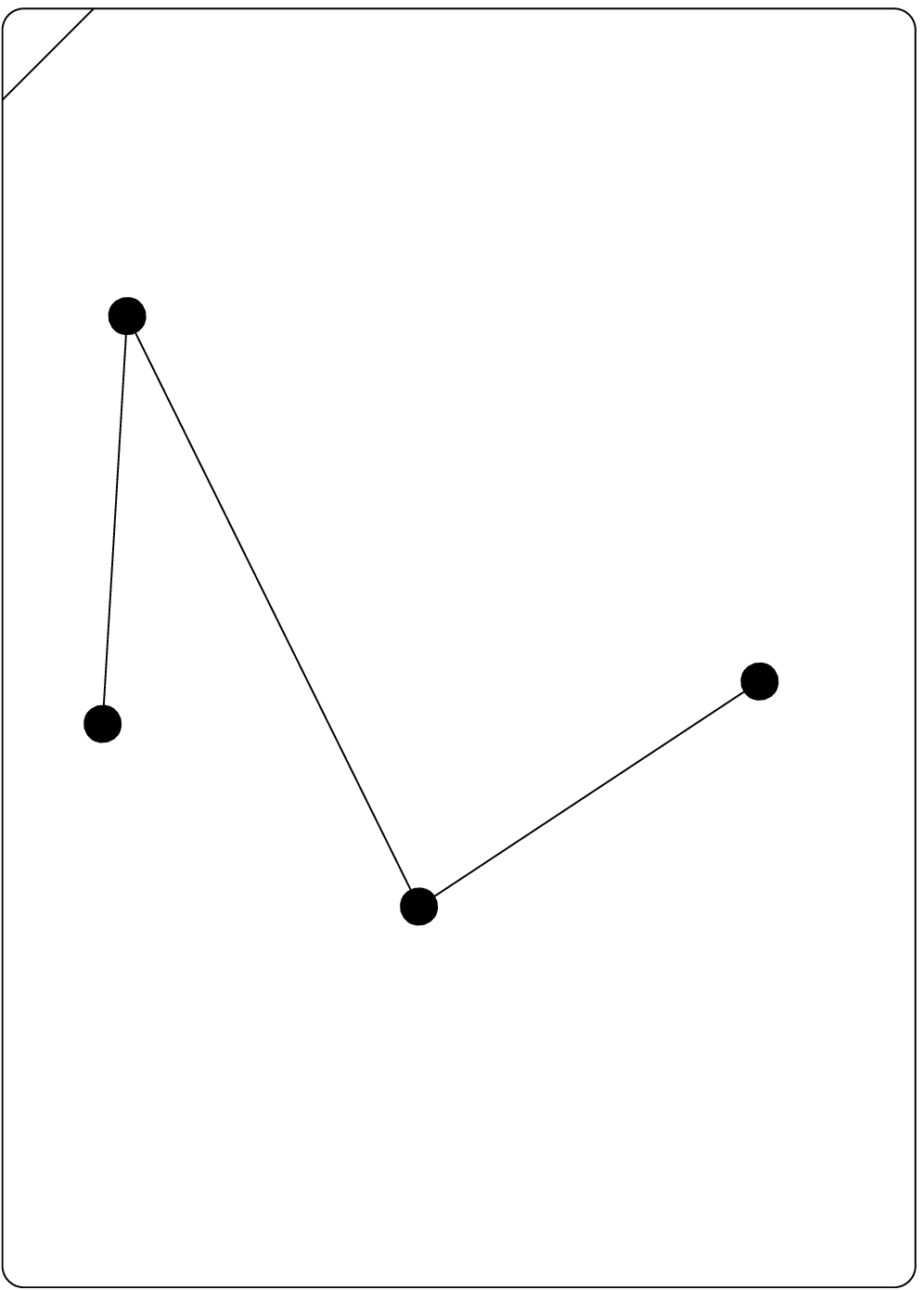}\graph{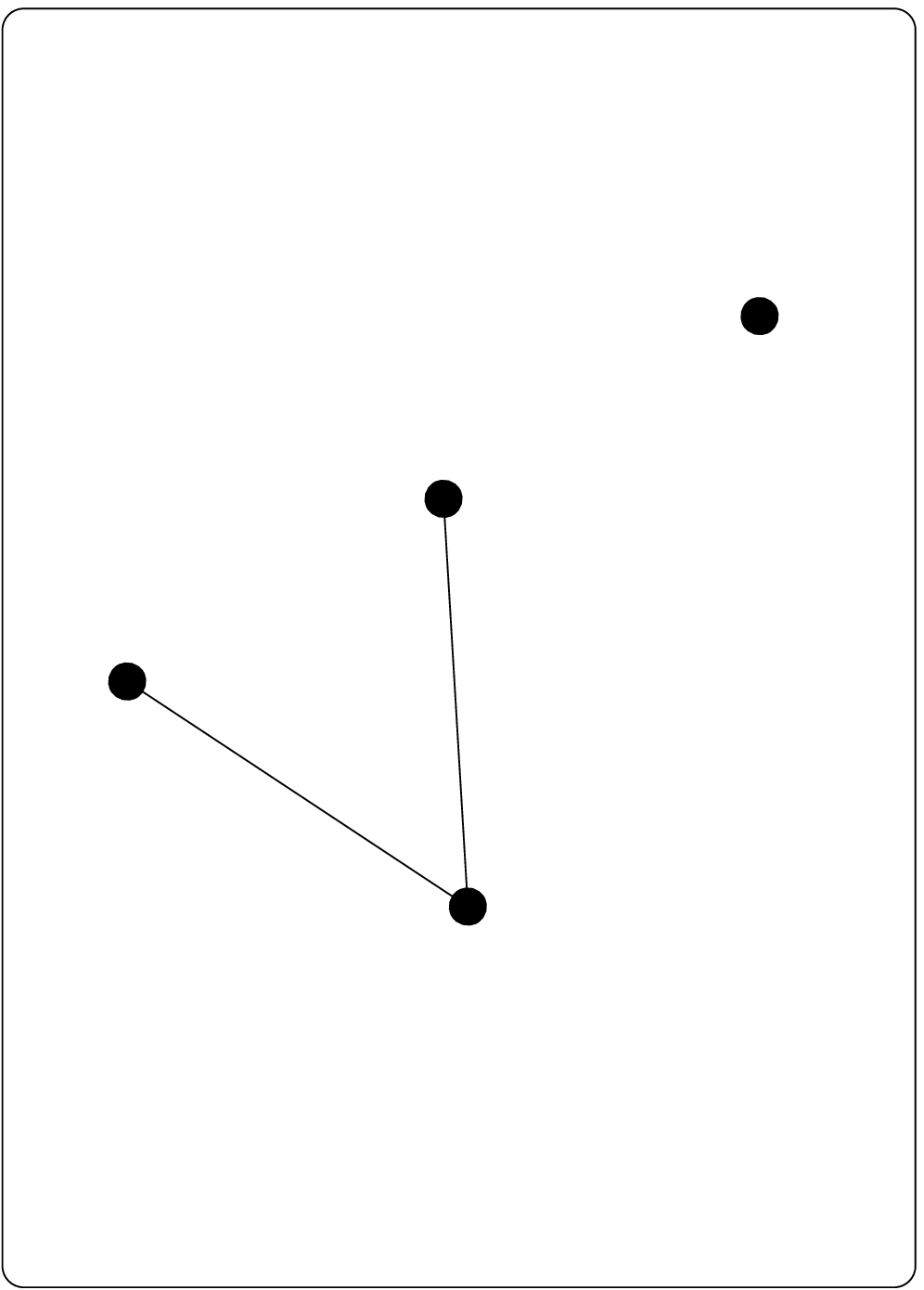}\graph{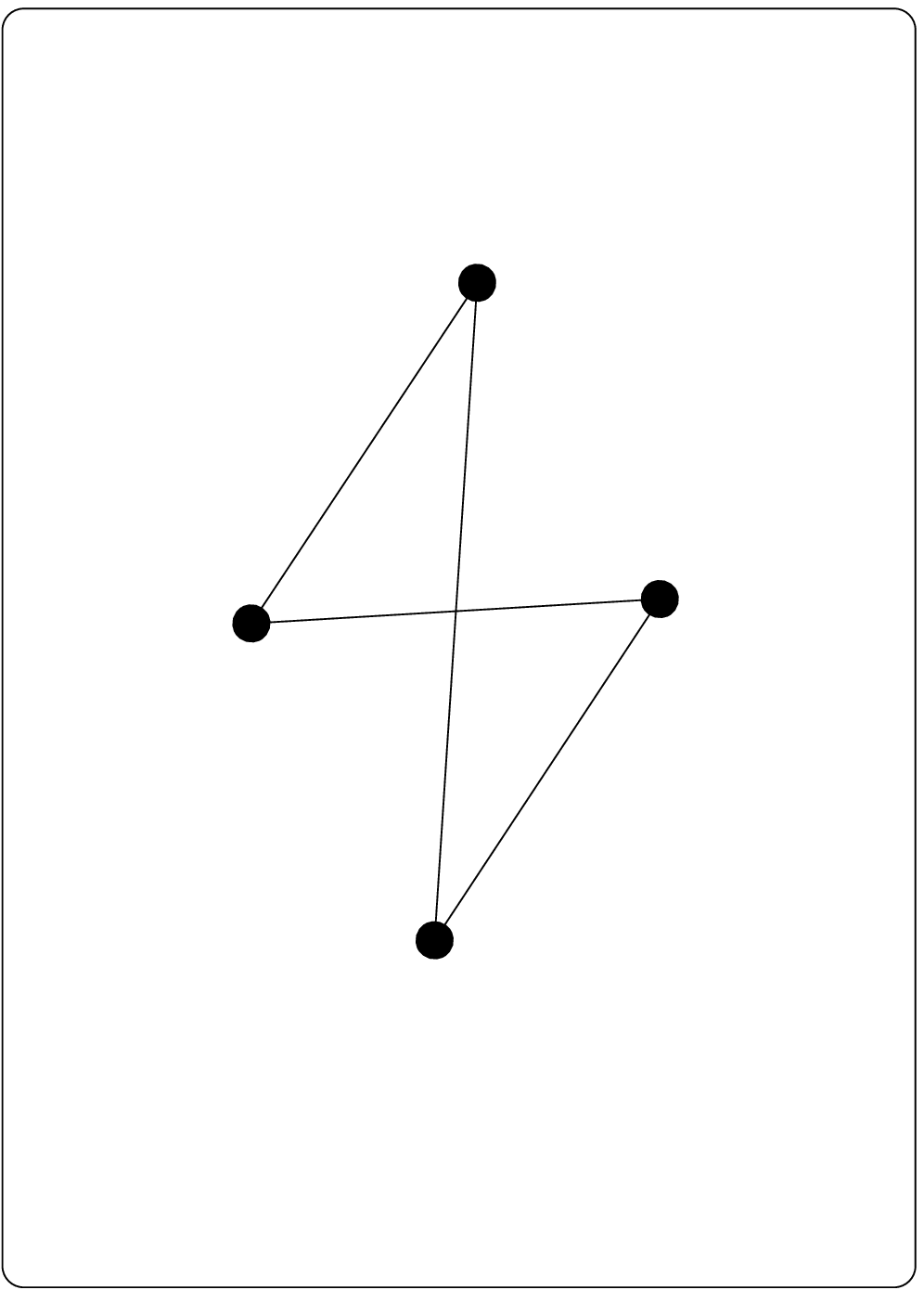}
  \caption{Un graphe simple et le jeu de cartes associé}
\end{figure}
Considérons un graphe simple (pas
de boucles, pas d'arêtes multiples) non étiqueté comme celui à $5$
sommets de la figure. Les cartes en dessous sont obtenues en
supprimant à chaque fois un unique sommet du graphe original. Les
graphes sont considérés à isomorphie près; en particulier, on ne tient
pas compte ici de la disposition géométrique des graphes. De ce fait
les cartes 1 et 3 sont considérées comme identiques. L'ordre des cartes
n'est pas significatif, mais on tient cependant compte de leurs
répétitions.
\begin{question}
  Est-il possible de retrouver le graphe de départ en ne connaissant
  que son jeu de cartes?
\end{question}
C'est un problème de \emph{reconstruction}, similaire à ce que l'on
fait en tomographie (par ex. scanner médical): reconstruire
complètement un objet à partir d'un certain nombre de vues
partielles. Avant même d'attaquer l'aspect algorithmique «comment
reconstruire?» il faut déjà répondre à la question «est-il possible de
reconstruire?». Autrement dit, est-ce que le jeu contient suffisamment
d'informations pour déterminer entièrement le graphe. Si c'est le cas,
le graphe est dit \emph{reconstructible}.
\renewcommand{\graph}[1]{\raisebox{.2ex}{\includegraphics[height=.8ex]{#1}}}
Les deux graphes simples à deux sommets (\graph{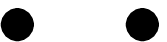} et \graph{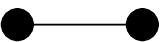}) ont
le même jeu, et ne sont donc pas reconstructibles.  Y en a-t-il
d'autres?
\begin{conjecture}[Ulam 1941~\cite{Ulam.CMP}]
  \label{conjecture.ulam}
  Tous les graphes simples à au moins trois sommets sont reconstructibles.
\end{conjecture}
McKay a vérifié~\cite{McKay.1997} cette conjecture sur ordinateur pour
tous les graphes à au plus 11 sommets (il y en a plus d'un milliard!).
C'est maintenant un vieux problème qui a attiré les meilleurs
chercheurs de la théorie des graphes et qui est à l'origine d'une vaste
littérature~\cite{Bondy.1991} avec une multitude de variantes
(reconstruction par sommets, par arêtes, etc.). Et pourtant il n'est
toujours pas résolu.

Pourquoi ce problème est-il important en général? Même s'il n'a pas
d'application directe, il concentre l'une des difficultés de la
reconstruction (ici l'isomorphie) dans le modèle le plus simple
imaginable. De ce fait, il participe à la taxonomie générale des
problèmes de reconstruction, avec pour objectif d'établir quelles sont
les difficultés intrinsèques du sujet. C'est un guide essentiel pour
le praticien dans le choix d'un bon modèle pour un problème donné de
reconstruction; si ce modèle contient le problème de Ulam, il peut
savoir immédiatement qu'il va au devant de difficultés.

Pourquoi ce problème était-il intéressant pour une thèse? Vu son
historique, il n'était évidemment pas question de l'aborder de
front.
En revanche, il fournit un excellent cas test pour l'utilisation de
nouveaux outils pour traiter de l'isomorphie. Je m'explique. Une des
techniques usuelles pour étudier des objets sous l'action d'un groupe
est de considérer les quantités qui restent invariantes sous cette
action.  Par exemple, le nombre d'arêtes d'un graphe ne change pas
lorsque l'on renumérote ses sommets.  Un des tous premiers résultats
de la théorie de la reconstruction est que le nombre d'arêtes d'un
graphe est un \emph{invariant} \emph{reconstructible}
(i.e. entièrement déterminé par le jeu du graphe). L'idée est très
simple; il suffit de moyenner le nombre d'arêtes sur les sous-graphes
du jeu pour obtenir le nombre total d'arêtes. Plus tard,
Tutte~\cite{Tutte.1979} a démontré que le déterminant (et plus
généralement son polynôme caractéristique) d'un graphe était aussi
reconstructible, et Pouzet avait noté que la preuve consistait à
démontrer que le déterminant s'exprimait par sommes et produits de
paramètres sur les graphes du jeu.

Jusqu'où peut-on espérer généraliser cette technique de preuve? Pour
tenter de répondre à cette question, il est naturel d'introduire
l'algèbre des \emph{invariants polynomiaux de graphes}, c'est-à-dire
l'algèbre $\Inv:=\CC[x_{\{i,j\}}]^{\sg}$ des polynômes en les $\binom n
2$ variables $(x_{\{i,j\}})_{i<j}$ qui sont invariants sous l'action
naturelle du groupe symétrique $\sg$.  Un tel invariant peut être
évalué sur un graphe simple en remplaçant chaque $x_{\{i,j\}}$ par $1$
s'il y a une arête entre $i$ et $j$ dans le graphe, et $0$ sinon. Par
un résultat général de théorie des invariants, les invariants
polynomiaux séparent les graphes à isomorphie près: deux graphes $g$
et $g'$ sont isomorphes si et seulement si ils donnent la même valeur
à tous les invariants polynomiaux; en fait, il suffit d'un nombre fini
de polynômes invariants pour séparer. Dans ce cadre, le nombre
d'arêtes (correspondant à l'invariant $\sum_{i<j} x_{\{i,j\}}$), ou le
déterminant sont des invariants \emph{algébriquement
  reconstructibles}: ils s'expriment par sommes et produit à partir
d'invariants appliqués aux graphes du jeu. Ce qui soulève
naturellement la question suivante.
\begin{question}[Pouzet~\cite{Pouzet.1977,Pouzet.1979}]
  \label{question.inv_alg_rec}
  Pour $n\ge 3$,  les polynômes invariants sont-ils tous
  algébriquement reconstructibles?
\end{question}
Une réponse positive entraînerait une réponse positive à la conjecture
de Ulam. Prudence, donc. Faute de moyens d'investigation, cette
question était restée vierge (à part pour $n=3$ pour lequel on se
ramène facilement aux polynômes symétriques usuels). L'apparition au
début des années 90 d'outils effectifs de calcul dans les
invariants~\cite{Kemper.Invar,Sturmfels.AIT} a motivé mon travail de
thèse: \emph{Étudier, tant d'un point de vue théorique que par
  l'exploration, l'algèbre des invariants sur les graphes (série de
  Hilbert, systèmes générateurs, etc.), et évaluer ce que la théorie
  des invariants peut dire sur les problèmes d'isomorphie et de
  reconstruction de graphe, et en particulier sur la
  question~\ref{question.inv_alg_rec}.}

\subsection{Reconstruction algébrique de graphes}

Commençons par l'aspect reconstruction. Pour une vue synthétique, voir
la figure~\ref{Rec.fig.conj_rec}.
\begin{sidewaysfigure}
  \centering
  \scalebox{.65}{\pgfdeclarelayer{background}
\pgfdeclarelayer{nodes}
\pgfsetlayers{main,nodes}
\begin{tikzpicture}[xscale=.4,yscale=.35]
  \begin{pgfonlayer}{nodes}
    \newcommand{\mynode}[1]{\begin{tabular}{@{}c@{}}#1\end{tabular}}
    \begin{pgfscope}
      \tikzstyle{every node}=[draw,inner sep = 1pt]
      \node(9) at (12,-22) [draw=red] {\mynode{Conj Pouzet (1977):\\Tout multigraphe alg.rec.\\Vérifié $n\leq 5$\\Faux pour $n\geq11$}};
      \node(12) at (12,-30) [draw=red] {\mynode{Pour $d$ assez grand\\Tout multigraphe alg. rec.}};
      \node(15) at (70,-22) [] {\mynode{Pour tout multigraphe $g$\\$\exp(g)$ rec.}};
      \node(18) at (70,-30) [] {\mynode{Pour $d$ assez grand\\Pour tout multigraphe $g$\\$\exp(g)$ rec.}};
      \node(21) at (88,-22) [] {\mynode{Tout graphe valué rec.}};
      \node(24) at (88,-44) [] {\mynode{Conj. Ulam (1941):\\Tout graphe simple rec.\\Vérifié $n\leq 11$ (McKay)}};
      \node(27) at (12,-8) [] {\mynode{Tout poly sur les $0$-réguliers alg. rec.}};
      \node(30) at (12,-14) [] {\mynode{Pour $d$ assez grand\\Tout poly sur les $0$-réguliers alg. rec.}};
      \node(33) at (70,-8) [] {\mynode{Tout poly sur les $0$-réguliers rec.}};
      \node(36) at (70,-14) [] {\mynode{Pour $d$ assez grand\\Tout poly sur les $0$-réguliers rec.}};
      \node(39) at (12,-44) [draw=red] {\mynode{Tout graphe simple alg. rec.\\Vérifié $n\leq6$\\Faux pour $n\geq13$}};
      \node(42) at (12,-52) [] {\mynode{Pour $d=n(n-1)/4$\\Tout graphe simple alg. rec.}};
      \node(45) at (32,-44) [draw=red] {\mynode{Pour le produit disjoint\\Tout graphe simple alg. rec.\\Faux pour $n\geq13$}};
      \node(48) at (12,-60) [] {\mynode{Tout arbre alg. rec.\\Vérifié $n\leq 13$}};
      \node(51) at (32,-52) [draw=red] {\mynode{Pour le produit disjoint\\Pour $d \geq n(n-1)/4$\\Tout graphe simple alg. rec.}};

      \node(54) at (51,-44) [] {\mynode{Pour le produit d'union\\Tout graphe simple alg. rec.}};

      \node(57) at (32,-60) [] {\mynode{Pour le produit disjoint\\Tout arbre alg. rec.}};
      \node(58) at (32,-66) [] {\mynode{Pour le produit disjoint\\Dans l'algèbre des forêts\\Tout arbre alg. rec.}};
      \node(60) at (70,-44) [] {\mynode{Pour tout graphe $g$ simple\\$G\mapsto s(g,G)$ rec.}};
      \node(63) at (51,-60) [] {\mynode{Conj Kocay (1982):\\Pour le produit d'union\\Tout arbre alg. rec.\\Vérifié $n\leq 7$ (Kocay)}};
      \node(64) at (51,-66) [draw=green, fill=green] {\mynode{Pour le produit d'union\\Dans l'algèbre des forêts\\Tout arbre alg. rec.}};

      \node(66) at (70,-60) [] {\mynode{Pour tout arbre $g$\\$G\mapsto s(g,G)$ rec.}};
      \node(69) at (88,-60) [draw=green, fill=green] {\mynode{Tout arbre rec.\\Kelly (1957)}};
      \node(72) at (70,-70) [] {\mynode{Pour tout arbre $g$\\$\exp(g)$ rec.}};
      \node(75) at (70,-38) [] {\mynode{Pour tout graphe $g$ simple\\$\exp(g)$ rec.}};
    \end{pgfscope}
    \tikzstyle{every path}=[thick, double distance=2pt]
    \draw [<->] (9) --  node[auto] {Borne + Div}  (12);
    \draw [<->] (15) --  node[auto] {Borne + Div ?}  (18);
    \draw [->] (12) --  (18);
    \draw [->] (9) --  (15);
    \draw [<->] (15) --  (21);
    \draw [->] (21) --  (24);
    \draw [<->] (27) --  node[auto] {Borne + Div}  (30);
    \draw [<->] (33) --  node[auto] {Borne + Div?}  (36);
    \draw [->] (27) --  (33);
    \draw [->] (9) --  (30);
    \draw [->] (15) --  (36);
    \draw [<->] (39) --  node[auto] {Div + Complément}  (42);
    \draw [->] (12) --  (39);
    \draw [->] (39) --  (45);
    \draw [<->] (45) --  node[auto] {Div}  (51);
    \draw [->] (45) --  (54);
    \draw [<->] (48) --  (57);
    \draw [<->] (60) -- node[auto] {Mnukhin} (24);
    \draw [<->] (54) -- node[auto] {Mnukhin} (60);
    \draw [->] (63) --  (66);
    \draw [->] (60) --  (66);
    \draw [->] (54) --  (63);
    \draw [->] (57) --  (63);
    \draw [->] (42) --  (48);
    \draw [->] (51) --  (57);
    \draw [->] (66) --  (69);
    \draw [->] (24) --  (69);
    \draw [->] (33) -- node[auto] {Sauf exceptions} +(25,0) |- (24);
    \draw [->] (72) --  (66);
    \draw [->] (48) |-  (72);
    \draw [->] (39) -- (18,-38) -- (75);
    \draw [->] (75) --  (60);
    \draw [->] (18) --  (75);
    \draw [<->] (57) --  (58);
    \draw [<->] (64) -|  (69);
    \draw [->] (63) --  (64);
    \draw [->] (58) --  (64);
  \end{pgfonlayer}{nodes}
\end{tikzpicture}}
  \caption[Récapitulatif des conjectures pour les différentes notions
    de reconstructibilité]{Récapitulatif des conjectures pour les différentes notions
    de reconstructibilité par sommets et de leurs relations. Nous ne connaissons
  pas le statut des réciproques non indiquées.  Abréviations:
  rec.=reconstructible; alg. rec.=algébriquement reconstructible;
  \mbox{exp(g)= polynôme invariant associé à g.}}
  \label{Rec.fig.conj_rec}
\end{sidewaysfigure}


J'ai montré que la réponse est positive pour $n\leq 5$ (et très
probablement pour $n=6$), j'ai quelque peu étendu la liste des
invariants classiques de graphes algébriquement reconstructibles, et
donné des propriétés générales sur l'algèbre des invariants
algébriquement reconstructibles. De ces dernières, on déduit que pour
$11\leq n\leq 18$ et très certainement au delà, la réponse à la
question~\ref{question.inv_alg_rec} est négative. Cette approche
est-elle donc vaine?  Sans le dire, nous avons fait ci-dessus un
choix: considérer l'algèbre des invariants polynomiaux sur $\CC$ afin
de pouvoir appliquer les résultats de la théorie des invariants. Cet
objet est plus gros que nécessaire; comme les graphes sont simples,
$x_{\{i,j\}}$ ne prend que les valeurs $0$ et $1$. Il aurait d'abord
été possible de travailler modulo $2$, mais alors la théorie des
invariants devient nettement plus ardue. Une autre option aurait été
de considérer à la place l'\emph{algèbre des graphes simples}, obtenue
en quotientant par $x_{\{i,j\}}^2=x_{\{i,j\}}$. C'est ce qu'avaient
fait avant moi Kocay~\cite{Kocay.1982} et Mnukhin~\cite{Mnukhin.1992}
(voir aussi~\cite{Cameron.1996}); c'est aussi la direction reprise par
la suite par Buchwalder et
Mikkonen~\cite{Mikkonen_Buchwalder.2007}. Au final, leurs résultats
sont pour l'instant de la même teneur que les miens: le problème est
ardu, et une fois obtenue la reconstruction de quelques invariants
explicites, des remarques simples autour des graphes non connexes et
des bornes passablement lâches, on ne peut guère aller au delà. Il y a
cependant deux différences importantes: d'une part, Mnukhin a démontré
que, dans l'algèbre des graphes simples, l'analogue de la
question~\ref{question.inv_alg_rec} \emph{est équivalent à la
  conjecture de Ulam}~\cite{Mnukhin.1992}. En revanche, on perd la
graduation, un outil essentiel en théorie des
invariants.

\begin{figure}[h]
  \begin{displaymath}
    \def\f#1{\multicolumn{1}{|c}{#1}}
    \def\l#1{\multicolumn{1}{c|}{#1}}
    \newcommand\myvcentermiddle[1]{%
      \setbox0\hbox{#1}%
      \dimen255=.5\dp0%
      \advance\dimen255 by -.5\ht0%
      \advance\dimen255 by .5ex%
      \ifvmode\hskip0em\fi\raise\dimen255\box0}
    \def\minigraph#1{\myvcentermiddle{\fbox{\includegraphics[width=.8cm]{#1}}}}
    \begin{array}{ccccccc}
      \minigraph{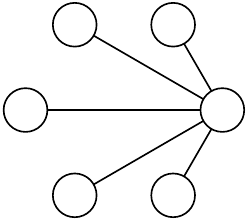} &
      \minigraph{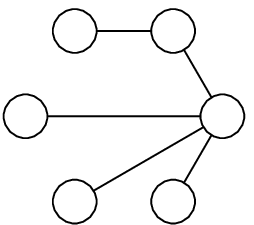} &
      \minigraph{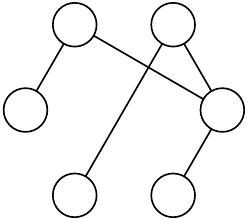} &
      \minigraph{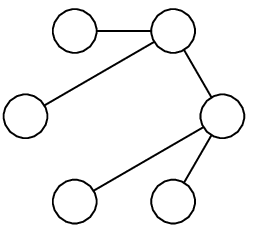} &
      \minigraph{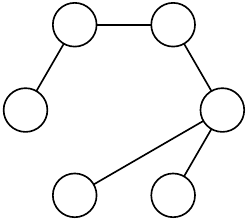} &
      \minigraph{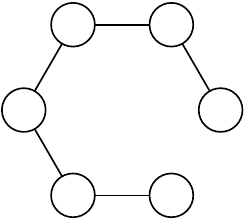} \\
      \f 5 & 1 & . & . & . & \l . &
      \minigraph{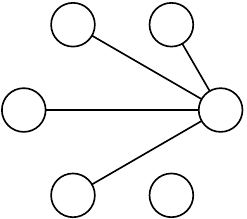} \\
      \f . & 1 & . & . & 1 & \l . &
      \minigraph{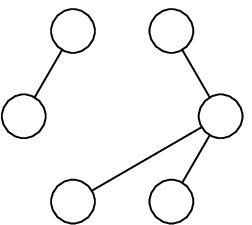} \\
      \f . & 3 & 2 & 4 & 1 & \l . &
      \minigraph{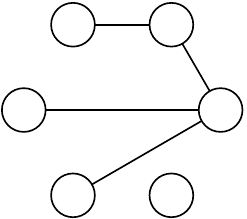} \\
      \f . & . & 2 & . & . & \l 2 &
      \minigraph{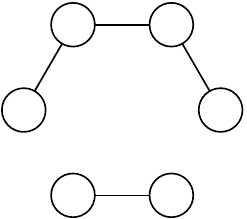} \\
      \f . & . & 1 & . & 2 & \l 2 &
      \minigraph{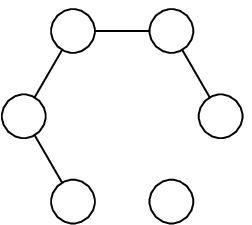} \\
      \f . & . & . & 1 & 1 & \l 1 &
      \minigraph{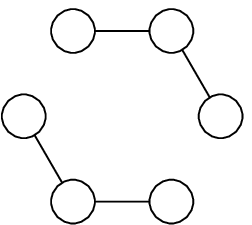} \\
    \end{array}
  \end{displaymath}
  \caption{Matrice d'incidence des arbres versus les forêts à $4$
    arêtes, pour $n=6$ sommets}
  \label{fig.incidence.arbres}
\end{figure}

Il reste un endroit pour lequel je suis convaincu que la
reconstruction algébrique a son mot à dire: une conjecture de Kocay
sur la reconstructibilité du nombre d'arbres couvrants d'un type
donné~\cite[Conjecture 5.1]{Kocay.1982}. On est ici au seuil
du connu, les arbres étant les plus petits graphes
connexes. D'ailleurs, les différentes variantes de l'algèbre
coïncident pour l'essentiel à cet endroit là. J'ai obtenu quelques
résultats partiels dans cette direction, et je ne résiste pas à
mentionner ici ma conjecture préférée issue de cette recherche:
\begin{conjecture}
  Soit $M_n$ la matrice d'incidence $(m_{f,a})_{f,a}$ dont les lignes
  sont indexées par les forêts étiquetées $f$ à $n$ sommets et $n-2$
  arêtes, et les colonnes sont indexées par les arbres $a$ à $n$
  sommets (et donc $n-1$ arêtes) avec $m_{f,a}=1$ si $f$ est un
  sous-graphe de $a$.
  Alors, $M_n$ est de rang maximal, ses lignes étant linéairement
  indépendantes.

  Même conclusion dans le cas non étiqueté, en prenant pour $m_{f,a}$
  le nombre d'occurrences de $f$ dans $a$ (voir
  figure~\ref{fig.incidence.arbres}). Ce deuxième point est un
  corollaire du premier.
\end{conjecture}
J'ai vérifié cette conjecture sur machine jusqu'à $n=19$ dans le cas
non étiqueté. La construction des matrices, très creuses, a été faite
en utilisant \texttt{Nauty} et un script \texttt{Perl}. Le rang a été
calculé par Jean-Guillaume Dumas à l'aide de
\texttt{Linbox}~\cite{Dumas_Villard.Rank.2002}. Pour $n=19$ cela donne
une matrice de dimension $241029\times 317955$, occupant environ
$20$~Mo de mémoire; le calcul a duré cinq jours sur un PC à $1$
GHz. Cette série de matrices a été utilisée comme banc d'essai de
matrices très creuses pour \texttt{Linbox}, confirmant en particulier
une analyse théorique sur des matrices aléatoires de la pertinence de
l'algorithme de Wiedemann pour ce type de
matrices~\cite{Dumas_Villard.Rank.2002,Duran_Saunders_Wan.2003}.

Dans le cas étiqueté, j'ai démontré une borne minimale sur le rang qui
permet de conclure jusqu'à $n=7$. Ironiquement, le calcul sur machine
s'arrête au même endroit (il y a $n^{n-2}$ arbres étiquetés!).

\subsection{Algèbre des invariants de graphes}

Les quelques résultats de reconstruction que j'ai obtenus
justifient-ils cinq ans d'efforts et 300 pages de thèse?  Comme je
l'ai dit, la conjecture de Ulam était surtout un cas test, un angle
d'attaque, pour étudier l'algèbre $\Inv$ des invariants de
graphes. Les résultats de cette étude sont présentés dans ma thèse,
dans~\cite{Thiery.IAGR}, et sont repris
dans~\cite{Kemper_Derksen.CIT.2002}. J'ai donné des propriétés
générales sur l'algèbre $\Inv$, ainsi que des résultats obtenus par
exploration informatique pour $n$ petit. Cela a suggéré plusieurs
conjectures (système de paramètres de petit degré, unimodalité) que
j'étudie en détail. Je fais aussi un tour d'horizon de variantes de
$\Inv$, obtenant par exemple un système générateur très simple pour le
corps des fractions invariantes et infirmant à l'occasion un lemme de
Grigoriev~\cite[Lemma~I]{Grigoriev.1979}.

L'exploration informatique s'est révélée beaucoup plus ardue que
prévu. Illustrons ce point; pour $n\leq3$, $\Inv[n]$ est une algèbre
libre, en fait une algèbre de polynômes symétriques. Le traitement de
$n=4$, à la main avait donné lieu à la
publication~\cite{Aslaksen_al.1996}. Les logiciels existants en
1999~\cite{Kemper.1999} permettaient aussi de traiter $n=4$ par le calcul
en une seconde, mais ne donnaient aucune information pour $n\geq5$. Le
système générateur minimal que j'ai obtenu pour $\Inv[5]$ est
constitué de $57$ polynômes en $10$ variables de degrés jusqu'à
$10$. Traiter $n=6$ complètement est encore hors de portée.

\subsection{Conclusions de cette étude}

$\Inv$ est très loin d'une algèbre libre. À la différence des polynômes
symétriques, elle semble dépourvue de description combinatoire riche
(multiples bases dont de Schur, liens avec la théorie des
représentations, etc.). Elle est aussi beaucoup trop grosse pour
espérer des applications (en particulier algorithmique) à des
problèmes d'isomorphie.

La théorie des invariants donne rapidement des informations
structurelles (engendrement fini, borne sur les degrés, structure de
Cohen-Macaulay, etc.), mais ces informations générales restent
très grossières en pratique. Ainsi, pour $\Inv[5]$, la borne théorique
est de $42$ au lieu de $10$.

Il apparaît clairement que les techniques actuelles de calcul
d'invariants~\cite{Kemper_Derksen.CIT.2002} buttent sur les limites
intrinsèques de l'élimination (bases de Gröbner et variantes), alors
même que les applications en combinatoire requièrent l'étude
d'exemples de taille plus importante.

\section{Théorie des invariants effective}
\label{section.invariantsPermutation}
\subsection{Invariants de groupes de permutations}

Les conclusions de mon étude des invariants de graphes m'ont amené à
développer des outils (bibliothèque \permuvar{}~\cite{Thiery.PDemo.2000} pour \mupad)
pour étudier les invariants de groupes de permutations et à
m'intéresser aux aspects effectifs de la théorie des invariants. 

D'une part, j'ai mis au point un nouvel algorithme de calcul de
systèmes générateurs de ces invariants~\cite{Thiery.CMGS.2001}, basé
sur des techniques d'élimination respectant les symétries (algorithme
de type F4~\cite{Faugere.1999} pour les bases
SAGBI-Gröbner~\cite{Miller.1996}). J'y reviendrai dans les
perspectives.

D'autre part, j'ai obtenu avec Stéphan Thomassé un résultat structurel
sur le comportement de ces invariants vis-à-vis de l'élimination que
je décris maintenant. Les bases
SAGBI~\cite{Kapur_Madlener.1989,Robbiano_Sweedler.1990} sont les
analogues, pour les sous-algèbres des anneaux de polynômes, des bases
de Gröbner pour les idéaux. Comme pour ces dernières, elles s'appuient
sur l'élimination vis-à-vis d'un ordre sur les termes.
Par exemple, la démonstration usuelle du théorème fondamental des
fonctions symétriques se fait par élimination, typiquement vis-à-vis
de l'ordre lexicographique. De fait, les polynômes symétriques
élémentaires $e_1,\dots,e_n$ forment une base SAGBI finie de l'anneau
des polynômes symétriques en $n$ variables. 

Contrairement aux bases de Gröbner, il n'y a pas de théorème de
finitude, et c'est une question ouverte de déterminer pour quelles
sous-algèbres et quels ordres on obtient une base SAGBI finie. Le
théorème suivant indique que, pour un groupe de permutation non
trivial, la base SAGBI est toujours infinie.
\begin{theorem}[T., Thomasse~\cite{Thiery_Thomasse.SAGBI.2002}]
  Soient $G$ un groupe de permutation agissant sur les variables
  $x_1,\dots,x_n$ et $<$ un ordre sur les termes quelconque de
  l'anneau des polynômes $\CC[x_1,\dots,x_n]$. Alors la base SAGBI de
  l'algèbre des invariants $\CC[x_1,\dots,x_n]^G$ vis-à-vis de $<$ est
  finie si et seulement si $G$ est le groupe symétrique (ou un produit
  direct de tels groupes).
\end{theorem}
Nous citons un commentaire du référé anonyme: \emph{"[this result] has
  been desired for a while; thus it can be said that this paper closes
  one of the chapters in the book of invariant theory [...] the proof
  is very beautiful"}.  J'ai depuis généralisé ce résultat au cas des
groupes monomiaux et à de nombreuses algèbres
d'âge~\cite{Pouzet_Thiery.AgeAlgebra.2005,Pouzet_Thiery.AgeAlgebra2}.

\subsection{Complexité d'évaluation des fonctions symétriques}
\label{subsection.sym.eval}

D'un autre côté, suite à une courte collaboration entre Florent Hivert
et moi-même d'une part et Pierrick Gaudry (LIX) et Éric Schost (STIX)
d'autre part, ces derniers ont pu utiliser des algorithmes sur les
fonctions symétriques pour accélérer notablement certains calculs sur
les courbes hyper-elliptiques~\cite{Gaudry_Schost.2004}. Cela nous a
amenés à entreprendre une étude de complexité précise de l'évaluation
de polynômes avec symétries dans le modèle SLP (Straight Line
Programm)~\cite{Gaudry_Schost_Thiery.2004}.  Nous avons ainsi montré
comment, connaissant le coût d'évaluation d'un polynôme symétrique
$P(x_1,\dots,x_n)$ décrit par un SLP, on peut donner le coût
d'évaluation de $P$ en fonction des valeurs des $n$ fonctions
symétriques élémentaires en ces variables. La même technique permet
d'obtenir le coût d'évaluation des coefficients de la décomposition
d'un polynôme $P$ dans une base des polynômes sur les polynômes
symétriques (Schur-Schubert, etc.).

\subsection{Perspectives: calcul d'invariants par transformée de Fourier}
\label{section.invariants.perspectives}

Les techniques d'évaluation que j'ai apprises lors de ce dernier projet
vont resservir pour le calcul d'invariants de groupes de
permutations. Je décris maintenant un projet que j'ai à l'esprit
depuis 2004, et qui fait l'objet de la thèse de mon étudiant Nicolas
Borie.

Le théorème fondamental de la théorie des invariants, démontré par
Hilbert, est l'existence de systèmes finis de générateurs pour
l'algèbre des invariants $\CC[x_1,\dots,x_n]$ d'un sous-groupe fini de
$GL(n)$. La démonstration originale était non constructive. La
théorie des invariants effective s'est développée depuis une quinzaine
d'années~\cite{Kemper.Invar,Kemper_Derksen.CIT.2002,King.2007}, avec pour objectif d'obtenir des
algorithmes (et des implantations) efficaces pour la théorie des
invariants; le problème typique étant le calcul d'un système minimal
de générateurs. L'application première est l'exploration informatique
d'exemples.

La stratégie usuelle utilise la décomposition de Hironaka de l'algèbre
comme module libre sur les invariants primaires
$\Theta_1,\dots,\Theta_n$ pour se ramener à calculer dans le quotient
$\CC[x_1,\dots,x_n]^G / \langle \Theta_1,\dots\Theta_n\rangle$ qui est
de dimension
finie. Pour calculer dans ce quotient, une option est de calculer une
base de Gröbner de l'idéal engendré par $\Theta_1,\dots,\Theta_n$ dans
$\CC[x_1,\dots,x_n]$. C'est par exemple ce qui est implanté dans
\magma. Pour des groupes de matrices avec $n$ petit, cela se
révèle très efficace. En revanche, dès que le nombre de variables
grandit (ce qui est souvent le cas pour des applications en
combinatoire; par exemple l'algèbre $\Inv[5]$ requiert $10$ variables),
ce calcul devient inabordable. Le problème central est que le calcul
de la base de Gröbner casse les symétries, et force à travailler dans
l'algèbre $\CC[x_1,\dots,x_n]$ tout entière, laquelle est de grande dimension dès
que l'on monte en degré.

Dans le cas des groupes de permutations, il est possible d'utiliser
une variante des bases de Gröbner qui préserve les symétries (bases
SAGBI-Gröbner)~\cite{Thiery.CMGS.2001}. Cela permet de calculer un
système générateur de $\Inv[5]$ en quelques minutes; mais
\emph{vérifier} que ce système est effectivement générateur reste
inabordable sans manipulations spécifiques. En effet le calcul complet
nécessite toujours de l'algèbre linéaire dans essentiellement toute
l'algèbre des invariants en degré $22$ (dimension $174403$ à comparer
à $20160075$ pour les polynômes, et $1$ pour le quotient). Il apparaît
ainsi clairement que les techniques actuelles de calcul d'invariants
butent sur les limites intrinsèques de l'élimination.

Aussi paraît-il judicieux d'introduire d'autres points de
vue. L'objectif est d'évaluer une nouvelle stratégie, dans le cas des
groupes de permutations, et de manière plus générale des sous-groupes
de groupes de réflexions.  L'idée est de spécialiser les variables aux
racines de l'unité (transformée de Fourier): cela élimine de facto
deux des obstacles principaux actuels: calculs de produits sur les
monômes (convolution sur le groupe) et calculs dans le quotient par
les invariants primaires (ici, les polynômes symétriques).  Ainsi, le
nombre de points d'évaluations est exactement la dimension du
quotient. Ainsi, pour le problème précédent, on est directement ramené
à des calculs dans une algèbre de dimension $30240$. Qui plus est,
cette algèbre est munie du produit de Hadamard qui est rapide et
préserve les structures creuses.

La première étape est de rédiger une démonstration complète de la
validité de la stratégie, et d'en obtenir une implantation
grossière. Ceci afin de tester concrètement l'approche par des bancs
d'essais comparatifs avec les implantations existantes. Il faudra
aussi comparer avec d'autres approches par
évaluation~\cite{Colin.TIE,Gaudry_Schost_Thiery.2004,Dahan_Schost_Wu.2008}.

Une fois la stratégie validée, le champ d'optimisations est très
ouvert: comment choisir des invariants dont la transformée de Fourier
est creuse, est-il judicieux de représenter les invariants par SLP,
peut-on réduire, par filtration, le nombre de points d'évaluations
lorsque l'on s'intéresse uniquement aux invariants d'un degré donné,
etc. Les progrès viendront principalement de l'étude théorique,
sachant que cette approche fait naturellement apparaître des objets
combinatoires intéressants comme les spécialisations principales des
polynômes de Schur et de Schubert sur un alphabet de la forme
$\frac1{1-q}$ (voir par
exemple~\cite{Littlewood.1950,Bandlow_DAdderio.2008}). Cette
spécialisation principale joue aussi un rôle naturel dans les
descriptions des polynômes de Macdonald et de Schubert par leurs
propriétés
d'évaluation~\cite{Lascoux.2008.MacdonaldAndSchubertForDummies}.

\section{Profil et algèbres d'âge des structures relationnelles}
\label{section.age}

Le travail que je décris dans cette section est un élargissement
naturel de l'utilisation d'outils provenant de la théorie des
invariants pour traiter algébriquement d'autres problèmes
d'isomorphisme en combinatoire. Réalisé en collaboration avec Maurice
Pouzet, il s'inscrit dans la lignée de ses travaux sur le profil des
structures relationnelles (voir~\cite{Pouzet.2006.SurveyProfile}
et~\cite{Pouzet.2008.IntegralDomain} pour des articles de synthèse).

\subsection{Âge et profil d'une structure relationnelle}

Une \emph{structure relationnelle} est une paire
$R:=(E,(\rho_{i})_{i\in I})$, où $E$ est un ensemble (le
\emph{domaine} de $R$) et $\rho_i$ est une famille de relations
$m_i$-aires sur $E$.  Typiquement $R$ est un graphe simple: l'ensemble
de ses sommets est donné par $E$ et l'ensemble de ses arêtes est
décrit par une unique relation $\rho_1$ binaire ($m_1=2$) et
symétrique. Nous prendrons comme exemple la somme directe $3K_\infty$
de trois graphes complets infinis.

Sur tout sous-ensemble $A$ de $E$, $R$ induit par restriction une
sous-structure relationnelle sur $A$. Les notions
d'\emph{isomorphisme}, et de \emph{type d'isomorphie}, sont définies
naturellement.  L'ensemble $\age(R)$ des types d'isomorphie des
restrictions finies de $R$, appelé \emph{âge de $R$}, a été introduit
par Rolland Fraïssé (voir~\cite{Fraisse.TR.2000}). Le \emph{profil} de
$R$ est la fonction $\profile_R$ qui compte pour chaque entier $n$ le
nombre $\profile_R(n)$ de types d'isomorphie de sous-structures de $R$
induites sur les ensembles à $n$ éléments~\cite[Exercise~8
p.~113]{Fraisse.CLM1.1971},~\cite{Pouzet.TR.1978}.

Dans $3K_\infty$, les restrictions sont de nouveau des sommes directes
de trois graphes complets au plus. Un type d'isomorphie de taille $n$
peut donc être décrit par une partition de l'entier $n$ à trois parts
au plus.  Ainsi, la série génératrice de $\profile_{3K_\infty}(n)$ est
donnée par:
\begin{equation}
  \sum_{n\in \NN} \profile_{3K_\infty}(n) Z^n = \frac{1}{(1-Z)(1-Z^2)( 1-Z^3)}\,.
\end{equation}

Si $I$ est fini, $\profile_R(n)$ est nécessairement fini. Afin de
modéliser des exemples venant de l'algèbre ou de la théorie des
groupes, il est cependant nécessaire d'autoriser des ensembles d'index
$I$ infini. Le profil étant fini dans ces exemples, \emph{nous faisons
  toujours l'hypothèse que $E$ est infini et que le profil
  $\profile_R(n)$ est fini}.

\subsection{Algèbre d'un âge}

La construction de l'algèbre d'âge de Peter Cameron suit un paradigme
classique de réalisation d'algèbre sur les mots, à cela près qu'elle
se réalise sur les ensembles~\cite{Cameron.1997}. On considère le
sous-espace des combinaisons linéaires formelles (éventuellement
infinies, mais de degré borné) de sous-ensembles de $R$, que l'on muni
d'un produit commutatif gradué en étendant par linéarité le produit
d'union disjoint sur les ensembles: $AB=A\uplus B$ si $A\cap
B=\emptyset$ et $AB=0$ sinon. L'algèbre d'âge a alors pour base les
«sommes sur orbites» $m_{\overline A} = \sum_{A\in \overline A} A$, où
$\overline A$ parcourt les éléments de l'âge $\age(R)$. Les propriétés
de la restriction garantissent qu'il s'agit bien d'une sous-algèbre;
celle-ci est bien entendue graduée connexe. De plus, par construction, la
série génératrice $\sum_n \profile_R(n) Z^n$ est la série de Hilbert
de $\agealgebra(R)$.

Par cette construction, les structures relationnelles deviennent un
modèle combinatoire pour les algèbres commutatives graduées. Ce modèle
est riche, permettant de réaliser, outre les invariants de groupes de
permutations, de nombreuses algèbres combinatoires commutatives au
cœur de travaux récents: en premier plan les \emph{polynômes
  quasi-symétriques} qui sont au centre de la théorie des algèbres de
Hopf combinatoires~\cite{Aguiar_Bergeron_Sottile.2003} et de
nombreuses variantes, comme par exemple les \emph{polynômes
  quasi-symétriques de graphes} que j'ai introduits et étudiés en 2004
avec J.-C.~Novelli et
J.-Y.~Thibon~\cite{Novelli_Thibon_Thiery.2004}. Mais aussi, par exemple,
l'\emph{algèbre des arbres planaires} de
Gerritzen~\cite{Drensky_Gerritzen.2004,Gerritzen.2004,Gerritzen.2.2004}.

La stratégie est de mettre à jour des liens entre propriétés
combinatoires de la structure relationnelle et propriétés algébriques
de l'algèbre d'âge. Par exemple, Maurice Pouzet a démontré, suite à une
conjecture de Peter Cameron~\cite{Cameron.1997} dans le cadre des
groupes, que l'algèbre d'âge est essentiellement toujours
intègre~\cite{Pouzet.2008.IntegralDomain}.  Il s'ensuit alors, même si
cela n'est pas la démonstration la plus courte, que le profil est une
fonction non décroissante.

\subsection{Décomposition monomorphe finie}

Nous nous intéressons principalement à un nouveau cadre, celui des
structures relationnelles admettant une \emph{décomposition monomorphe
  finie}.  Cette condition est suffisamment large pour couvrir la
plupart des exemples mentionnés ci-dessus (sont exclus l'algèbre des
arbres ainsi que les fonctions quasi-symétriques et leurs variantes sur un
alphabet infini, car leur croissance est exponentielle). Dans ce cas:
\begin{proposition}[Pouzet, T.~\cite{Pouzet_Thiery.AgeAlgebra.2005,Pouzet_Thiery.AgeAlgebra2}]
  Soit $R$ une structure relationnelle admettant une décomposition
  monomorphe finie. Alors, l'algèbre d'âge $\agealgebra(R)$ se plonge
  dans les polynômes en un nombre fini de variables (ou un quotient
  trivial de ceux-ci).
\end{proposition}

Notre premier résultat principal donne des informations très précises
sur le profil.
\begin{theorem}[Pouzet, T.~\cite{Pouzet_Thiery.AgeAlgebra.2005,Pouzet_Thiery.AgeAlgebra1,Pouzet_Thiery.AgeAlgebra2}]
  Soit $R$ une structure relationnelle admettant une
  \emph{décomposition monomorphe finie} avec $k$ composantes infinies.
  Alors, l'âge est un langage rationnel, et le profil de la relation
  est une fraction rationnelle de la forme suivante:
  \begin{equation}
    \frac{P(Z)}{(1-Z)\cdots(1-Z^k)}\,,
  \end{equation}
  avec $P\in\ZZ[Z]$ et $P(1)\ne0$. En particulier $\profile_R(n) \approx n^{k-1}$.
\end{theorem}
Ce théorème confirme en particulier la
conjecture~\ref{conjecture.profil.rationel} dans ce cadre.

Le résultat précédent serait essentiellement trivial si l'algèbre
d'âge était toujours finiment engendrée. C'est loin d'être le cas, et
notre deuxième résultat principal est une caractérisation combinatoire
de ce fait.
\begin{theorem}[Pouzet, T.~\cite{Pouzet_Thiery.AgeAlgebra.2005,Pouzet_Thiery.AgeAlgebra1,Pouzet_Thiery.AgeAlgebra2}]
  Soit $R$ une structure relationnelle admettant une
  \emph{décomposition monomorphe finie}. Alors, l'algèbre d'âge est finiment
  engendrée si et seulement si la décomposition monomorphe est
  \emph{récursivement minimale}. Dans ce cas, l'algèbre est un module
  de type fini sur une sous-algèbre jouant un rôle similaire à celle
  des polynômes symétriques.
\end{theorem}

Plus généralement, notre objectif est d'étudier jusqu'à quel point il
est possible de généraliser chacun des théorèmes et outils que j'ai
utilisés en théorie des invariants (voir table~\ref{table.overview}).

\begin{sidewaystable}[h]
  \Small
  \begin{bigcenter}
    \makeatletter
\newcommand\myvcentermiddle[1]{%
  \setbox0\hbox{#1}%
  \dimen255=.5\dp0%
  \advance\dimen255 by -.5\ht0%
  \advance\dimen255 by .5ex%
  \ifvmode\hskip0em\fi\raise\dimen255\box0}
\makeatother
%
%
\newcommand{\includegraphicsrelation}[1]{
  \myvcentermiddle{\includegraphics[height=10ex]{PICTURES/#1}}}
\newcommand{\tabtitle}[2]{\parbox[m]{#1}{\begin{center}#2\end{center}}}
\newcommand{\polyring}[2]{#1^{\makebox[0cm][l]{$\scriptscriptstyle\
      \leftarrow\in#2$}}}
\renewcommand{\footnoterule}{}
\renewcommand\thefootnote{\alph{footnote}}
\begin{tabular}{|m{26ex}|*{9}{@{\hspace{0.5ex}}c@{\hspace{0.5ex}}|}}
  \hline
  & 
    \tabtitle{12ex}{Structure \mbox{relationnelle}} &
      \tabtitle{15ex}{\begin{center}\mbox{isomorphismes} locaux\end{center}} &
        \tabtitle{8ex}{Série de Hilbert} &
          \tabtitle{14ex}{Engendrement fini} &
            \tabtitle{10ex}{Borne sur le degré} &
              \tabtitle{11ex}{dimension de Krull} &
                \tabtitle{8ex}{Sym-module} &
                  \tabtitle{9ex}{Cohen-Macaulay} &
                    \tabtitle{8ex}{SAGBI finie}
\\\hline
  $|X|<\infty$ &&&
        $\frac{\polyring{P(Z)}{\ZZ[Z]}}{(1-Z)\cdots(1-Z^{|X_\infty|})}$ &
          jamais\footnotemark[1]\footnotetext[1]{\tiny sauf mention
          contraire explicite dessous} &
            $=\infty$\footnotemark[1] &
	      $\leq |X_\infty|$ &
                non\footnotemark[2]\footnotetext[2]{\tiny Ici, "non" signifie
      "pas toujours": il y a des exemples et des contre-exemples} &
		  non\footnotemark[2] &
		    non\footnotemark[2]\\\hline
  Optimalement héréditaire &&&
        $\frac{\polyring{P(Z)}{\ZZ[Z]}}{(1-Z)\cdots(1-Z^{|X_\infty|})}$ &
          oui &
            $<\infty$ &
	      $|X_\infty|$ &
                presque &
		  non\footnotemark[2] &
                    non\footnotemark[2]\\\hline
  Formes préservées &&&
        $\frac{\polyring{P(Z)}{\ZZ[Z]}}{(1-Z)\cdots(1-Z^{|X_\infty|})}$ &
          oui &
            $<\infty$ &
	      $|X_\infty|$ &
                oui &
		  non\footnotemark[2] &
                    non\footnotemark[2]\\\hline
  Polynômes $r$-quasi-symétriques~\cite{Hivert.RQSym.2004} &
    \includegraphicsrelation{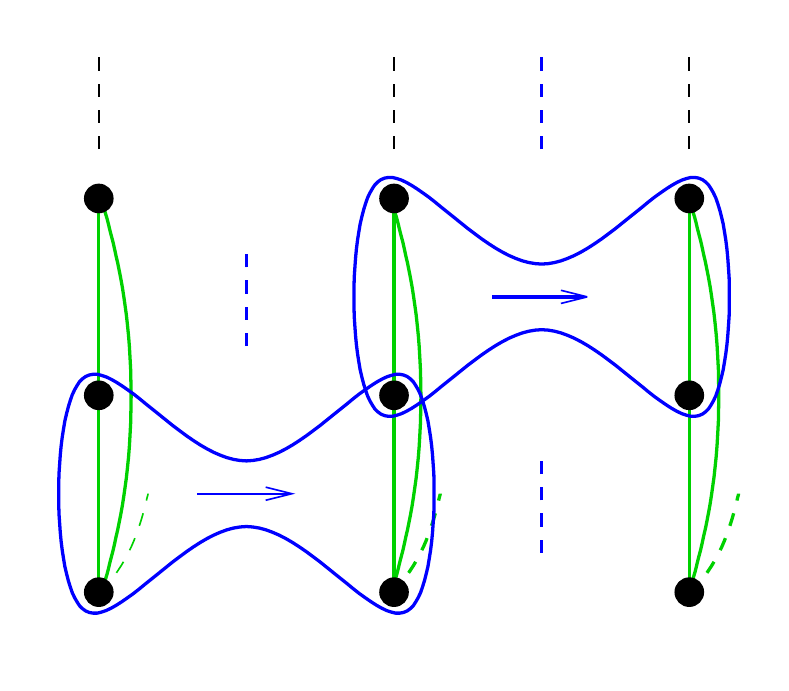} &&
        $\frac{\polyring{P(Z)}{\NN[Z]}}{(1-Z)\cdots(1-Z^{|X|})}$ &
          oui &
            $\leq\frac{|X|(|X|+2r-1)}{2}$ &
	      $|X|$ &
                oui &
		  oui &
                    non\\\hline
  Invariants d'un groupoide\vfil\mbox{de permutations $G$} &&
      $G\wr \sg[\NN]$&
        $\frac{\polyring{P(Z)}{\ZZ[Z]}}{(1-Z)\cdots(1-Z^{|X|})}$ &
          oui &
            $\leq\frac{|X|(|X|+1)}{2}$ &
	      $|X|$ &
                oui &
		  non\footnotemark[2] &
                    \tabtitle{10ex}{jamais\footnotemark[1]}\\\hline
  Exemple non Cohen-Macaulay&
    \includegraphicsrelation{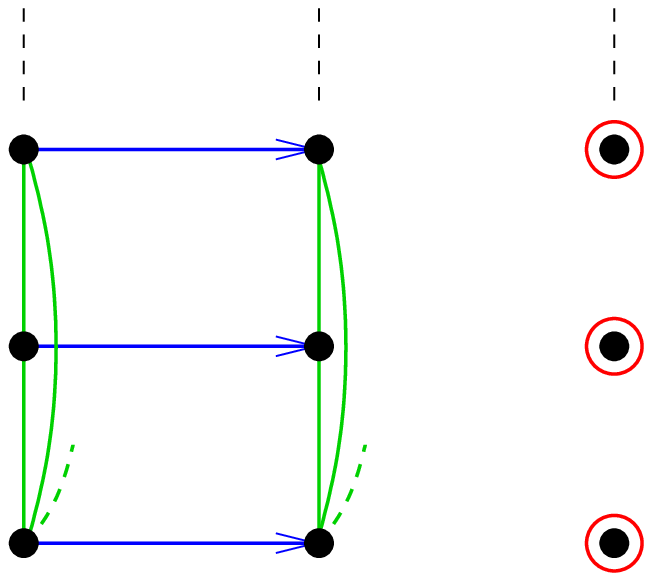} &
      $\langle 1\mapsto 2\rangle \wr \sg[\NN]$ &
        $\frac{1+Z^2+Z^3- Z^4}{(1-Z)^2(1-Z^2)}$ &
          oui &
            $2$ &
	      $|X|$ &
                oui &
		  non &
                    non \\\hline
  Polynômes quasi-symétriques~\cite{Gessel.QSym.1984} &
    \includegraphicsrelation{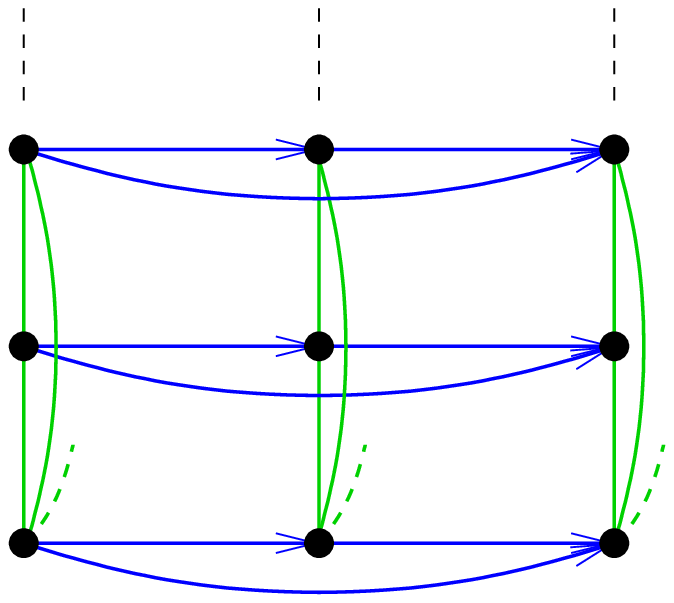} &
      $\operatorname{Inc} \wr \sg[\NN]$ &
        $\frac{\polyring{P(Z)}{\NN[Z]}}{(1-Z)\cdots(1-Z^{|X|})}$ &
          oui&
            $\leq\frac{|X|(|X|+1)}{2}$ &
	      $|X|$ &
                oui &
		  oui~\cite{Garsia_Wallach.2003} &
		    non\\\hline
  Invariants d'un groupe \vfil\mbox{de permutations $G$} &&
      $G\wr \sg[\NN]$&
        $\frac{\polyring{P(Z)}{\in\NN[Z]}}{(1-Z)\cdots(1-Z^{|X|})}$ &
          oui &
            $\leq\frac{|X|(|X|-1)}{2}$ &
	      $|X|$ &
                oui &
		  oui &
		    \tabtitle{10ex}{jamais\footnotemark[1] \cite{Thiery_Thomasse.SAGBI.2002}}\\\hline
  Polynômes symétriques & 
    \includegraphicsrelation{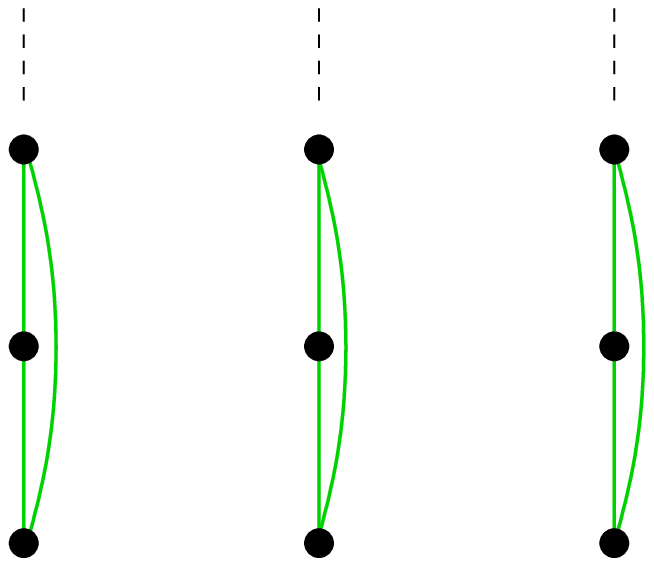} &
      $\sg\wr \sg[\NN]$&
        $\frac{1}{(1-Z)\cdots(1-Z^{|X|})}$ &
          oui &
            $|X|$ &
	      $|X|$ &
                oui &
		  oui &
		    oui\\\hline
  Polynômes &
    \includegraphicsrelation{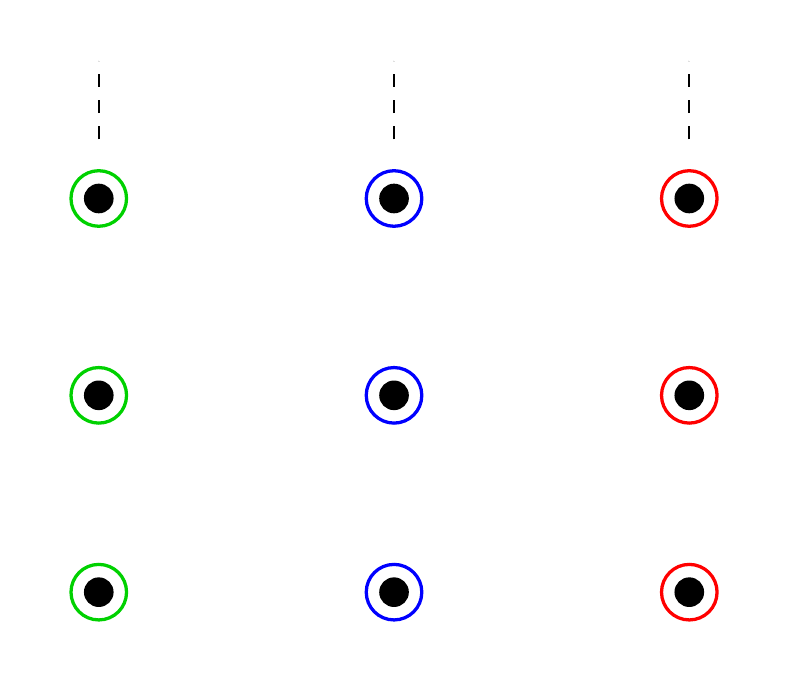} &
      $\id\wr \sg[\NN]$&
        $\frac{(1+Z)\cdots(1+Z+\dots+Z^{|X|})}{(1-Z)\cdots(1-Z^{|X|})}$ &
          oui &
            $1$ &
	      $|X|$ &
                oui &
		  oui &
		    oui\\\hline
\end{tabular}


  \end{bigcenter}
  \caption{Synthèse de nos résultats sur les algèbres d'âge}
  \label{table.overview}
\end{sidewaystable}

Les démonstrations reposent essentiellement sur des généralisations
des techniques d'algèbres de Stanley-Reisner utilisées
dans~\cite{Garsia_Stanton.1984,Thiery.AIG.2000} pour étudier les
invariants de groupes de permutations, sur les ordres d'élimination et
sur la théorie de Ramsey.

Cette recherche est basée sur l'exploration d'une multitude
d'exemples. C'est toute la richesse des structures relationnelles: il
y a une grande souplesse et, selon les contraintes que l'on se fixe,
on peut construire toutes sortes d'exemples exotiques. Le prix à payer
est qu'il n'y a pas de bonne structure de données générique pour
représenter une structure relationnelle; en dehors de cas
particuliers, le calcul sur machine est impuissant.  De ce fait, en
dehors de quelques calculs élémentaires de séries, l'exploration a été
réalisée entièrement au tableau noir.

\subsection{Perspectives}

Une des approches favorites du Phalanstère est de réaliser les
algèbres de Hopf étudiées comme quotients ou sous-algèbres de
l'algèbre des mots non-commutatifs. La combinatoire sous-jacente
devient alors habituellement simple, ce qui permet de donner des
démonstrations élémentaires de la plupart des propriétés
algébriques. Cette approche se complète bien avec l'approche
opéradique de Jean-Louis Loday consistant en particulier à casser les
opérations (produit, coproduits) en plusieurs sous-opérations (algèbre
dendriforme ou tridendriforme), de façon à faire apparaître les
algèbres comme provenant de l'action d'une opérade libre sur un petit
nombre de générateurs. 

Nous essayons avec Jean-Christophe Novelli d'appliquer ces deux
approches à l'algèbre des arbres planaires de Lothar
Gerritzen~\cite{Drensky_Gerritzen.2004,Gerritzen.2004,Gerritzen.2.2004};
cet exemple est intéressant, car il est à la fois combinatoirement
très proche des nôtres (j'ai montré que la base de cette algèbre est
en bijection naturelle avec le quotient des fonctions de parking par
les relations hypoplaxiques), tout en ayant des propriétés algébriques
singulières (le seul coproduit connu n'est pas coassociatif).  La
réalisation que j'ai déjà obtenue en terme d'algèbre d'âge (autrement
dit sur les ensembles) semble un bon premier pas vers cette
réalisation sur les mots.

Cela mène naturellement à des questions sur les algèbres d'âge, et en
tout premier: quelles conditions doit-on imposer sur la structure
relationnelle pour pouvoir définir naturellement un coproduit
coassociatif sur l'algèbre d'un âge, et en faire ainsi une algèbre de
Hopf?

Un autre problème ouvert important, et difficile, est de caractériser
sous quelles conditions ces algèbres sont de Cohen-Macaulay. En effet,
s'il est connu depuis longtemps que les invariants de groupes de
permutations sont de Cohen-Macaulay en toute caractéristique, la
démonstration pour les polynômes quasi-symétriques, même sur les
rationnels, est récente~\cite{Garsia_Wallach.2003}. Existe-t-il une
explication unifiée à ces deux phénomènes au niveau de la combinatoire
des structures relationnelles?

\chapter{Combinatoire pour la théorie des représentations}
\label{chapter.theorieDesRepresentations}

\begin{quotation}
  \href{http://www.quotedb.com/quotes/1360}{«Make everything as simple as possible, but not simpler.»}

  \hfill Albert Einstein
\end{quotation}

Ce chapitre présente mes travaux en théorie des représentations. Les
objets en jeu sont les tours d'algèbres de dimension finie, les
groupes de Coxeter ou de Weyl et leurs algèbres de Hecke (affine), les
algèbres de Hopf et groupes quantiques. Le leitmotiv est la recherche
de modèles combinatoires simples (mais cependant riches!)  pour
décrire ces structures algébriques et leur représentations: tableaux,
systèmes de racines, modèles d'alcôves, et graphes de modules
combinatoires et en particulier graphes cristallins.
Le pari est que la plupart des problèmes (mais pas tous!) ne sont
difficiles qu'en apparence; la clef est alors de trouver le bon point
de vue, le bon modèle dans lequel la démonstration devient courte et
élémentaire.
%
%
Ces modèles permettent aussi une description constructive, voire
effective, des structures algébriques, permettant leur exploration
informatique. En retour, l'exploration informatique joue un rôle
inestimable pour essayer rapidement de nombreux points de vue,
jusqu'à trouver le bon.

\section{Algèbres de Hecke groupes}
\label{section.heckegroupe}

Cette section présente le sujet central de mes recherches en théorie
des représentations: l'algèbre de Hecke groupe d'un groupe de Coxeter,
obtenue par recollement de son algèbre de groupe et sa $0$-algèbre de
Hecke. Comme je l'ai mentionné dans l'introduction, l'intérêt que j'ai
porté à cette algèbre avec mes collaborateurs avait une double
motivation: comprendre les liens entre les représentations de la
$0$-algèbre de Hecke et de l'algèbre de Hecke affine suggérés par la
combinatoire sous-jacente commune des classes de descentes
(problème~\ref{problem.0hecke.affine}), et construire de nouveaux
exemples d'algèbres (si possible de Hopf) comme groupes de
Grothendieck des caractères de tours d'algèbres. À cela se rajoutait
la curiosité devant l'existence d'une structure très riche là où nous
ne l'attendions pas forcément. Enfin, c'était l'occasion pour moi de
comprendre en profondeur des outils (systèmes de racines, algèbres de
Hecke affines) qui me seront utiles pour d'autres projets.

\subsection{Algèbre de Hecke groupe d'un groupe de Coxeter}
\label{subsection.coxeter.groups}

Commençons par quelques préliminaires pour définir l'algèbre de Hecke
groupe d'un groupe de Coxeter $\W$. Notons $(s_i)_{i\in S}$ les
réflexions simples de $W$, et $\Wmax$ son élément maximal (lorsque
$\W$ est fini). On réalise la représentation régulière à droite de
$\W$ en faisant agir les opérateurs $s_i$ sur l'espace vectoriel $\kW$
par $w.s_i := ws_i$%
\footnote{La plupart des actions étant ici à droite, le produit $fg$
  de deux opérateurs $f$ et $g$ dénote sa composition de gauche à
  droite, de sorte que $x.fg=(x.f).g$.}.
De la sorte, on voit $\Wa$ comme sous-algèbre de $\End(\kW)$. De
même, la $0$-algèbre de Hecke peut être réalisée comme la sous-algèbre
de $\End(\kW)$ engendrée par les opérateurs $\pi_i$ définis par $w.\pi_i
:= ws_i$ si $i$ n'est pas une descente de $w$, et $w.\pi_i=w$
sinon. Une variante serait de prendre les opérateurs $\opi_i$ qui
suppriment des descentes. Les opérateurs $\pi_i$ satisfont les
relations de tresse et sont des projecteurs. Ils engendrent un monoïde
de taille $|W|$ qui forme une base de la $0$-algèbre de Hecke.

En type $A$, les opérateurs $s_i$ agissent par transposition sur les
positions, tandis que les opérateurs $\pi_i$ et $\opi_i$ peuvent être
interprétés respectivement comme opérateurs élémentaires d'antitri et
de tri à bulle:
\begin{equation}
    \begin{tikzpicture}[join=bevel,baseline=(current bounding box.east)]
    \def\r#1{{\textcolor{DarkRed}{#1}}}
    \def\b#1{{\textcolor{DarkBlue}{#1}}}
    \matrix [row sep=1cm]{
      \node (sa2)   {$1\r{38}5$};&[-1.6ex] \node (sa6)   {$4\b{62}7$};&[2cm] 
      \node (pia2)  {$1\r{38}5$};&[-1.6ex] \node (pia6)  {$4\b{62}7$};&[2cm] 
      \node (opia2) {$1\r{38}5$};&[-1.6ex] \node (opia6) {$4\b{62}7$};\\
      \node (sb2)   {$1\r{83}5$};&         \node (sb6)   {$4\b{26}7$};&
      \node (pib2)  {$1\r{83}5$};&         \node (pib6)  {$4\b{62}7$};&
      \node (opib2) {$1\r{38}5$};&         \node (opib6) {$4\b{26}7$};\\
      \node (sc2)   {$1\r{38}5$};&         \node (sc6)   {$4\b{62}7$};&
      \node (pic2)  {$1\r{83}5$};&         \node (pic6)  {$4\b{62}7$};&
      \node (opic2) {$1\r{38}5$};&         \node (opic6) {$4\b{26}7$};\\
    };
    \draw[->,DarkRed]  (sa2)   -- node[left]  {$s_2$}    (sb2);
    \draw[->,DarkBlue] (sa6)   -- node[right] {$s_6$}    (sb6);
    \draw[->,DarkRed]  (pia2)  -- node[left]  {$\pi_2$}  (pib2);
    \draw[->,DarkBlue] (pia6)  -- node[right] {$\pi_6$}  (pib6);
    \draw[->,DarkRed]  (opia2) -- node[left]  {$\opi_2$} (opib2);
    \draw[->,DarkBlue] (opia6) -- node[right] {$\opi_6$} (opib6);
    \draw[->,DarkRed]  (sb2)   -- node[left]  {$s_2$}    (sc2);
    \draw[->,DarkBlue] (sb6)   -- node[right] {$s_6$}    (sc6);
    \draw[->,DarkRed]  (pib2)  -- node[left]  {$\pi_2$}  (pic2);
    \draw[->,DarkBlue] (pib6)  -- node[right] {$\pi_6$}  (pic6);
    \draw[->,DarkRed]  (opib2) -- node[left]  {$\opi_2$} (opic2);
    \draw[->,DarkBlue] (opib6) -- node[right] {$\opi_6$} (opic6);
  \end{tikzpicture}
\end{equation}

De manière générale, l'algèbre de (Iwahori)-Hecke générique
$\heckeW{q_1,q_2}$, engendrée par des opérateurs $T_i$ satisfaisant
les relations de tresses ainsi que la relation quadratique
$(T_i-q_1)(T_i-q_2)=0$, peut être construite par interpolation par:
$T_i := (q_1+q_2) \pi_i - q_1 s_i$.

De la sorte, on a réalisé simultanément sur le même espace les
représentations régulières à droite de toutes les algèbres de Hecke de
$\W$, ce qui permet de les recoller.
\begin{definition}[Algèbre de Hecke groupe~\cite{Hivert_Thiery.HeckeGroup.2007}]
  L'algèbre de Hecke groupe est la sous-algèbre de $\End(\kW)$
  engendrée par les opérateurs $s_i$ et $\pi_i$, pour $i$ dans $S$.
\end{definition}

Cette définition originelle n'est pas très informative. Elle a
cependant le mérite de permettre quelques calculs sur machine. Ceux-ci
suivis d'une recherche sur l'encyclopédie des suites
d'entiers~\cite{Sloane} suggèrent une combinatoire sous-jacente forte:
la dimension de cette algèbre compterait le nombre de paires
d'éléments de $W$ sans descentes communes (suite A000275 en type $A$:
$1,3,19,211,\ldots$), tandis que la dimension du quotient semi-simple
par le radical serait donnée par la somme des carrés des tailles des
classes de descentes (suite A060350: $1,2,10,88,\ldots$).

Afin de démontrer ces propriétés, nous avons besoin d'une description
plus conceptuelle de cette algèbre. Commençons par le rang $1$. On
montre facilement que la sous-algèbre parabolique
$\CC[\pi_i,\opi_i,s_i]$ de $\heckeWW$ est de dimension $3$ et que les
relations sont données par:
\begin{equation}
  \label{equation.sigmapi}
  \begin{gathered}
    \begin{aligned}
      s_i\pi_i &= \pi_i\ ,   \qquad &  s_i\opi_i &= \opi_i\ , \\
      \opi_i\pi_i   &= \pi_i\ ,   &  \pi_i\opi_i    &= \opi_i\ , \\
      \pi_is_i &= \opi_i\ ,  &  \opi_is_i &= \pi_i\ , \\
    \end{aligned}\\
    \pi_i+\opi_i=1+s_i\,.
  \end{gathered}
\end{equation}
En particulier, on peut prendre comme générateurs n'importe quel choix
de deux opérateurs dans $\{\pi_i, \opi_i, s_i\}$. De plus, on peut
prendre $\{1, s_i, \pi_i\}$ comme base. Enfin, on note que la droite
$\langle 1-s_i \rangle$ de $\kW$ est stabilisée par l'action de
l'algèbre.

Notre premier résultat est une généralisation de ces remarques, qui,
comme prévu, fait intervenir les classes de descentes
\begin{theorem}[H.,T.~\cite{Hivert_Thiery.HeckeGroup.2007}]
  Une base de $\heckeWW$ est donnée par
  \begin{displaymath}
    \{\,w \pi_{w'} \suchthat w,w'\in W \text{ et }\; \Des(w) \cap \Rec(w') = \emptyset\}\,,
  \end{displaymath}
  où $\Rec(w)$ et $\Des(w)$ désignent respectivement l'ensemble des
  descentes à gauche et à droite d'un élément $w$ de $W$.

  $\heckeWW$ est l'algèbre des opérateurs de $\End(\kW)$ préservant
  les antisymétries à gauche. Sa transposée $\heckeWW^*$ est l'algèbre
  des opérateurs de $\End(\kW)$ préservant les symétries à gauche.
\end{theorem}
La forme de la base suggère une démonstration par règle de
redressement. De fait, nous avons une conjecture sur les relations de
cette algèbre, mais nous n'avons pas de preuve de terminaison pour le
système de réécriture associé (mais cela a-t-il un intérêt?).

En fait, ce théorème donne avant tout une définition alternative plus
conceptuelle de l'algèbre $\heckeWW$. C'est la découverte de cette
définition, \emph{via} l'exploration informatique, qui a permis de
progresser. La démonstration du théorème dans sa globalité est en
effet élémentaire, grâce à deux propriétés de triangularité en
dualité: l'une du côté de la base, et l'autre du côté des relations
linéaires imposées par la préservation des antisymétries.

Cela illustre une idée chère à Alain Lascoux: étudier une algèbre
\emph{via} ses représentations concrètes comme algèbre d'opérateurs
plutôt que \emph{via} une présentation par générateurs et relations.

\subsection{Théorie des représentations}

Une fois trouvée la bonne description de l'algèbre de Hecke groupe
$\heckeWW$, les suggestions de l'exploration informatique permettent
de dérouler sa théorie des représentations; cette théorie est très
uniforme et complètement indépendante du type.

Étant donné un sous-ensemble $I$ de $S$, notons
\begin{displaymath}
  P_I := \left\{ v \in \kW \ | \ s_iv=-v, \forall i\in I \right\}
\end{displaymath}
le sous-espace des vecteurs de $\kW$ antisymétriques à gauche pour
tout $i$ dans $I$. Par construction, $P_I$ est un module pour
$\heckeWW$. Sa dimension est donnée par la taille de la classe de
descente ${}_S^I\!W|$. Par inclusion, la famille $(P_I)_{I\subset S}$
forme un treillis de modules pour $\heckeWW$ anti-isomorphe au
treillis booléen ($I\subset J \Longrightarrow P_J \subset P_I$). 

Il est utile pour la suite de choisir une base de $\kW$ compatible par
restriction avec chaque $P_I$. On peut prendre par exemple:
\begin{displaymath}
  \left\{ v_w := \sum_{w' \in
      W_{S\backslash\Rec(w)}} (-1)^{l(w')} w'w \quad |\quad {w \in W}
  \right\}\,,
\end{displaymath}
qui a le bon nombre de vecteurs avec les bonnes antisymétries. Pour
être snob, on pourrait aussi prendre la base de Kazhdan-Lusztig.
\begin{proposition}[H.T.~\cite{Hivert_Thiery.HeckeGroup.2007}]
  \label{proposition.heckeGroupe.digraphe}
  Une base de $\heckeWW$ est donnée par $\{ e_{w,w'} \suchthat \Rec(w)
  \subset \Rec(w') \}$, où les $e_{w,w'}$ dénotent les unités
  matricielles de $\End(\kW)$ vis-à-vis de la base $v_w$.
\end{proposition}
Cette proposition réalise l'algèbre de Hecke groupe comme
\emph{algèbre d'un digraphe}. Rappelons que l'algèbre d'un digraphe
$g$ est l'algèbre dont la base $\{v_{e,f}\}$ est indexée par les
couples $e,f$ de sommets de $g$ tels qu'il existe un chemin de $e$ à
$f$, et dont le produit est donné par $v_{e,f}v_{e',f'}=\delta_{f,e'}
v_{e,f'}$. C'est le quotient naturel de l'algèbre des chemins lorsque
l'on ne conserve que l'information sur les extrémités des chemins.

Le théorème suivant, décrivant la théorie des représentations de
l'algèbre de Hecke groupe $\heckeWW$, est essentiellement un
corollaire de la proposition~\ref{proposition.heckeGroupe.digraphe}.
\begin{theorem}[Hivert, T.~\cite{Hivert_Thiery.HeckeGroup.2007}]\ \\
  \label{theorem.heckegroup.reptheo}
  \vspace{-1ex}
  \begin{enumerate}[(i)]
  \item La famille $(e_{w,w})$ forme une décomposition maximale de
    l'identité en idempotents orthogonaux;
  \item L'algèbre $\heckeWW$ est Morita équivalente à l'algèbre du
    treillis booléen;
  \item Les $P_I$ sont les modules projectifs indécomposables;
  \item Les modules simples sont obtenus par quotient des modules
    projectifs $S_I := P_I \ /\ \sum_{J\supset I} P_J$;
    leurs éléments sont anti-symétriques à gauche pour $i\in I$ et
    symétriques à gauche pour $i\notin I$. 

    Par restriction, ils donnent exactement:
    \begin{itemize}
    \item Les représentations de Young de forme ruban de $\W$;
    \item Les modules projectifs de $\heckeW{0}$.
    \end{itemize}
  \end{enumerate}
\end{theorem}

En type $A$, on peut calculer explicitement les règles d'induction et
de restriction pour la tour d'algèbre $(\heckesg)_n$. Les structures
d'algèbres et de cogèbres correspondantes sur groupes de Grothendieck
obtenus redonnent des bases connues et des nouvelles bases des
fonctions symétriques non commutatives. Cependant ces structures
d'algèbres et de cogèbres ne sont pas compatibles, de sorte que l'on
obtient pas de nouvelle algèbre de Hopf comme nous l'espérions à
l'origine~\cite{Hivert_Thiery.HeckeGroup.2007}.

\subsection{Algèbres de Hecke groupe et algèbres de Hecke affines}

Le dernier point du théorème~\ref{theorem.heckegroup.reptheo} établit
un lien clair entre les représentations de l'algèbre de Hecke groupe,
et celles de la $0$-algèbre de Hecke. Pour clore le
problème~\ref{problem.0hecke.affine}, il reste à établir un lien entre
l'algèbre de Hecke groupe et l'algèbre de Hecke affine.  C'est l'objet
du résultat suivant: une description alternative, dans le cas des
groupes de Weyl, de l'algèbre de Hecke groupe comme quotient naturel
de l'algèbre de Hecke affine. Il s'ensuit que les modules simples de
l'algèbre de Hecke groupe sont aussi les modules simples de la
spécialisation centrale principale de l'algèbre de Hecke affine.
\begin{theorem}
  \label{theorem.quotient}
  Soit $\W$ un groupe de Weyl affine (éventuellement tordu), et $\clW$
  le groupe de Weyl classique associé.  Soit $\cl: \heckeW{q_1,q_2}\to
  \heckeWW[\clW]$ le morphisme défini par l'action de niveau zéro de
  $\W$ sur $\clW$. Supposons que $q_1,q_2\ne0$ et que $q:=\q$ n'est
  pas une racine $k$-ième de l'unité avec $k\leq
  2\height(\theta^\vee)$). Alors, le morphisme $\cl$ est surjectif et
  fait de l'algèbre de Hecke groupe $\heckeWW[\clW]$ un quotient de
  l'algèbre de Hecke affine $\heckeW{q_1,q_2}$.

  De plus, le morphisme $\cl$ factorise par la spécialisation centrale
  principale de l'algèbre de Hecke affine.
\end{theorem}
La signification de l'action de niveau zéro sera précisée plus loin;
quant à la borne $\height(\theta^\vee)$, nous nous contenterons de
préciser qu'elle est linéaire en le rang de $W$ avec une petite
constante. Par ailleurs, le rôle particulier des racines de l'unité
n'est pas surprenant dans le contexte des algèbres de Hecke.

La démonstration de ce théorème pour $q$ générique repose sur un lemme
combinatoire que nous allons d'abord présenter en type
$A$. L'identification ultérieure de la représentation de niveau zéro
avec une certaine représentation calibrée de série principale permet
de réduire aux petites racines de l'unité les valeurs de $q$
exceptionnelles.

\subsubsection{Type A: transitivité du tri à bulle circulaire}

Nous avons vu que les opérateurs $\pi_i$ agissent par antitri à bulle
élémentaire. En particulier, partant d'une permutation quelconque, par
exemple $51432$, on peut par tri à bulle la transformer en la
permutation maximale $54321$. Cela revient à descendre dans le
permutohèdre. Par contre, l'opération inverse est impossible; on ne
peut pas remonter.

\begingroup
\newcommand{\affinepermutation}[1]{
  \begin{tikzpicture}
    \foreach \theta / \l in { #1 } {
      \node at (90+360/5-\theta*360/5:1) {$\l$};
    }
  \end{tikzpicture}
}
\def\r{\textcolor{DarkRed}}
Écrivons maintenant la permutation $54321$ sur un cercle:
\begin{displaymath}
  \affinepermutation{1/\r5, 2/  4, 3/  3, 4/  2, 5/  1}
\end{displaymath}
Cela introduit naturellement une nouvelle position où l'on peut agir,
entre la dernière lettre et la première. Y appliquant la même règle
que pour les autres positions, nous obtenons:
\begin{displaymath}
  \affinepermutation{1/  1, 2/  4, 3/  3, 4/  2, 5/\r5}
\end{displaymath}
\endgroup
Notons donc $\pi_0$ l'opérateur correspondant. 
Il est clair que l'action de $\pi_0$ tend à faire remonter les
permutations dans le permutohèdre. Peut-on toujours remonter complètement?
\begin{lemma}[H.T. 2005~\cite{Hivert_Schilling_Thiery.HeckeGroupAffine.2008}]
  \label{lemme.tri.A}
  Les opérateurs $\pi_0,\dots,\pi_{n-1}$ agissent transitivement sur
  le groupe symétrique $\sg$.
\end{lemma}
La démonstration de ce lemme repose sur un algorithme de tri à bulle
circulaire récursif\footnote{voir
  \url{http://inst-mat.utalca.cl/fpsac2008/talks/Hivert-Schilling-Thiery.pdf}
  pour une animation expliquant son fonctionnement}.

\subsubsection{Application du lemme combinatoire  en type $A$}

Quel rapport avec l'algèbre de Hecke groupe? Les opérateurs
$\pi_0,\dots,\pi_n$ agissant sur $\sg$ satisfont les relations de la
$0$-algèbre de Hecke affine $\affinehecke{0}$, dont le diagramme de
Dynkin est un cercle. Ils définissent en fait un morphisme $\cl$ de
$\affinehecke{0}$ dans l'algèbre de Hecke groupe. Ce morphisme est-il
surjectif? L'algèbre de Hecke groupe agissant transitivement sur
$\sg$, le lemme~\ref{lemme.tri.A} est une condition nécessaire. Nous
avons montré, \emph{via} la construction d'une base triangulaire
appropriée de $\heckesg$, qu'elle est en fait suffisante.
\begin{theorem}[H.T. 2005~\cite{Hivert_Schilling_Thiery.HeckeGroupAffine.2008}]
  L'action des opérateurs $\pi_0,\dots,\pi_n$ sur $\sg$ définit un
  morphisme surjectif de la $0$-algèbre de Hecke affine
  $\affinehecke{0}$ sur l'algèbre de Hecke groupe $\heckesg$.

  De ce fait, le morphisme de la $q$-algèbre de Hecke affine
  $\affinehecke{q}$ dans l'algèbre de Hecke groupe $\heckesg$ défini
  par interpolation naturelle est surjectif pour $q$ suffisamment
  générique.
\end{theorem}

\subsubsection{Cadre géométrique}

L'action des opérateurs $\pi_0,\dots,\pi_n$ en type $A$ a un pendant
géométrique qui permet de la définir pour tout type. Soit $\W$ un
groupe de Weyl. Il est commode de travailler dans l'espace des copoids
$\coweightspace$.  À chaque racine $\alpha$ est associée une coracine
$\coroot\in\coweightspace$ et un hyperplan $H_\alpha\subset
\coweightspace$ qui coupe $\coweightspace$ en deux demi-espaces
$H_\alpha^+$ et $H_\alpha^-$. La coracine et l'hyperplan définissent
une réflexion $s_\alpha$. Ils définissent aussi une projection
demi-linéaire $\pi_\alpha$ qui fixe $H^-$, et envoie $H^+$ sur $H^-$
par la réflexion $s_\alpha$. Fixons un choix de racines simples
$\alpha_i$. À chaque chambre de l'arrangement d'hyperplans est associé
naturellement un élément de $W$, de telle sorte que l'action des
opérateurs $s_i:=s_{\alpha_i}$ et $\pi_i := \pi_{\alpha_i}$ sur les
chambres est cohérente avec l'action combinatoire de ces mêmes
opérateurs sur $W$.

\begin{figure}[h]
  \begin{tikzpicture}[x={(-1cm,1cm)}, y={(1cm,1cm)}]
  \tikzstyle{point}=[circle,draw,fill=black,inner sep=0mm, minimum size=1mm]
  \tikzstyle{ref}=[inner sep=0mm, minimum size=0mm]
  \tikzstyle{alcove}=[DarkRed]
  \tikzstyle{s0}=[black]
  \tikzstyle{s1}=[DarkBlue]
  
  \node[ref] (tl) at (4,-2)  {};
  \node (tr) at (-2,4)  {};
  \node (bl) at (2,-4)  {};
  \node (br) at (-4,2)  {};

  \node[ref] (L0) at ( 1, 0) {};
  \node[ref] (L1) at ( 0, 1) {};
  \node[ref] (A0) at ( 2,-2) {};
  \node[ref] (A1) at (-2, 2) {};

  \draw[-] (intersection of A0--A1 and bl--tl)[very thin] --
           (intersection of A0--A1 and br--tr) node [below]{$\coweightspace^0$};
  \draw[-,alcove] (intersection of L0--L1 and bl--tl)[very thin] --
           (intersection of L0--L1 and br--tr) node [below]{$\coweightspace^1$};

  \node[point] (L) at (.5,.5)  [label=below:$\rhoc$] {};

  \node[alcove] at (1.15,0.15) {$\phantom{{}^1}0^1$};

  \draw[-] (intersection of 0,0--L0 and bl--br)[very thin,s1] -- 
           (intersection of 0,0--L0 and tl--tr)
           node [above right]{$H_{\alpha_{1,0}}=H_{\alpha_1}$};
  \draw[-] (intersection of 0,0--L1 and bl--br)[very thin,s0] -- 
           (intersection of 0,0--L1 and tl--tr)
           node [above right]{$H_{\alpha_{1,1}}=H_{\alpha_0}$};
  \foreach \x/\c in { -2/s1, -1/s0, 2/s1, 3/s0} {
    \node (start) at (intersection cs: first line  = {(0,0)--(1-\x,\x)},
                     second line = {(tl)--(bl)}) {};
    \node (end)   at (intersection cs: first line  = {(0,0)--(1-\x,\x)},
                     second line = {(tr)--(br)}) {};
    \draw[-,\c] (start) -- (end)[very thin] node [right]{$H_{\alpha_1,\x}$};
  }

  \foreach \x / \l in { -1/$s_0(C)$, 0/$C$, 1/$s_1(C)$ }
    \node at ( .9-1.6*\x,.9+1.6*\x) {\l};

  \foreach \x/\c in { -2/s0, -1/s1, 0/s0, 1/s1, 2/s0, 3/s1} {
    \draw[-,\c] (-.05+\x,0.95-\x) -- (.05+\x,1.05-\x);
  }
  \foreach \x / \l in { -2/$s_0s_1(A)$, -1/$s_0(A)$, 0/$A$,
    1/$s_1(A)$, 2/$s_1s_0(A)$ } {
    \node[alcove] at ( .6-1.15*\x,.6+1.15*\x) {\l};
    \node[point,alcove] at (.5-\x,.5+\x) {};
  }

  \foreach \x / \l in { -1.5, -0.5, 0.5, 1.5 }{
    \node at ( 0.8-1.83*\x,1  +1.77*\x) {+};
    \node at ( 1  -1.77*\x,0.8+1.83*\x) {-};
}

  \draw[->,thick,black] (0,0) -- (L0) node [below left ]{$\Lambdac_0$};
  \draw[->,thick,DarkBlue] (0,0) -- (L1) node [below right]{$\Lambdac_1$};

  \draw[->,thick] (0,0) -- (A0) node [below]{$\coroot_0$}; 
  \draw[->,thick,DarkBlue] (0,0) -- (A1) node [below]{$\coroot_1$}; 
\end{tikzpicture}

  \caption{Réalisation du modèle d'alcôves au niveau $1$ de l'espace
    $\coweightlattice$ des copoids en type $A_1^{(1)}$}
  \label{figure.alcoves}
\end{figure}
Supposons maintenant que $W$ soit un groupe de Weyl affine. De ce
fait, les coracines sont toutes dans un même hyperplan
$\coweightspace^0$ (voir figure~\ref{figure.alcoves} pour le type
$A_1^{(1)}$). Cet hyperplan n'a qu'un nombre fini de chambres, qui
sont en correspondance avec un groupe de Weyl classique $\clW$.
L'action $\cl$ du groupe de Weyl affine $W$ sur $\clW$ est appelée
usuellement \emph{action de niveau zéro} (en référence au niveau
$\ell$ des hyperplans affines $\coweightspace^\ell$ parallèles à
$\coweightspace^0$). L'image $\cl(W)$ de $W$ par l'action est
simplement le groupe classique $\clW$.

Cette même construction géométrique définit aussi une action de niveau
zéro des opérateurs $\pi_0,\dots,\pi_n$ de la $0$-algèbre de Hecke
$\heckeW{0}$ sur le groupe de Weyl classique $\clW$. Mais
contrairement à ce qui se passe pour le groupe, $\cl(\heckeW{0})$ est
plus gros que $\heckeW[\clW]{0}$, car l'opérateur $\pi_0$ ne
s'exprime pas en fonction de $\pi_1,\dots,\pi_n$. \emph{La
dégénérescence de l'algèbre de Hecke affine
\emph{via} l'action de niveau zéro est non triviale}.

On retrouve alors le même lemme combinatoire qu'en type $A$.
\begin{theorem}[S. T. 2008~\cite{Hivert_Schilling_Thiery.HeckeGroupAffine.2008}]
  \label{lemme.tri}
  Les opérateurs $\pi_0,\dots,\pi_n$ agissent transitivement sur le
  groupe de Weyl classique $\clW$.
\end{theorem}
Nous avons donné au cas par cas des algorithmes récursifs de
tri-antitri pour les types classiques dans le même esprit qu'en type
$A$. Nous avons aussi vérifié sur ordinateur, pour tous les types
exceptionnels, l'existence d'un algorithme utilisant le même schéma de
récurrence. Pour $E_7$ et $E_8$, il a fallu utiliser astucieusement la
structure des classes à droite; vérifier directement la forte
connexité du graphe de l'action des opérateurs n'était évidemment pas
souhaitable (696 729 600 sommets). Enfin, nous avons donné une
démonstration géométrique indépendante du type. Les idées
sous-jacentes s'inspirent de notes privées de
Kashiwara~\cite{Kashiwara.2008} sur les représentations de dimensions
finies des groupes quantiques, réinterprétées dans le contexte des
chemins d'alcôves. La figure~\ref{figure.antisorting} illustre cette
démonstration pour tous les groupes de Weyl de rang $2$.

On réobtient comme conséquence de ce théorème un fait connu de
Kashiwara:
\begin{corollary}
  Les graphes cristallins affines finis (tels que ceux étudiés plus
  loin dans la section~\ref{section.cristaux}) sont fortement
  connexes.
\end{corollary}

\begin{figure}
  $\begin{array}{ccc}
    \begin{tikzpicture}[>=latex,join=bevel,scale=.5,baseline=(current bounding box.east)]
\tiny%
  \node (N_1) at (50bp,211bp) [draw=none] {$123$};
  \node (N_2) at (10bp,143bp) [draw=none] {$132$};
  \node (N_3) at (90bp,143bp) [draw=none] {$213$};
  \node (N_4) at (10bp, 75bp) [draw=none] {$312$};
  \node (N_5) at (90bp, 75bp) [draw=none] {$231$};
  \node (N_6) at (50bp,  7bp) [draw=none] {$321$};

  \draw [->,DarkRed]  (N_1) -- node [left]  {$2$} (N_2);
  \draw [->,DarkBlue] (N_1) -- node [right] {$1$} (N_3);

  \draw [->,DarkBlue] (N_2) -- node [left]  {$1$} (N_4);
  \draw [->,DarkRed]  (N_3) -- node [right] {$2$} (N_5);

  \draw [->,DarkRed]  (N_4) -- node [left]  {$2$} (N_6);
  \draw [->,DarkBlue] (N_5) -- node [right] {$1$} (N_6);

  \draw [->]      (N_5) -- node [left =1mm] {$0$} (N_2);
  \draw [->]      (N_4) -- node [right=1mm] {$0$} (N_3);

  \draw [->]      (N_6) -- 
                            (N_1);
\end{tikzpicture} &
    \begin{tikzpicture}[>=latex,join=bevel,scale=.5,baseline=(current bounding box.east)]
\tiny%
  \node (N_4) at (50bp,287bp) [draw,draw=none] {$12$};

  \node (N_6) at (10bp,218bp) [draw,draw=none] {$21$};
  \node (N_2) at (90bp,218bp) [draw,draw=none] {$1\underline{2}$};

  \node (N_8) at (10bp,148bp) [draw,draw=none] {$2\underline{1}$};
  \node (N_5) at (90bp,148bp) [draw,draw=none] {$\underline{2}1$};

  \node (N_3) at (10bp,78bp) [draw,draw=none] {$\underline{1}2$};
  \node (N_7) at (90bp,78bp) [draw,draw=none] {$\underline{2}\underline{1}$};

  \node (N_1) at (50bp,8bp) [draw,draw=none] {$\underline{1}\underline{2}$};

  \draw [->,DarkBlue] (N_4) -- node [left]  {$1$} (N_6);
  \draw [->,DarkRed]  (N_4) -- node [right] {$2$} (N_2);

  \draw [->,DarkRed]  (N_6) -- node [left]  {$2$} (N_8);
  \draw [->,DarkBlue] (N_2) -- node [right] {$1$} (N_5);

  \draw [->,DarkBlue] (N_8) -- node [left]  {$1$} (N_3);
  \draw [->,DarkRed]  (N_5) -- node [right] {$2$} (N_7);

  \draw [->]      (N_5) -- node [above] {$0$} (N_6);

  \draw [->,DarkRed]  (N_3) -- node [left]  {$2$} (N_1);
  \draw [->,DarkBlue] (N_7) -- node [right] {$1$} (N_1);
  \draw [->]      (N_3) -- node [left]  {$0$} (N_4);
  \draw [->]      (N_7) -- node [above] {$0$} (N_8);

  \draw [->]      (N_1) -- node [right] {$0$} (N_2);
\end{tikzpicture}&
    \scalebox{0.75}{\begin{tikzpicture}[>=latex,join=bevel,scale=.5,baseline=(current bounding box.east)]
\tiny%
  \node (N_1) at (50bp,376bp) [draw,draw=none] {$1$};

  \node (N_2) at (10bp,314bp) [draw,draw=none] {$$};
  \node (N_3) at (90bp,314bp) [draw,draw=none] {$$};

  \node (N_4) at (10bp,252bp) [draw,draw=none] {$$};
  \node (N_11) at (90bp,252bp) [draw,draw=none] {$$};

  \node (N_7) at (10bp,190bp) [draw,draw=none] {$$};
  \node (N_5) at (90bp,190bp) [draw,draw=none] {$$};

  \node (N_9) at (10bp,128bp) [draw,draw=none] {$$};
  \node (N_6) at (90bp,128bp) [draw,draw=none] {$$};

  \node (N_10) at (10bp,66bp) [draw,draw=none] {$$};
  \node (N_8) at (90bp,66bp) [draw,draw=none] {$$};

  \node (N_12) at (50bp,4bp) [draw,draw=none] {$\Wmax$};

  \draw [->,DarkRed] (N_1) -- node [left] {$2$} (N_2);
  \draw [->,DarkBlue] (N_1) -- node [right] {$1$} (N_3);

  \draw [->,DarkBlue] (N_2) -- node [left] {$1$} (N_4);
  \draw [->,DarkRed] (N_3) -- node [right] {$2$} (N_11);

  \draw [->,DarkRed] (N_4) -- node [left] {$2$} (N_7);
  \draw [->,DarkBlue] (N_11) -- node [right] {$1$} (N_5);

  \draw [->,DarkBlue] (N_7) -- node [left] {$1$} (N_9);
  \draw [->,DarkRed] (N_5) -- node [right] {$2$} (N_6);
  \draw [->] (N_5) -- node [above=2mm]{$0$} (N_4);

  \draw [->,DarkRed] (N_9) -- node [left] {$2$} (N_10);
  \draw [->,DarkBlue] (N_6) -- node [right] {$1$} (N_8);
  \draw [->] (N_9) -- node [right] {$0$} (N_3);
  \draw [->] (N_6) -- node [above=2mm] {$0$} (N_7);

  \draw [->,DarkBlue] (N_10) -- node [left] {$1$} (N_12);
  \draw [->,DarkRed] (N_8) -- node [right] {$2$} (N_12);
  \draw [->] (N_10) -- node [right]{$0$} (N_11);
  \draw [->] (N_8) -- node [above=3mm] {$\,0$} (N_1);

  \draw [->] (N_12) -- node [below=3mm] {$\!0$} (N_2);

\end{tikzpicture}}\\\\
    %
    \scalebox{.9}{\begin{tikzpicture}[baseline=(current bounding box.east)]
\draw[->,  color = black,] (0.0,0.0) -- (0.0,-1.732050808) node[at end, auto=right] {$\alpha^\vee_{0}$};
\draw[->,  color = DarkBlue,] (0.0,0.0) -- (1.5,0.8660254038) node[at end, auto=right] {$\alpha^\vee_{1}$};
\draw[->,  color = DarkRed,] (0.0,0.0) -- (-1.5,0.8660254038) node[at end, auto=right] {$\alpha^\vee_{2}$};
\draw[ color = black, very thick,](0.5,0.8660254038) -- (-0.5,0.8660254038);
\draw[ color = DarkBlue, ,](0.0,0.0) -- (-0.5,0.8660254038);
\draw[ color = DarkRed, very thin,](0.0,0.0) -- (0.5,0.8660254038);
\draw[ color = black, very thick,](0.5,0.8660254038) -- (1.0,0.0);
\draw[ color = DarkBlue, ,](0.0,0.0) -- (1.0,0.0);
\draw[ color = DarkRed, very thin,](0.0,0.0) -- (0.5,0.8660254038);
\draw[ color = black, very thick,](0.5,-0.8660254038) -- (1.0,0.0);
\draw[ color = DarkBlue, ,](0.0,0.0) -- (1.0,0.0);
\draw[ color = DarkRed, very thin,](0.0,0.0) -- (0.5,-0.8660254038);
\draw[ color = black, very thick,](0.5,-0.8660254038) -- (-0.5,-0.8660254038);
\draw[ color = DarkBlue, ,](0.0,0.0) -- (-0.5,-0.8660254038);
\draw[ color = DarkRed, very thin,](0.0,0.0) -- (0.5,-0.8660254038);
\draw[ color = black, very thick,](-1.0,0.0) -- (-0.5,0.8660254038);
\draw[ color = DarkBlue, ,](0.0,0.0) -- (-0.5,0.8660254038);
\draw[ color = DarkRed, very thin,](0.0,0.0) -- (-1.0,0.0);
\draw[ color = black, very thick,](-1.0,0.0) -- (-0.5,-0.8660254038);
\draw[ color = DarkBlue, ,](0.0,0.0) -- (-0.5,-0.8660254038);
\draw[ color = DarkRed, very thin,](0.0,0.0) -- (-1.0,0.0);
\draw[ color = black, very thick,](0.5,2.598076211) -- (-0.5,2.598076211);
\draw[ color = DarkBlue, ,](0.0,1.732050808) -- (-0.5,2.598076211);
\draw[ color = DarkRed, very thin,](0.0,1.732050808) -- (0.5,2.598076211);
\draw[ color = black, very thick,](0.5,2.598076211) -- (1.0,1.732050808);
\draw[ color = DarkBlue, ,](0.0,1.732050808) -- (1.0,1.732050808);
\draw[ color = DarkRed, very thin,](0.0,1.732050808) -- (0.5,2.598076211);
\draw[ color = black, very thick,](0.5,0.8660254038) -- (1.0,1.732050808);
\draw[ color = DarkBlue, ,](0.0,1.732050808) -- (1.0,1.732050808);
\draw[ color = DarkRed, very thin,](0.0,1.732050808) -- (0.5,0.8660254038);
\draw[ color = black, very thick,](0.5,0.8660254038) -- (-0.5,0.8660254038);
\draw[ color = DarkBlue, ,](0.0,1.732050808) -- (-0.5,0.8660254038);
\draw[ color = DarkRed, very thin,](0.0,1.732050808) -- (0.5,0.8660254038);
\draw[ color = black, very thick,](-1.0,1.732050808) -- (-0.5,2.598076211);
\draw[ color = DarkBlue, ,](0.0,1.732050808) -- (-0.5,2.598076211);
\draw[ color = DarkRed, very thin,](0.0,1.732050808) -- (-1.0,1.732050808);
\draw[ color = black, very thick,](-1.0,1.732050808) -- (-0.5,0.8660254038);
\draw[ color = DarkBlue, ,](0.0,1.732050808) -- (-0.5,0.8660254038);
\draw[ color = DarkRed, very thin,](0.0,1.732050808) -- (-1.0,1.732050808);
\draw[->,  color = purple,thick] (0.0,1.154700538) -- (0.0,0.5773502692);
\end{tikzpicture}}&
    \scalebox{1.2}{\begin{tikzpicture}[baseline=(current bounding box.east)]
\draw[->,  color = black,] (0.0,0.0) -- (0.0,-1.0) node[at end, auto=right] {$\alpha^\vee_{0}$};
\draw[->,  color = DarkBlue,] (0.0,0.0) -- (-1.0,1.0) node[at end, auto=left] {$\alpha^\vee_{1}$};
\draw[->,  color = DarkRed,] (0.0,0.0) -- (1.0,0.0) node[at end, auto=right] {$\alpha^\vee_{2}$};
\draw[ color = black, very thick,](0.0,0.5) -- (0.5,0.5);
\draw[ color = DarkBlue, ,](0.0,0.0) -- (0.5,0.5);
\draw[ color = DarkRed, very thin,](0.0,0.0) -- (0.0,0.5);
\draw[ color = black, very thick,](0.0,0.5) -- (-0.5,0.5);
\draw[ color = DarkBlue, ,](0.0,0.0) -- (-0.5,0.5);
\draw[ color = DarkRed, very thin,](0.0,0.0) -- (0.0,0.5);
\draw[ color = black, very thick,](0.0,-0.5) -- (0.5,-0.5);
\draw[ color = DarkBlue, ,](0.0,0.0) -- (0.5,-0.5);
\draw[ color = DarkRed, very thin,](0.0,0.0) -- (0.0,-0.5);
\draw[ color = black, very thick,](0.0,-0.5) -- (-0.5,-0.5);
\draw[ color = DarkBlue, ,](0.0,0.0) -- (-0.5,-0.5);
\draw[ color = DarkRed, very thin,](0.0,0.0) -- (0.0,-0.5);
\draw[ color = black, very thick,](-0.5,0.0) -- (-0.5,0.5);
\draw[ color = DarkBlue, ,](0.0,0.0) -- (-0.5,0.5);
\draw[ color = DarkRed, very thin,](0.0,0.0) -- (-0.5,0.0);
\draw[ color = black, very thick,](-0.5,0.0) -- (-0.5,-0.5);
\draw[ color = DarkBlue, ,](0.0,0.0) -- (-0.5,-0.5);
\draw[ color = DarkRed, very thin,](0.0,0.0) -- (-0.5,0.0);
\draw[ color = black, very thick,](0.5,0.0) -- (0.5,0.5);
\draw[ color = DarkBlue, ,](0.0,0.0) -- (0.5,0.5);
\draw[ color = DarkRed, very thin,](0.0,0.0) -- (0.5,0.0);
\draw[ color = black, very thick,](0.5,0.0) -- (0.5,-0.5);
\draw[ color = DarkBlue, ,](0.0,0.0) -- (0.5,-0.5);
\draw[ color = DarkRed, very thin,](0.0,0.0) -- (0.5,0.0);
\draw[ color = black, very thick,](0.0,1.5) -- (0.5,1.5);
\draw[ color = DarkBlue, ,](0.0,1.0) -- (0.5,1.5);
\draw[ color = DarkRed, very thin,](0.0,1.0) -- (0.0,1.5);
\draw[ color = black, very thick,](0.0,1.5) -- (-0.5,1.5);
\draw[ color = DarkBlue, ,](0.0,1.0) -- (-0.5,1.5);
\draw[ color = DarkRed, very thin,](0.0,1.0) -- (0.0,1.5);
\draw[ color = black, very thick,](0.0,0.5) -- (0.5,0.5);
\draw[ color = DarkBlue, ,](0.0,1.0) -- (0.5,0.5);
\draw[ color = DarkRed, very thin,](0.0,1.0) -- (0.0,0.5);
\draw[ color = black, very thick,](0.0,0.5) -- (-0.5,0.5);
\draw[ color = DarkBlue, ,](0.0,1.0) -- (-0.5,0.5);
\draw[ color = DarkRed, very thin,](0.0,1.0) -- (0.0,0.5);
\draw[ color = black, very thick,](-0.5,1.0) -- (-0.5,1.5);
\draw[ color = DarkBlue, ,](0.0,1.0) -- (-0.5,1.5);
\draw[ color = DarkRed, very thin,](0.0,1.0) -- (-0.5,1.0);
\draw[ color = black, very thick,](-0.5,1.0) -- (-0.5,0.5);
\draw[ color = DarkBlue, ,](0.0,1.0) -- (-0.5,0.5);
\draw[ color = DarkRed, very thin,](0.0,1.0) -- (-0.5,1.0);
\draw[ color = black, very thick,](0.5,1.0) -- (0.5,1.5);
\draw[ color = DarkBlue, ,](0.0,1.0) -- (0.5,1.5);
\draw[ color = DarkRed, very thin,](0.0,1.0) -- (0.5,1.0);
\draw[ color = black, very thick,](0.5,1.0) -- (0.5,0.5);
\draw[ color = DarkBlue, ,](0.0,1.0) -- (0.5,0.5);
\draw[ color = DarkRed, very thin,](0.0,1.0) -- (0.5,1.0);
\draw[ color = black, very thick,](1.0,1.5) -- (1.5,1.5);
\draw[ color = DarkBlue, ,](1.0,1.0) -- (1.5,1.5);
\draw[ color = DarkRed, very thin,](1.0,1.0) -- (1.0,1.5);
\draw[ color = black, very thick,](1.0,1.5) -- (0.5,1.5);
\draw[ color = DarkBlue, ,](1.0,1.0) -- (0.5,1.5);
\draw[ color = DarkRed, very thin,](1.0,1.0) -- (1.0,1.5);
\draw[ color = black, very thick,](1.0,0.5) -- (1.5,0.5);
\draw[ color = DarkBlue, ,](1.0,1.0) -- (1.5,0.5);
\draw[ color = DarkRed, very thin,](1.0,1.0) -- (1.0,0.5);
\draw[ color = black, very thick,](1.0,0.5) -- (0.5,0.5);
\draw[ color = DarkBlue, ,](1.0,1.0) -- (0.5,0.5);
\draw[ color = DarkRed, very thin,](1.0,1.0) -- (1.0,0.5);
\draw[ color = black, very thick,](0.5,1.0) -- (0.5,1.5);
\draw[ color = DarkBlue, ,](1.0,1.0) -- (0.5,1.5);
\draw[ color = DarkRed, very thin,](1.0,1.0) -- (0.5,1.0);
\draw[ color = black, very thick,](0.5,1.0) -- (0.5,0.5);
\draw[ color = DarkBlue, ,](1.0,1.0) -- (0.5,0.5);
\draw[ color = DarkRed, very thin,](1.0,1.0) -- (0.5,1.0);
\draw[ color = black, very thick,](1.5,1.0) -- (1.5,1.5);
\draw[ color = DarkBlue, ,](1.0,1.0) -- (1.5,1.5);
\draw[ color = DarkRed, very thin,](1.0,1.0) -- (1.5,1.0);
\draw[ color = black, very thick,](1.5,1.0) -- (1.5,0.5);
\draw[ color = DarkBlue, ,](1.0,1.0) -- (1.5,0.5);
\draw[ color = DarkRed, very thin,](1.0,1.0) -- (1.5,1.0);
\draw[ color = black, very thick,](0.0,2.5) -- (0.5,2.5);
\draw[ color = DarkBlue, ,](0.0,2.0) -- (0.5,2.5);
\draw[ color = DarkRed, very thin,](0.0,2.0) -- (0.0,2.5);
\draw[ color = black, very thick,](0.0,2.5) -- (-0.5,2.5);
\draw[ color = DarkBlue, ,](0.0,2.0) -- (-0.5,2.5);
\draw[ color = DarkRed, very thin,](0.0,2.0) -- (0.0,2.5);
\draw[ color = black, very thick,](0.0,1.5) -- (0.5,1.5);
\draw[ color = DarkBlue, ,](0.0,2.0) -- (0.5,1.5);
\draw[ color = DarkRed, very thin,](0.0,2.0) -- (0.0,1.5);
\draw[ color = black, very thick,](0.0,1.5) -- (-0.5,1.5);
\draw[ color = DarkBlue, ,](0.0,2.0) -- (-0.5,1.5);
\draw[ color = DarkRed, very thin,](0.0,2.0) -- (0.0,1.5);
\draw[ color = black, very thick,](-0.5,2.0) -- (-0.5,2.5);
\draw[ color = DarkBlue, ,](0.0,2.0) -- (-0.5,2.5);
\draw[ color = DarkRed, very thin,](0.0,2.0) -- (-0.5,2.0);
\draw[ color = black, very thick,](-0.5,2.0) -- (-0.5,1.5);
\draw[ color = DarkBlue, ,](0.0,2.0) -- (-0.5,1.5);
\draw[ color = DarkRed, very thin,](0.0,2.0) -- (-0.5,2.0);
\draw[ color = black, very thick,](0.5,2.0) -- (0.5,2.5);
\draw[ color = DarkBlue, ,](0.0,2.0) -- (0.5,2.5);
\draw[ color = DarkRed, very thin,](0.0,2.0) -- (0.5,2.0);
\draw[ color = black, very thick,](0.5,2.0) -- (0.5,1.5);
\draw[ color = DarkBlue, ,](0.0,2.0) -- (0.5,1.5);
\draw[ color = DarkRed, very thin,](0.0,2.0) -- (0.5,2.0);
\draw[ color = black, very thick,](1.0,2.5) -- (1.5,2.5);
\draw[ color = DarkBlue, ,](1.0,2.0) -- (1.5,2.5);
\draw[ color = DarkRed, very thin,](1.0,2.0) -- (1.0,2.5);
\draw[ color = black, very thick,](1.0,2.5) -- (0.5,2.5);
\draw[ color = DarkBlue, ,](1.0,2.0) -- (0.5,2.5);
\draw[ color = DarkRed, very thin,](1.0,2.0) -- (1.0,2.5);
\draw[ color = black, very thick,](1.0,1.5) -- (1.5,1.5);
\draw[ color = DarkBlue, ,](1.0,2.0) -- (1.5,1.5);
\draw[ color = DarkRed, very thin,](1.0,2.0) -- (1.0,1.5);
\draw[ color = black, very thick,](1.0,1.5) -- (0.5,1.5);
\draw[ color = DarkBlue, ,](1.0,2.0) -- (0.5,1.5);
\draw[ color = DarkRed, very thin,](1.0,2.0) -- (1.0,1.5);
\draw[ color = black, very thick,](0.5,2.0) -- (0.5,2.5);
\draw[ color = DarkBlue, ,](1.0,2.0) -- (0.5,2.5);
\draw[ color = DarkRed, very thin,](1.0,2.0) -- (0.5,2.0);
\draw[ color = black, very thick,](0.5,2.0) -- (0.5,1.5);
\draw[ color = DarkBlue, ,](1.0,2.0) -- (0.5,1.5);
\draw[ color = DarkRed, very thin,](1.0,2.0) -- (0.5,2.0);
\draw[ color = black, very thick,](1.5,2.0) -- (1.5,2.5);
\draw[ color = DarkBlue, ,](1.0,2.0) -- (1.5,2.5);
\draw[ color = DarkRed, very thin,](1.0,2.0) -- (1.5,2.0);
\draw[ color = black, very thick,](1.5,2.0) -- (1.5,1.5);
\draw[ color = DarkBlue, ,](1.0,2.0) -- (1.5,1.5);
\draw[ color = DarkRed, very thin,](1.0,2.0) -- (1.5,2.0);
\draw[ color = black, very thick,](2.0,2.5) -- (2.5,2.5);
\draw[ color = DarkBlue, ,](2.0,2.0) -- (2.5,2.5);
\draw[ color = DarkRed, very thin,](2.0,2.0) -- (2.0,2.5);
\draw[ color = black, very thick,](2.0,2.5) -- (1.5,2.5);
\draw[ color = DarkBlue, ,](2.0,2.0) -- (1.5,2.5);
\draw[ color = DarkRed, very thin,](2.0,2.0) -- (2.0,2.5);
\draw[ color = black, very thick,](2.0,1.5) -- (2.5,1.5);
\draw[ color = DarkBlue, ,](2.0,2.0) -- (2.5,1.5);
\draw[ color = DarkRed, very thin,](2.0,2.0) -- (2.0,1.5);
\draw[ color = black, very thick,](2.0,1.5) -- (1.5,1.5);
\draw[ color = DarkBlue, ,](2.0,2.0) -- (1.5,1.5);
\draw[ color = DarkRed, very thin,](2.0,2.0) -- (2.0,1.5);
\draw[ color = black, very thick,](1.5,2.0) -- (1.5,2.5);
\draw[ color = DarkBlue, ,](2.0,2.0) -- (1.5,2.5);
\draw[ color = DarkRed, very thin,](2.0,2.0) -- (1.5,2.0);
\draw[ color = black, very thick,](1.5,2.0) -- (1.5,1.5);
\draw[ color = DarkBlue, ,](2.0,2.0) -- (1.5,1.5);
\draw[ color = DarkRed, very thin,](2.0,2.0) -- (1.5,2.0);
\draw[ color = black, very thick,](2.5,2.0) -- (2.5,2.5);
\draw[ color = DarkBlue, ,](2.0,2.0) -- (2.5,2.5);
\draw[ color = DarkRed, very thin,](2.0,2.0) -- (2.5,2.0);
\draw[ color = black, very thick,](2.5,2.0) -- (2.5,1.5);
\draw[ color = DarkBlue, ,](2.0,2.0) -- (2.5,1.5);
\draw[ color = DarkRed, very thin,](2.0,2.0) -- (2.5,2.0);
\draw[->,  color = purple,] (0.875,1.625) -- (0.875,1.375);
\draw[->,  color = purple,] (0.875,1.375) -- (0.625,1.125);
\draw[->,  color = purple,] (0.625,1.125) -- (0.375,1.125);
\draw[->,  color = purple,] (0.375,1.125) -- (0.375,0.875);
\draw[->,  color = purple,] (0.375,0.875) -- (0.125,0.625);
\draw[->,  color = purple,] (0.125,0.625) -- (0.125,0.375);
\end{tikzpicture}}&
    \scalebox{.9}{\input{Fig/minimalAntisortingAlcoveWalk-ambient-G2}}\\
    %
    \scalebox{.9}{\begin{tikzpicture}[baseline=(current bounding box.east)]
\draw[->,  color = black,] (0.0,0.0) -- (0.0,-1.732050808) node[at end, auto=right] {$\alpha^\vee_{0}$};
\draw[->,  color = DarkBlue,] (0.0,0.0) -- (1.5,0.8660254038) node[at end, auto=right] {$\alpha^\vee_{1}$};
\draw[->,  color = DarkRed,] (0.0,0.0) -- (-1.5,0.8660254038) node[at end, auto=right] {$\alpha^\vee_{2}$};
\draw[ color = black, very thick,](0.5,0.8660254038) -- (-0.5,0.8660254038);
\draw[ color = DarkBlue, ,](0.0,0.0) -- (-0.5,0.8660254038);
\draw[ color = DarkRed, very thin,](0.0,0.0) -- (0.5,0.8660254038);
\draw[ color = black, very thick,](0.5,0.8660254038) -- (1.0,0.0);
\draw[ color = DarkBlue, ,](0.0,0.0) -- (1.0,0.0);
\draw[ color = DarkRed, very thin,](0.0,0.0) -- (0.5,0.8660254038);
\draw[ color = black, very thick,](0.5,-0.8660254038) -- (1.0,0.0);
\draw[ color = DarkBlue, ,](0.0,0.0) -- (1.0,0.0);
\draw[ color = DarkRed, very thin,](0.0,0.0) -- (0.5,-0.8660254038);
\draw[ color = black, very thick,](0.5,-0.8660254038) -- (-0.5,-0.8660254038);
\draw[ color = DarkBlue, ,](0.0,0.0) -- (-0.5,-0.8660254038);
\draw[ color = DarkRed, very thin,](0.0,0.0) -- (0.5,-0.8660254038);
\draw[ color = black, very thick,](-1.0,0.0) -- (-0.5,0.8660254038);
\draw[ color = DarkBlue, ,](0.0,0.0) -- (-0.5,0.8660254038);
\draw[ color = DarkRed, very thin,](0.0,0.0) -- (-1.0,0.0);
\draw[ color = black, very thick,](-1.0,0.0) -- (-0.5,-0.8660254038);
\draw[ color = DarkBlue, ,](0.0,0.0) -- (-0.5,-0.8660254038);
\draw[ color = DarkRed, very thin,](0.0,0.0) -- (-1.0,0.0);
\draw[->,  color = purple,thick] (0.0,-0.5773502692) -- (0.0,0.5773502692);
\end{tikzpicture}}&
    \scalebox{1.2}{\begin{tikzpicture}[baseline=(current bounding box.east)]
\draw[->,  color = black,] (0.0,0.0) -- (0.0,-1.0) node[at end, auto=right] {$\alpha^\vee_{0}$};
\draw[->,  color = DarkBlue,] (0.0,0.0) -- (-1.0,1.0) node[at end, auto=left] {$\alpha^\vee_{1}$};
\draw[->,  color = DarkRed,] (0.0,0.0) -- (1.0,0.0) node[at end, auto=right] {$\alpha^\vee_{2}$};
\draw[ color = black, very thick,](0.0,0.5) -- (0.5,0.5);
\draw[ color = DarkBlue, ,](0.0,0.0) -- (0.5,0.5);
\draw[ color = DarkRed, very thin,](0.0,0.0) -- (0.0,0.5);
\draw[ color = black, very thick,](0.0,0.5) -- (-0.5,0.5);
\draw[ color = DarkBlue, ,](0.0,0.0) -- (-0.5,0.5);
\draw[ color = DarkRed, very thin,](0.0,0.0) -- (0.0,0.5);
\draw[ color = black, very thick,](0.0,-0.5) -- (0.5,-0.5);
\draw[ color = DarkBlue, ,](0.0,0.0) -- (0.5,-0.5);
\draw[ color = DarkRed, very thin,](0.0,0.0) -- (0.0,-0.5);
\draw[ color = black, very thick,](0.0,-0.5) -- (-0.5,-0.5);
\draw[ color = DarkBlue, ,](0.0,0.0) -- (-0.5,-0.5);
\draw[ color = DarkRed, very thin,](0.0,0.0) -- (0.0,-0.5);
\draw[ color = black, very thick,](-0.5,0.0) -- (-0.5,0.5);
\draw[ color = DarkBlue, ,](0.0,0.0) -- (-0.5,0.5);
\draw[ color = DarkRed, very thin,](0.0,0.0) -- (-0.5,0.0);
\draw[ color = black, very thick,](-0.5,0.0) -- (-0.5,-0.5);
\draw[ color = DarkBlue, ,](0.0,0.0) -- (-0.5,-0.5);
\draw[ color = DarkRed, very thin,](0.0,0.0) -- (-0.5,0.0);
\draw[ color = black, very thick,](0.5,0.0) -- (0.5,0.5);
\draw[ color = DarkBlue, ,](0.0,0.0) -- (0.5,0.5);
\draw[ color = DarkRed, very thin,](0.0,0.0) -- (0.5,0.0);
\draw[ color = black, very thick,](0.5,0.0) -- (0.5,-0.5);
\draw[ color = DarkBlue, ,](0.0,0.0) -- (0.5,-0.5);
\draw[ color = DarkRed, very thin,](0.0,0.0) -- (0.5,0.0);
\draw[->,  color = purple,] (-0.125,-0.375) -- (-0.125,0.375);
\draw[->,  color = purple,] (-0.125,0.375) -- (-0.375,0.125);
\draw[->,  color = purple,] (-0.375,0.125) -- (0.375,0.125);
\draw[->,  color = purple,] (0.375,0.125) -- (0.375,-0.125);
\draw[->,  color = purple,] (0.375,-0.125) -- (0.125,-0.375);
\draw[->,  color = purple,] (0.125,-0.375) -- (0.125,0.375);
\end{tikzpicture}}&
    \scalebox{.9}{\begin{tikzpicture}[baseline=(current bounding box.east)]
\draw[->,  color = black,] (0.0,0.0) -- (0.0,-1.732050808) node[at end, auto=right] {$\alpha^\vee_{0}$};
\draw[->,  color = DarkBlue,] (0.0,0.0) -- (-1.5,0.8660254038) node[at end, auto=right] {$\alpha^\vee_{1}$};
\draw[->,  color = DarkRed,] (0.0,0.0) -- (3.0,0.0) node[at end, auto=right] {$\alpha^\vee_{2}$};
\draw[ color = black, very thick,](0.0,0.8660254038) -- (0.5,0.8660254038);
\draw[ color = DarkBlue, ,](0.0,0.0) -- (0.5,0.8660254038);
\draw[ color = DarkRed, very thin,](0.0,0.0) -- (0.0,0.8660254038);
\draw[ color = black, very thick,](0.0,0.8660254038) -- (-0.5,0.8660254038);
\draw[ color = DarkBlue, ,](0.0,0.0) -- (-0.5,0.8660254038);
\draw[ color = DarkRed, very thin,](0.0,0.0) -- (0.0,0.8660254038);
\draw[ color = black, very thick,](0.0,-0.8660254038) -- (0.5,-0.8660254038);
\draw[ color = DarkBlue, ,](0.0,0.0) -- (0.5,-0.8660254038);
\draw[ color = DarkRed, very thin,](0.0,0.0) -- (0.0,-0.8660254038);
\draw[ color = black, very thick,](0.0,-0.8660254038) -- (-0.5,-0.8660254038);
\draw[ color = DarkBlue, ,](0.0,0.0) -- (-0.5,-0.8660254038);
\draw[ color = DarkRed, very thin,](0.0,0.0) -- (0.0,-0.8660254038);
\draw[ color = black, very thick,](0.75,0.4330127019) -- (0.5,0.8660254038);
\draw[ color = DarkBlue, ,](0.0,0.0) -- (0.5,0.8660254038);
\draw[ color = DarkRed, very thin,](0.0,0.0) -- (0.75,0.4330127019);
\draw[ color = black, very thick,](0.75,0.4330127019) -- (1.0,0.0);
\draw[ color = DarkBlue, ,](0.0,0.0) -- (1.0,0.0);
\draw[ color = DarkRed, very thin,](0.0,0.0) -- (0.75,0.4330127019);
\draw[ color = black, very thick,](-0.75,0.4330127019) -- (-0.5,0.8660254038);
\draw[ color = DarkBlue, ,](0.0,0.0) -- (-0.5,0.8660254038);
\draw[ color = DarkRed, very thin,](0.0,0.0) -- (-0.75,0.4330127019);
\draw[ color = black, very thick,](-0.75,0.4330127019) -- (-1.0,0.0);
\draw[ color = DarkBlue, ,](0.0,0.0) -- (-1.0,0.0);
\draw[ color = DarkRed, very thin,](0.0,0.0) -- (-0.75,0.4330127019);
\draw[ color = black, very thick,](-0.75,-0.4330127019) -- (-1.0,0.0);
\draw[ color = DarkBlue, ,](0.0,0.0) -- (-1.0,0.0);
\draw[ color = DarkRed, very thin,](0.0,0.0) -- (-0.75,-0.4330127019);
\draw[ color = black, very thick,](-0.75,-0.4330127019) -- (-0.5,-0.8660254038);
\draw[ color = DarkBlue, ,](0.0,0.0) -- (-0.5,-0.8660254038);
\draw[ color = DarkRed, very thin,](0.0,0.0) -- (-0.75,-0.4330127019);
\draw[ color = black, very thick,](0.75,-0.4330127019) -- (1.0,0.0);
\draw[ color = DarkBlue, ,](0.0,0.0) -- (1.0,0.0);
\draw[ color = DarkRed, very thin,](0.0,0.0) -- (0.75,-0.4330127019);
\draw[ color = black, very thick,](0.75,-0.4330127019) -- (0.5,-0.8660254038);
\draw[ color = DarkBlue, ,](0.0,0.0) -- (0.5,-0.8660254038);
\draw[ color = DarkRed, very thin,](0.0,0.0) -- (0.75,-0.4330127019);
\draw[->,  color = purple,] (-0.25,-0.7216878365) -- (-0.25,0.7216878365);
\draw[->,  color = purple,] (-0.25,0.7216878365) -- (-0.5,0.5773502692);
\draw[->,  color = purple,] (-0.5,0.5773502692) -- (-0.75,0.1443375673);
\draw[->,  color = purple,] (-0.75,0.1443375673) -- (-0.75,-0.1443375673);
\draw[->,  color = purple,] (-0.75,-0.1443375673) -- (-0.5,-0.5773502692);
\draw[->,  color = purple,] (-0.5,-0.5773502692) -- (0.75,0.1443375673);
\draw[->,  color = purple,] (0.75,0.1443375673) -- (0.75,-0.1443375673);
\draw[->,  color = purple,] (0.75,-0.1443375673) -- (0.5,-0.5773502692);
\draw[->,  color = purple,] (0.5,-0.5773502692) -- (0.25,-0.7216878365);
\draw[->,  color = purple,] (0.25,-0.7216878365) -- (0.25,0.7216878365);
\end{tikzpicture}}\\
    %
    \dynkinAIIa & \dynkinCIIa & \dynkinGIIa\\[4mm]
    %
    \widetilde A_2 = A_2^{(1)}
    & \widetilde C_2 = C_2^{(1)} & \widetilde G_2 = G_2^{(1)}
  \end{array}$
  \\\bigskip
  \parbox{\textwidth}{
   \textbf{En haut:} Graphe de l'action de niveau zéro de
    $\pi_0,\pi_1,\dots,\pi_n$ sur le groupe de  Weyl classique
    $\clW$ (notation par permutations signées avec $\underline2:=-2$).

    \textbf{Milieu:} 
    Les alcôves dans l'espace ambiant, avec un plus court chemin
    d'alcôve descendant depuis une alcôve $w(A)$ dans la chambre
    dominante telle que $\cl(w)=\Wmax$ jusqu'à l'alcôve fondamentale
    $A$. 

    \textbf{En bas:} Le graphe du haut peut être réalisé
    géométriquement par le \emph{tore de Steinberg}, quotient des
    alcôves par les translations, ou de manière équivalente par
    identification des faces opposées du polygone fondamental. 
    Le chemin d'alcôve de la figure du milieu devient alors un chemin
    de retour depuis la chambre anti-dominante $\Wmax(A)$ vers la
    chambre dominante $A$.}

  \caption[Transitivité de l'action de niveau zéro de la $0$-algèbre de
    Hecke affine]{Transitivité de l'action de niveau zéro de la $0$-algèbre de
    Hecke affine $\heckeW{0}$ sur le groupe de Weyl classique $\clW$.}
  \label{figure.antisorting}
\end{figure}

\subsubsection{Représentations de série principale de l'algèbre de Hecke affine}
\label{section.calibree}

La fin de la démonstration du théorème~\ref{theorem.quotient}, pour
$q$ générique, est une généralisation directe du type $A$. Pour
réduire aux petites racines de l'unité les valeurs de $q$
exceptionnelles, nous avons utilisé, sur la suggestion d'Arun Ram, une
autre approche.

Le point de départ est que $\Wmax$ dans $\kclW$ est un vecteur propre
pour le tore commutatif $\CC[Y^{\alpha_i}]$ de l'algèbre de Hecke
affine engendré par les opérateurs de Cherednic $Y^{\alpha_i}$. On
peut alors utiliser une construction classique, due elle aussi à
Cherednic, qui permet de construire de nouveaux vecteurs propres grâce
aux \emph{opérateurs d'entrelacement} $\tau_i$ (des déformations des
$T_i$ qui commutent presque avec les $Y_j$). Cela revient à utiliser
un graphe de Yang-Baxter pour un bon choix de paramètres spectraux
(voir~\cite[section~10.7]{Lascoux.2003.CBMS}).

Nous avons alors montré que, lorsque $q$ n'est pas une petite racine
de l'unité, les valeurs propres sont suffisamment différentes
(représentation \emph{calibrée}) pour garantir que l'on a diagonalisé
simultanément l'action des $Y_i$ sur $\kW$. Plus précisément, la
représentation de niveau zéro de l'algèbre de Hecke affine est un cas
particulier de \emph{représentation de série principale} $M(t)$ (voir
par exemple~\cite[section 2.5]{Ram.2003}), pour le caractère
$t:Y^{\lc}\mapsto q^{-\height(\lc)}$. La vérification de la
surjectivité du morphisme $\cl$ se fait alors grâce à un simple calcul
de dimension; celui-ci relie le nombre de vecteurs propres où
s'annulent les opérateurs $\tau_i^2$ avec la combinatoire des
descentes de $\clW$.

\subsection{Exploration informatique}

Dans cette recherche, l'ordinateur a été principalement un outil
d'exploration: où y a-t-il de la structure? Quelles conjectures faire?
Sur quelles propriétés s'appuyer? Quel est le bon point de vue?  Par
exemple, la vérification de l'existence d'un algorithme de tri-antitri
récursif pour tous les types exceptionnels a fortement motivé la
recherche d'une démonstration géométrique. Au final, la plupart des
démonstrations sont élémentaires. Par exemple, il n'est pas difficile de
dérouler la théorie des représentations de l'algèbre de Hecke groupe,
une fois que l'on a vu que la combinatoire sous-jacente est celle des
descentes dans le groupe de Coxeter.  Ce travail a été aussi pour moi
l'occasion, et c'était l'un des objectifs, de comprendre et d'implanter
les systèmes de racines, groupes de Coxeter et de Weyl. Des outils
préexistaient pour les premiers, principalement dans le cas fini, par
exemple dans \gap ou \maple. L'algorithmique en est fortement
inspirée. La conception en revanche est complètement nouvelle. Elle
permet de manipuler simultanément et de manière naturelle les
différentes réalisations classiques du réseau des racines ou des
poids. De plus, la majorité du code est complètement générique; il ne
dépend que des données de la matrice de Cartan. De ce fait, il
s'applique aussi, lorsque cela fait sens, aux cas affine ou de
Kac-Moody, et à terme au cas non cristallographique. Cela a fourni une
base solide pour l'implantation des cristaux, chemins d'alcôves et
algèbres de Hecke affines.

\subsection{Perspectives}
La structure riche de l'algèbre de Hecke groupe est maintenant bien
comprise. L'existence de trois définitions équivalentes (par
générateurs, comme algèbre d'opérateurs préservant certaines
(anti)symétries, par quotient de l'algèbre de Hecke affine) en fait un
objet naturel, qui éclaire les liens entres les représentations de la
$0$-algèbre de Hecke et de l'algèbre de Hecke affine. En dehors de
cela, est-elle utile? Nous pensons qu'elle est susceptible de fournir
un modèle combinatoire simple (\emph{via} les classes de descentes de
$\clW$) pour mieux comprendre certaines représentations de dimension
$|\clW|$ de l'algèbre de Hecke affine (coinvariants, harmoniques pour
l'algèbre de Steenrod, etc.).

Nous présentons ici plusieurs projets de recherche en cours qui dans
cette direction, ou tendent à généraliser la structure de l'algèbre de
Hecke groupe à des algèbres et monoïdes proches.

\subsubsection{Algèbres de Hecke affines aux racines de l'unité}
L'énoncé du théorème~\ref{theorem.quotient} soulève immédiatement le
problème suivant.
\begin{problem}
  Déterminer l'ensemble des racines de l'unité $q$ pour lesquelles le
  morphisme \mbox{$\cl:\heckeW{q}\mapsto\heckeWW[\clW]$} n'est pas
  surjectif.
\end{problem}
Au printemps dernier, j'ai donné l'exploration informatique de ce
problème comme sujet de Master 2 à Nicolas Borie. Les calculs sont
lourds; par exemple, traiter $\sg$ nécessite de faire de l'algèbre
linéaire sur des vecteurs qui sont des matrices $n!\times n!$, le tout
sur une extension algébrique des rationnels. L'essentiel de son
travail a été d'élaborer des stratégies mathématiques pour exploiter
au mieux les résultats partiels.

À l'instant, il semble que la borne du théorème~\ref{theorem.quotient}
soit assez précise: la plupart des petites racines de l'unité, mais
pas toutes, semblent donner lieu à un morphisme non surjectif.  C'est
en particulier toujours le cas pour $q=-1$. Parallèlement au sujet
principal de thèse que je lui ai confié, Nicolas Borie va continuer à
étudier ce problème. La résolution de difficultés techniques dans
\sage devrait permettre de mener les calculs suffisamment loin pour
établir une conjecture exacte.  D'un autre côté, le cas $q=-1$ doit
pouvoir être démontré, et cela devrait donner une idée plus fine de la
difficulté du cas général.

\subsubsection{Monoïde des fonctions décroissantes de type $C$}

Lorsque nous avions étudié la théorie des représentations de
$\heckesg$, nous disposions depuis peu d'un outil, mis au point par
Florent Hivert, calculant automatiquement la théorie des
représentations des premiers étages d'une tour d'algèbres. Pour
expérimenter, nous avions alors considéré quelques tours d'algèbres
jouets comme l'algèbre $\ndfa n$ (resp. $\ndpfa n$) du monoïde des
fonctions (resp. fonctions de parking) croissantes. À notre grande
surprise, celles-ci se sont naturellement intégrées dans un
diagramme. Cela nous a permis de définir des représentations de
$(\heckesg)_n$ et $(H_n(q))_n$ sur les puissances extérieures de la
représentation naturelle, et de retrouver comme cas particulier la
tour d'algèbres de Temperley-Lieb $(\operatorname{TL}_n)_n$:
\medskip\smallskip
\begin{equation}
  \vcenter{
  \xymatrix@R=1cm@C=0.52cm{
    \hecke{-1}        \ar@{->>}[d] \ar@{^(->}@/^4ex/[rrrr] &
    \hecke{0}         \ar@{->>}[d] \ar@{^(->}@/^3ex/[rrr]  &
    \hecke{1}=\sga n  \ar@{->>}[d] \ar@{^(->}@/^2ex/[rr]   &
    \hecke{q}         \ar@{->>}[d] \ar@{^(->}    [r]       &
    \heckesg          \ar@{->>}[d]\\
    \operatorname{TL}_n                \ar@{^(->}@/_4ex/[rrrr]&
    \ndpfa n                                         \ar@{^(->}@/_3ex/[rrr] &
    \sga n     \circlearrowleft\! \bigwedge\!\CC^n \ar@{^(->}@/_2ex/[rr]  &
    \hecke[n]{q}\circlearrowleft\! \bigwedge\!\CC^n \ar@{^(->}    [r]      &
    \ndfa n
  }}
\end{equation}
\medskip

Avec Tom Denton, doctorant à UC Davis sous la direction de Anne
Schilling, nous essayons actuellement de généraliser ce schéma à
d'autres types. Cela a déjà fait apparaître des fonctions de parking
croissantes signées de type C, dont le monoïde est auto-injectif, et
dont la dimension semble être donnée par la suite A086618 de
l'encyclopédie des suites d'entiers~\cite{Sloane}.

\subsubsection{Monoïde des opérateurs de tri-antitri}

Avec Anne Schilling, Jean-Christophe Novelli et Florent Hivert, nous
avons entrepris l'étude du monoïde $M(W) = \langle
\pi_i,\opi_i\rangle$ (et non de l'algèbre) engendré par les opérateurs
$\pi_i$ et $\opi_i$. À cela deux motivations. D'une part, Philippe
Gaucher nous a contactés après être avoir reconnu
dans~\cite{Hivert_Thiery.HeckeGroup.2007} les premiers termes de la
suite $|M(\sg)|$. Selon lui, une meilleure compréhension de ce monoïde
pourrait avoir des applications dans la modélisation de processus
concurrents par des méthodes de topologie
algébrique~\cite{Gaucher.2008}.

D'autre part, il semblerait que ce monoïde ait une structure
intrinsèque, même si sa taille n'est pas connue. En effet, le calcul
dans les petits cas de sa théorie des représentations indique que les
modules simples seraient indexés par les éléments du groupe de Coxeter
et que, en type $A$, la somme des dimensions des modules simples
serait donnée par la suite A006245~\cite{Sloane} qui compte les
réseaux de tris réduits. Une compréhension complète de sa théorie des
représentations devrait permettre d'obtenir une formule par sommation
pour sa taille $|M(W)|$.  Elle pourrait aussi suggérer une description
plus conceptuelle de ce monoïde, prérequise pour mieux le comprendre.

\subsubsection{Polynômes de Macdonald non symétriques}

Les opérateurs d'entrelacement ont été originellement introduits pour
construire les polynômes de Macdonald non symétriques, ces derniers
étant définis comme vecteurs propres simultanés des opérateurs
$Y^{\alpha_i}$ de Cherednick de l'algèbre de Hecke affine. Cela
suggère bien entendu d'interpréter la base de vecteurs propres de
$\kW$ obtenus dans la section~\ref{section.calibree} comme polynômes
de Macdonald. Pour que cette interprétation soit naturelle, il
faudrait, en type $A$, réaliser $\ksg n$ comme module quotient de
l'algèbre des polynômes sous l'action de l'algèbre de Hecke affine par
différences divisées isobares. Il est probable que cette première
étape découle directement de la construction de la représentation
polynomiale de l'algèbre de Hecke doublement affine (coinvariants).

Il reste à comprendre comment agissent dans cette réalisation les
différents opérateurs de l'algèbre de Hecke groupe, à commencer par
les permutations, mais surtout $T_0$.  En effet, la construction des
polynômes de Macdonald par opérateurs d'entrelacement
$\tau_1,\dots,\tau_n$ ne donne comme sous-produit qu'une description
simple de l'action de $T_1,\dots,T_n$ sur les polynômes de Macdonald,
mais pas de $T_0$. En revanche, dans l'action de l'algèbre de Hecke
groupe sur $\kW$, les opérateurs $T_0,T_1,\dots,T_n$ jouent des rôles
très symétriques, et cela pourrait éclairer le rôle de $T_0$.

\section{Opérateurs de promotion sur les graphes cristallins affines}
\label{section.cristaux}

\begin{figure}[h]
  \centering
  \scalebox{0.6}{\input{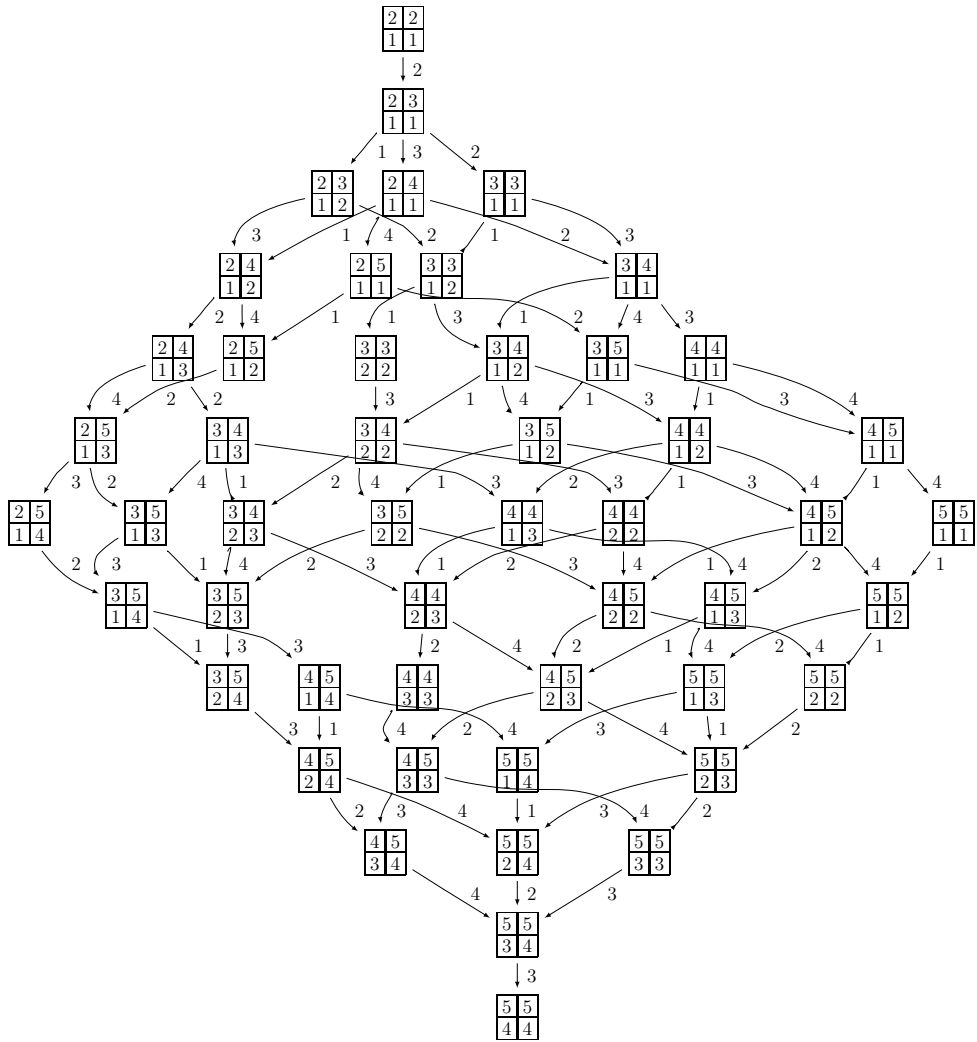}}
  \caption[Le graphe cristallin de forme $(2,2)$ en type $A_4$]
  {Le graphe cristallin de forme $(2,2)$ en type $A_4$; une
    flèche $f_i$ remplace un $i$ par un $i+1$, selon une règle
    contrainte par la préservation des conditions de croissance au
    sens large le long des lignes et au sens strict le long des colonnes}
  \label{figure.cristal}
\end{figure}

En parallèle avec notre travail sur les relations entre algèbres de
Hecke groupe et algèbres de Hecke affine, mon séjour à l'Université de
Californie à Davis a été l'occasion d'une collaboration avec Jason
Bandlow et Anne Schilling sur les \emph{graphes cristallins}. Ce sont
des graphes orientés, comme celui de la figure~\ref{figure.cristal},
qui ont été introduits en théorie des représentations par Masaki
Kashiwara afin d'encoder l'essentiel de la structure des modules pour
les groupes quantiques $U_q(G)$ lorsque $q$ tend vers $0$. Ces graphes
jouent en particulier un rôle important en physique mathématique, en
lien avec les modèles intégrables sur réseaux (voir, par
exemple~\cite{Hatayama_Kuniba_Okado_Takagi_Tsuboi.2001}).  Les arêtes
indexées par $i$ décrivent l'action des opérateurs descendants $f_i$,
les opérateurs montants $e_i$ étant les inverses locaux.

La théorie pour les modules de plus haut poids (ce qui est toujours le
cas lorsque $G$ est une algèbre de Lie de dimension finie) est
essentiellement complète: les modules possèdent toujours une base
cristalline, on a des modèles combinatoires (par exemple les tableaux
en type A, les chemins de Littlemann en d'autres types), ainsi qu'une
caractérisation locale complète des graphes qui sont
cristallins~\cite{Stembridge.2003}. En particulier, tout graphe
cristallin $B(\lambda)$ peut être obtenu par produit tensoriel de
graphes élémentaires correspondant aux poids fondamentaux. En
revanche, pour les algèbres de Lie affine (ou plus généralement de
Kac-Moody), de nombreuses questions restent ouvertes comme: quels
modules admettent une base cristalline, quels graphes sont
cristallins, etc. Par exemple, une partie des travaux récents d'Anne
Schilling concerne la construction de modèles combinatoires pour les
modules $B^{r,s}$ de
Kirillov-Reshetikhin~\cite{Fourier_Schilling_Shimozono.2007,
  Schilling.2007,Okado_Schilling.2008}, qui sont en lien avec les
modèles intégrables sur réseaux en physique mathématique, et semblent
jouer le rôle des graphes cristallins élémentaires.
\begin{conjecture}[{Kashiwara~\cite[Introduction]{Kashiwara.2005}}]
  \label{conjecture.Kashiwara}
  Tout bon cristal affine (i.e. provenant d'un $U_q(\widehat{sl_n})$
  module) de dimension finie est un produit tensoriel de cristaux de
  Kirillov-Reshetikhin.
\end{conjecture}
Notre travail va en direction de cette conjecture, en type $A$,
sachant que jusqu'ici seul le cas de $\widehat{sl_2}$ a été validé
par un résultat analogue de Chari et Pressley pour les
$U_q(\widehat{sl_2})$-modules~\cite{Chari_Pressley.1995.QuantumAffineAlgebrasAndTheirRepresentations}.

En type affine $A_n^{(1)}$, il existe un modèle combinatoire simple
pour les graphes cristallins des modules
$B^{r,s}$~\cite{Shimozono.2002}. On considère le graphe cristallin
usuel $B(s^r)$ de type $A_n$ et de plus haut poids la partition
$(s^r)$. Les sommets sont indexés par les tableaux semi-standard de
forme rectangle, sur l'alphabet $(1,\dots,n+1)$. Pour $i=1,\dots,n$,
les opérateurs $f_i$ transforment une lettre $i$ en une lettre $i+1$ et
réciproquement pour les opérateurs $e_i$. La règle exacte est
complètement contrainte par les conditions de croissance sur les
tableaux semi-standard. Pour étendre ce graphe en un graphe affine, il
faut définir l'action des opérateurs $f_0$. Pour cela, on utilise la
symétrie en rotation du diagramme de Dynkin en type $A_n^{(1)}$ (ici
$A_n^{(5)}$):
\begin{displaymath}
  \dynkinAVa
\end{displaymath}
Cette symétrie se traduit combinatoirement par l'opérateur de
promotion $\pr$ sur les tableaux semi-standard, introduit par
Schützenberger~\cite{Schuetzenberger.1972} à l'aide du jeu de taquin
dans le cadre général des ordres partiels. Cela provient du fait que
le jeu de taquin et les opérateurs cristallins commutent. L'opérateur
de promotion contraint complètement l'opération affine $f_0$ qui est
définie par $f_0 := \pr^{-1} \circ f_1 \circ \pr$. L'opérateur de
promotion est d'ordre $n+1$ sur les tableaux de forme $\lambda$ si et
seulement si $\lambda$ est de forme rectangle~\cite{Haiman.1992}; cela
explique le rôle spécial de ces formes.

Notre objectif est de généraliser cela aux produits tensoriels.
\begin{conjecture}[Bandlow, Schilling, T.~\cite{Bandlow_Schilling_Thiery.2008.Promotion}]
  \label{conjecture.promotion.unique}
  Il existe un unique opérateur de promotion, et donc une unique bonne
  structure affine de type $A_n^{(1)}$, sur les \emph{produits
    tensoriels} de graphes cristallins de type $A_n$ et de forme
  rectangle.
\end{conjecture}
Dans~\cite{Bandlow_Schilling_Thiery.2008.Promotion} nous démontrons
que cette conjecture est vraie pour un produit tensoriel à deux termes
en type $A_n^{(1)}$ pour $n\geq 2$, ainsi que pour de nombreux
exemples.

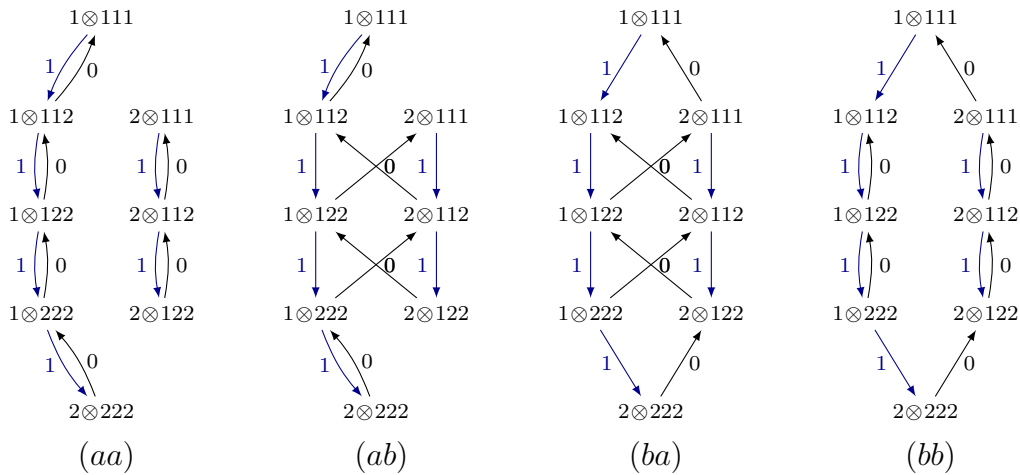
\begin{figure}
  \begin{bigcenter}
    \begin{tabular}{c@{\qquad}c@{\qquad}c@{\qquad}c}
      \begin{tikzpicture}[>=latex,join=bevel,scale=.5]
\tiny%
  \node (N_1) at (50bp,306bp) [draw,draw=none] {$1 \!\otimes\!111$};
  \node (N_2) at (5bp,232bp) [draw,draw=none] {$1 \!\otimes\!112$};
  \node (N_3) at (5bp,158bp) [draw,draw=none] {$1 \!\otimes\!122$};
  \node (N_4) at (95bp,232bp) [draw,draw=none] {$2 \!\otimes\!111$};
  \node (N_5) at (95bp,158bp) [draw,draw=none] {$2 \!\otimes\!112$};
  \node (N_6) at (5bp,84bp) [draw,draw=none] {$1 \!\otimes\!222$};
  \node (N_7) at (95bp,84bp) [draw,draw=none] {$2 \!\otimes\!122$};
\node (N_8) at (50bp,10bp) [draw,draw=none] {$2 \!\otimes\!222$};
  \draw [->,DarkBlue] (N_6) to [bend right=10] node [left] {$1$} (N_8);
  \draw [<-] (N_6) to [bend left=10] node [right] {$0$} (N_8);
  \draw [->,DarkBlue] (N_1) to [bend right=10] node [left] {$1$} (N_2);
  \draw [<-] (N_1) to [bend left=10] node [right] {$0$} (N_2);
  \draw [->,DarkBlue] (N_3) to [bend right=10] node [left] {$1$} (N_6);
  \draw [<-] (N_3) to [bend left=10] node [right] {$0$} (N_6);
  \draw [->,DarkBlue] (N_4) to [bend right=10] node [left] {$1$} (N_5);
  \draw [<-] (N_4) to [bend left=10] node [right] {$0$} (N_5);
  \draw [->,DarkBlue] (N_5) to [bend right=10] node [left] {$1$} (N_7);
  \draw [<-] (N_5) to [bend left=10] node [right] {$0$} (N_7);
  \draw [->,DarkBlue] (N_2) to [bend right=10] node [left] {$1$} (N_3);
  \draw [<-] (N_2) to [bend left=10] node [right] {$0$} (N_3);
\end{tikzpicture} &
      \begin{tikzpicture}[>=latex,join=bevel,scale=.5]
\tiny%
  \node (N_1) at (50bp,306bp) [draw,draw=none] {$1 \!\otimes\!111$};
  \node (N_2) at (5bp,232bp) [draw,draw=none] {$1 \!\otimes\!112$};
  \node (N_3) at (5bp,158bp) [draw,draw=none] {$1 \!\otimes\!122$};
  \node (N_4) at (95bp,232bp) [draw,draw=none] {$2 \!\otimes\!111$};
  \node (N_5) at (95bp,158bp) [draw,draw=none] {$2 \!\otimes\!112$};
  \node (N_6) at (5bp,84bp) [draw,draw=none] {$1 \!\otimes\!222$};
  \node (N_7) at (95bp,84bp) [draw,draw=none] {$2 \!\otimes\!122$};
\node (N_8) at (50bp,10bp) [draw,draw=none] {$2 \!\otimes\!222$};

  \draw [<-] (N_5) to [] node [right] {$0$} (N_6);

  \draw [->,DarkBlue] (N_6) to [bend right=10] node [left] {$1$} (N_8);
  \draw [<-] (N_6) to [bend left=10] node [right] {$0$} (N_8);

  \draw [->,DarkBlue] (N_1) to [bend right=10] node [left] {$1$} (N_2);
  \draw [<-] (N_1) to [bend left=10] node [right] {$0$} (N_2);

  \draw [->,DarkBlue] (N_3) to [] node [left] {$1$} (N_6);
  \draw [->,DarkBlue] (N_4) to [] node [left] {$1$} (N_5);
  \draw [<-] (N_3) to [] node [right] {$0$} (N_7);
  \draw [->,DarkBlue] (N_5) to [] node [left] {$1$} (N_7);
  \draw [->,DarkBlue] (N_2) to [] node [left] {$1$} (N_3);
  \draw [<-] (N_2) to [] node [right] {$0$} (N_5);
  \draw [<-] (N_4) to [] node [right] {$0$} (N_3);
\end{tikzpicture} &
      \begin{tikzpicture}[>=latex,join=bevel,scale=.5]
\tiny%
  \node (N_1) at (50bp,306bp) [draw,draw=none] {$1 \!\otimes\!111$};
  \node (N_2) at (5bp,232bp) [draw,draw=none] {$1 \!\otimes\!112$};
  \node (N_3) at (5bp,158bp) [draw,draw=none] {$1 \!\otimes\!122$};
  \node (N_4) at (95bp,232bp) [draw,draw=none] {$2 \!\otimes\!111$};
  \node (N_5) at (95bp,158bp) [draw,draw=none] {$2 \!\otimes\!112$};
  \node (N_6) at (5bp,84bp) [draw,draw=none] {$1 \!\otimes\!222$};
  \node (N_7) at (95bp,84bp) [draw,draw=none] {$2 \!\otimes\!122$};
\node (N_8) at (50bp,10bp) [draw,draw=none] {$2 \!\otimes\!222$};
  \draw [<-] (N_5) to [] node [right] {$0$} (N_6);
  \draw [<-] (N_7) to [] node [right] {$0$} (N_8);
  \draw [->,DarkBlue] (N_6) to [] node [left] {$1$} (N_8);
  \draw [->,DarkBlue] (N_1) to [] node [left] {$1$} (N_2);
  \draw [->,DarkBlue] (N_3) to [] node [left] {$1$} (N_6);
  \draw [->,DarkBlue] (N_4) to [] node [left] {$1$} (N_5);
  \draw [<-] (N_3) to [] node [right] {$0$} (N_7);
  \draw [->,DarkBlue] (N_5) to [] node [left] {$1$} (N_7);
  \draw [->,DarkBlue] (N_2) to [] node [left] {$1$} (N_3);
  \draw [<-] (N_1) to [] node [right] {$0$} (N_4);
  \draw [<-] (N_2) to [] node [right] {$0$} (N_5);
  \draw [<-] (N_4) to [] node [right] {$0$} (N_3);
\end{tikzpicture} &
      \begin{tikzpicture}[>=latex,join=bevel,scale=.5]
\tiny%
  \node (N_1) at (50bp,306bp) [draw,draw=none] {$1 \!\otimes\!111$};

  \node (N_2) at (5bp,232bp) [draw,draw=none] {$1 \!\otimes\!112$};
  \node (N_4) at (95bp,232bp) [draw,draw=none] {$2 \!\otimes\!111$};

  \node (N_3) at (5bp,158bp) [draw,draw=none] {$1 \!\otimes\!122$};
  \node (N_5) at (95bp,158bp) [draw,draw=none] {$2 \!\otimes\!112$};

  \node (N_6) at (5bp,84bp) [draw,draw=none] {$1 \!\otimes\!222$};
  \node (N_7) at (95bp,84bp) [draw,draw=none] {$2 \!\otimes\!122$};

\node (N_8) at (50bp,10bp) [draw,draw=none] {$2 \!\otimes\!222$};
  \draw [<-] (N_7) to [] node [right] {$0$} (N_8);
  \draw [->,DarkBlue] (N_6) to [] node [left] {$1$} (N_8);
  \draw [->,DarkBlue] (N_1) to [] node [left] {$1$} (N_2);

  \draw [->,DarkBlue] (N_3) to [bend right=10] node [left] {$1$} (N_6);
  \draw [<-] (N_3) to [bend left=10] node [right] {$0$} (N_6);

  \draw [->,DarkBlue] (N_4) to [bend right=10] node [left] {$1$} (N_5);
  \draw [<-] (N_4) to [bend left=10] node [right] {$0$} (N_5);

  \draw [->,DarkBlue] (N_5) to [bend right=10] node [left] {$1$} (N_7);
  \draw [<-] (N_5) to [bend left=10] node [right] {$0$} (N_7);

  \draw [->,DarkBlue] (N_2) to [bend right=10] node [left] {$1$} (N_3);
  \draw [<-] (N_2) to [bend left=10] node [right] {$0$} (N_3);

  \draw [<-] (N_1) to [] node [right] {$0$} (N_4);
\end{tikzpicture} \\
      $(aa)$ & $(ab)$ & $(ba)$ & $(bb)$
    \end{tabular}
  \end{bigcenter}
  \caption[Graphes cristallins affines associés au graphe
    cristallin classique $B(1)\otimes B(3)$ de type $A_1$]{Les quatre graphes cristallins affines associés au graphe
    cristallin classique $B(1)\otimes B(3)$ de type $A_1$. Le graphe
    cristallin $B^{1,1}\otimes B^{3,1}$ correspond à $(bb)$. Les
    autres ne proviennent pas de $U_q(\widehat{sl_2})$-modules.}
  \label{fig:promotionsForA11}
\end{figure}

Notre démonstration, passablement technique, repose principalement sur
la combinatoire des tableaux (jeu de taquin, etc.) avec une analyse
fine de la structure du graphe cristallin lorsque le nombre de lignes
est petit et plusieurs inductions pour s'y ramener. L'existence ne
pose pas de problème. La difficulté est de garantir l'unicité.  D'une
part, la définition de bon graphe cristallin affine n'est pas encore
complètement établie dans la communauté; ainsi il faut rajouter des
hypothèses supplémentaires en type $A_1^{(1)}$ (voir
figure~\ref{fig:promotionsForA11}).
D'autre part, il faut imposer une condition de connexité; celle-ci
n'est pas pratique à manipuler car c'est une condition globale, alors
que toutes les autres opérations cristallines sont par nature
locales. Enfin, et c'est ce qui pour le moment nous empêche de
généraliser à trois facteurs ou plus, l'unicité n'est en général vraie
qu'\emph{à isomorphie près}. De fait, on peut présenter cette
conjecture comme un \emph{problème de reconstruction}: chaque
puissance $\pr^k$ de l'opérateur de promotion appliquée au graphe
classique donne une vue partielle du graphe (toutes les arêtes sauf
$k$), à isomorphie près (on ne sait pas quel sommet est envoyé où);
cette collection de vues détermine-t-elle entièrement le graphe
affine à l'isomorphie près? Il reste à évaluer ce qu'apporte ce
point de vue, sachant que les graphes en jeu ont une structure très
particulière.

Une fois n'est pas coutume, je suis rentré dans ce projet \emph{via}
la technique, en aidant à la conception de la bibliothèque sur les
graphes cristallins (voir section~\ref{section.combinat.demo.cristaux}
pour quelques petits exemples de calculs). Puis Jason Bandlow et Anne Schilling ont eu
besoin pour leur exploration informatique de rechercher
systématiquement tous les opérateurs de promotion, et j'ai pris en
charge l'algorithmique et la programmation de cette part. Dans l'état
actuel des connaissances, cette recherche ne peut se faire que par
recherche exhaustive, mais en utilisant des heuristiques de coupe de
branches (séparation et évaluation), avec propagation de contraintes
liées aux propriétés des opérateurs cristallins et exploitation
partielle des symétries.
Par exemple, pour le cas difficile de $B(1) ^{\otimes4}$ en type
$A_2$, où l'espace de recherche est a priori de taille
$144473849856000$, avec $2!3!3!=72$ symétries, l'algorithme explore
$115193$ branches en $5$ heures et $26$ minutes (sur un PC Linux à $2$
GHz ), utilisant $16$Mo de mémoire. Au final, on obtient $8$
opérateurs de promotion connexes isomorphes: $9$ symétries sur les
$72$ ont pu être exploitées pour réduire l'espace de recherche. Encore
une fois, l'explosion combinatoire est importante; le rôle de ce
travail a été de rendre \emph{possible} (éventuellement au prix de
quelques semaines de calcul) l'étude d'exemples non triviaux, sans
lesquels nous n'aurions pas établi la
conjecture~\ref{conjecture.promotion.unique}.  Par la suite, et une
fois rentré dans le sujet, j'ai apporté quelques clefs essentielles de
la démonstration. Étant familier avec les graphes et les problèmes
d'isomorphie, j'ai de plus pu simplifier ou abstraire certains
passages techniques.

\section{Algèbres de Kac et treillis de sous-facteurs}
\label{section.kac}

\subsection{Sous-facteurs}

\newcommand{\II}{\operatorname{II}}

Ma collaboration avec Marie-Claude David a commencé à mon arrivée à
Orsay en 2004, 
grâce aux mots clefs «tours d'algèbres», «diagrammes de Bratelli» et
«calculs dans les algèbres de Hopf». Marie-Claude David s'intéresse
aux inclusions de \emph{facteurs} de type $\II_1$ obtenus via l'action
d'une algèbre de Kac (Hopf-von Neumann) de dimension finie.

Un facteur est une algèbre de von Neumann dont le centre est trivial.
Les facteurs ont été classifiés par Murray et von Neumann (1943),
Connes (1976), et finalement Haagerup. En dimension finie $n$, il
existe à isomorphisme près un unique facteur: $M_n(\CC)$, et un
uniquement plongement de ce facteur dans lui-même.  De même, il existe
un unique facteur $R$ dit hyperfini de type $\II_1$; on peut le
construire comme limite d'une tour d'algèbres $M_n(\CC)$ emboîtées, ou
comme algèbre du groupe $S_\infty$ des permutations de $\NN$ à support
fini. En revanche, ce facteur de type $\II_1$ admet de multiples
plongements dans lui-même.

L'étude des inclusions $N\subset M$ d'indice fini, où $N$ et $M$ sont
des facteurs hyperfinis de type $\II_1$ et $M$ est un module de
dimension finie sur $N$ a connu un essor
considérable~\cite{Jones_Sunder.1997} depuis les articles
fondateurs~\cite{Jones.1983.Subfactors,Jones.1985.KnotInvariant}. Dans
ceux-ci, Vaughan Jones développe les outils fondamentaux (théorie de
l'indice, construction de base, etc.) et en dérive un nouvel
invariant polynomial pour les noeuds (et une médaille Fields). En
particulier, l'inclusion est décrite par une tour d'algèbres de
dimensions finies (la tour dérivée). Souvent, les premières algèbres
(profondeur finie), voire la première (profondeur deux) sont
suffisantes pour caractériser l'inclusion. \emph{Via} le diagramme de
Brattelli de cette tour d'algèbre on peut construire un invariant de
l'inclusion qui est un diagramme de Dynkin.

Le rôle de cette tour d'algèbre est similaire à celui du groupe
d'automorphisme pour une inclusion de corps. Cela donne lieu à une
correspondance de Galois pour le treillis des facteurs intermédiaires
$N \subset P \subset M$. En profondeur deux, le rôle du groupe est
joué par une \emph{algèbre de Kac} de dimension finie (une algèbre de
Hopf involutive semi-simple), ou plus généralement un \emph{groupoïde
  quantique}; celui du treillis des sous-groupes est joué par le
treillis des \emph{sous-algèbres coidéales} (ou pour faire court
\emph{coidéaux})~\cite{Nikshych_Vainerman.2000.2}. À son tour, chaque
coidéal est caractérisé par un élément singulier, son \emph{projecteur
  de Jones}, de sorte que l'on est ramené à étudier un treillis
d'idempotents particuliers de l'algèbre de Kac.

L'apparition des algèbres de Kac dans ce contexte est naturelle. Elles
ont en effet été introduites précisément pour donner un cadre commun
contenant à la fois les algèbres de groupes (qui sont les algèbres de
Kac cocommutatives) et leurs duales (qui sont les algèbres de Kac
commutatives) (voir par exemple~\cite{Enock_Schwartz.1992}), sachant
que dans les deux cas, les coidéaux de l'algèbre correspondent aux
sous-groupes.

\subsection{Études d'exemples d'algèbres de Kac}

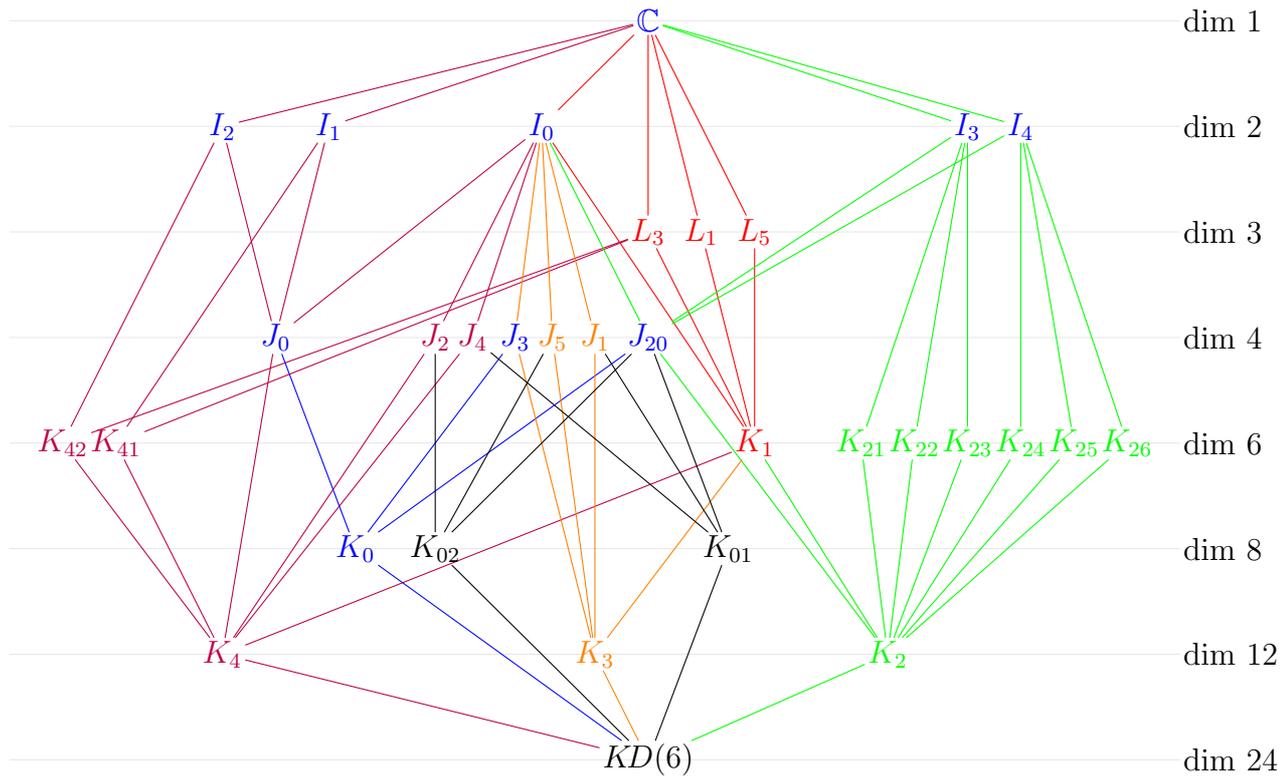
\begin{figure}[h]
  \let\C=\CC
  \begin{bigcenter}
    \pgfdeclarelayer{background}
\pgfdeclarelayer{nodes}
\pgfsetlayers{background,main,nodes}
\begin{tikzpicture}[xscale=.7,yscale=.7]
  \tikzstyle{every node}=[inner sep=1pt]
  \begin{pgfonlayer}{nodes}
    \node(C)	at (6,8)[text=blue] {$\CC$};
    \node(I2)	at (-2,6) [text=blue] {$I_2$};
    \node(I1)	at (0,6) [text=blue] {$I_1$};
    \node(I0)	at (4,6) [text=blue] {$I_0$};
    \node(I3)	at (12,6) [text=blue] {$I_3$};
    \node(I4)	at (13,6) [text=blue] {$I_4$};
    \node(L3)	at (6,4) [text=red] {$L_3$};
    \node(L1)	at (7,4) [text=red] {$L_1$};
    \node(L5)	at (8,4) [text=red] {$L_5$};
    \node(J0)	at (-1,2) [text=blue] {$J_0$};
    \node(J5)	at (2,2) [text=purple] {$J_2$};
    \node(J6)	at (2.7,2) [text=purple] {$J_4$};
    \node(J3)	at (4.2,2) [text=orange] {$J_5$};
    \node(J2)	at (3.5,2) [text=blue] {$J_3$};
    \node(J4)	at (5,2) [text=orange] {$J_1$};
    \node(J1)	at (6,2) [text=blue] {$J_{20}$};
    \node(Kp3)	at (-5,0) [text=purple] {$K_{42}$};
    \node(Kp4)	at (-4,0) [text=purple] {$K_{41}$};
    \node(K1)	at (8,0) [text=red] {$K_1$};
    \node(Kp21)	at (10,0) [text=green] {$K_{21}$};
    \node(Kp22)	at (11,0) [text=green] {$K_{22}$};
    \node(Kp23)	at (12,0) [text=green] {$K_{23}$};
    \node(Kp24)	at (13,0) [text=green] {$K_{24}$};
    \node(Kp25)	at (14,0) [text=green] {$K_{25}$};
    \node(Kp26)	at (15,0) [text=green] {$K_{26}$};
    \node(K0)	at (0.5,-2) [text=blue] {$K_0$};
    \node(Kp135)at (2,-2) {$K_{02}$};
    \node(Kp146)at (7.5,-2) {$K_{01}$};
    \node(K4)	at (-2,-4) [text=purple] {$K_4$};
    \node(K3)	at (5,-4) [text=orange] {$K_3$};
    \node(K2)	at (10.5,-4) [text=green] {$K_2$};
    \node(KD6)	at (6,-6) {$K\!D(6)$};
  \end{pgfonlayer}
  \tikzstyle{every path}=[purple]
  \draw (C) -- (I2);
  \draw (C) -- (I1);
  \draw (I2) -- (J0);
  \draw (I1) -- (J0);
  \draw (I0) -- (J0);
   \draw (I0) -- (J6);
  \draw (I0) -- (J5);
  \draw (L3) -- (Kp3);
  \draw (I2) -- (Kp3);
  \draw (I1) -- (Kp4);
  \draw (K4) -- (Kp3);
  \draw (K4) -- (Kp4);
  \draw (L3) -- (Kp4);
  \draw (J0) -- (K4);
  \draw (J6) -- (K4);
  \draw (J5) -- (K4);
  \draw (K1) -- (K4);
  \draw (K4) -- (KD6);
  \tikzstyle{every path}=[blue]
  \draw (J0) -- (K0);
  \draw (J2) -- (K0);
  \draw (J1) -- (K0);
  \draw (K0) -- (KD6);
  \tikzstyle{every path}=[red]
   \draw (C) -- (I0);
  \draw (C) -- (L3);
  \draw (C) -- (L1);
  \draw (C) -- (L5);
  \d2aw (C) -- (I0);
  \draw (I0) -- (K1);
  \draw (L3) -- (K1);
  \draw (L1) -- (K1);
  \draw (L5) -- (K1);
    \tikzstyle{every path}=[green]
   \draw (C) -- (I3);
  \draw (C) -- (I4);
   \draw (I0) -- (J1);
   \draw (I3) -- (J1);
  \draw (I4) -- (J1);
  \draw (I3) -- (Kp21);
  \draw (I3) -- (Kp22);
  \draw (I3) -- (Kp23);
  \draw (I4) -- (Kp24);
  \draw (I4) -- (Kp25);
  \draw (I4) -- (Kp26);
  \draw (K2) -- (J1);
  \draw (K2) -- (K1);
  \draw (K2) -- (Kp21);
  \draw (K2) -- (Kp22);
  \draw (K2) -- (Kp23);
  \draw (K2) -- (Kp24);
  \draw (K2) -- (Kp25);
  \draw (K2) -- (Kp26);
  \draw (K2) -- (KD6);
  \tikzstyle{every path}=[orange]
  \draw (I0) -- (J2);
  \draw (I0) -- (J3);
  \draw (I0) -- (J4);
  \draw (J4) -- (K3);
  \draw (J3) -- (K3);
  \draw (J2) -- (K3);
  \draw (K1) -- (K3);
  \draw (K3) -- (KD6);
  \tikzstyle{every path}=[black]
  \draw (J5) -- (Kp135);
  \draw (J3) -- (Kp135);
  \draw (J1) -- (Kp135);
  \draw (J6) -- (Kp146);
  \draw (J4) -- (Kp146);
  \draw (J1) -- (Kp146);
  \draw (KD6) -- (Kp135);
  \draw (KD6) -- (Kp146);
  \begin{pgfonlayer}{background}
    \newcommand{\dimline}[2]{\draw[color=black!10,very thin](-6,#1) -- (16,#1); \node[right] at (16,#1){dim $#2$};}
    \dimline{8}{1}
    \dimline{6}{2}
    \dimline{4}{3}
    \dimline{2}{4}
    \dimline{0}{6}
    \dimline{-2}{8}
    \dimline{-4}{12}
    \dimline{-6}{24}
  \end{pgfonlayer}
\end{tikzpicture}
  \end{bigcenter}
  \centering
  \caption{Exemple de treillis de coidéaux: $\KD(6)$}
  \label{figure.treillis.KD6}
\end{figure}

En dehors des treillis des groupes, il existe peu d'exemples
explicités de treillis non triviaux de sous-facteurs.
%
%
\FIXME{motiver plus l'étude d'exemples}
Marie-Claude David avait auparavant établi à la main le treillis des
coidéaux pour la plus petite algèbre de Kac non triviale (de dimension
$8$)~\cite{Kac.Paljutkin.1966}. Nous étudions
ensemble~\cite{David_Thiery.2008.Kac} deux familles d'exemples
construites par Leonid Vainerman par déformation d'algèbres de groupes
de dimension $4n$~\cite{Vainerman.1998} par des $2$-pseudo
cocycles. En dehors des formules pour les opérations de Kac
(coproduit, coinvolution, etc.), rien n'était connu sur ces algèbres,
pour la simple raison que les calculs nécessaires à leur exploration
sont très lourds. L'exploration informatique nous a permis de
conjecturer, puis de démontrer, que ces familles sont en fait
isomorphes pour $n$ pair. Nous avons aussi pu obtenir leurs groupes
d'automorphismes et montrer que, pour $n$ impair, les algèbres de ces
deux familles sont autoduales. Nous avons une description complète du
treillis en petite dimension (en gros $n\leq 7$, voir
figure~\ref{figure.treillis.KD6} pour un exemple en dimension $24$),
pour $n$ premier et, conjecturalement, pour $n$ impair, ainsi que des
résultats partiels pour $n$ pair.
Nous en dérivons, par le calcul des diagrammes de Bratelli, les
graphes principaux d'un certain nombre
d'inclusions.

\subsection{Exploration informatique}

Cette collaboration est un cas exemplaire de ma stratégie favorite.
Marie-Claude David a apporté un sujet intéressant, une culture générale, un
savoir calculer à la main impressionnant et, point de départ
indispensable, un problème concret de calcul. J'ai apporté un outil et
un savoir explorer dans un domaine où, contrairement à la combinatoire
algébrique, l'ordinateur n'est pas un outil traditionnel. J'étais donc
particulièrement attiré par les questions «qu'est-ce qui est
calculable?» «jusqu'où?», «dans quel but?». Comme d'habitude,
l'ordinateur a été un \emph{outil d'exploration} et en particulier un
\emph{outil créatif}, permettant avec les bonnes questions et les bons
algorithmes de deviner ex-nihilo des formules d'isomorphisme comme:
\begin{align*}
  \phi(\lambda(a)) &= a+1/4 (a-a^{-1}) (a^{n}-1) (1 - ba^m)\,,\\
  \phi(\lambda(b)) &=1/2[b(a^n+1)+\mathrm{i}(a^n-1)]a^m\,,
\end{align*}
À titre plus exceptionnel, dans ce domaine où les lemmes tendent à
être techniques et les formules intrinsèquement lourdes, il s'est
révélé être aussi un \emph{outil de preuve} quasiment indispensable,
avec l'élaboration de stratégies pour ramener le cas général à $n$
suffisamment petit pour pouvoir vérifier sur machine.
Même à la machine, les calculs sont difficiles. Ils nécessitent tout
d'abord une réflexion algorithmique. De fait, la construction des
treillis de coïdéaux est encore loin d'être automatisée. Ils
nécessitent aussi une infrastructure d'un haut niveau d'abstraction
(voir les exemples en section~\ref{section.combinat.demo.kac}).
L'effort de développement que j'ai consenti à cette occasion a permis
de compléter considérablement l'algorithmique et l'infrastructure pour
les algèbres de Hopf et la théorie des représentations des algèbres de
dimension finie dans \starcombinat; cela a déjà trouvé d'autres
applications.

\subsection{Perspectives}

L'étape suivante, que nous entreprenons avec Leonid Vainerman et son
étudiante en thèse Camille Mével, est l'étude des exemples plus
généraux venant de \emph{groupoïdes quantiques}, qui correspondent à
des inclusions non irréductibles.  Cela permettra d'obtenir, pour de
petites tailles, tous les treillis provenant d'inclusions d'indice
fini et de profondeur finie; celles-ci peuvent en effet toujours être
réalisées comme inclusions intermédiaires d'inclusions de profondeur
$2$, en général non irréductibles.

De manière générale, les domaines environnants sont pratiquement
vierges d'exploration informatique. Cela ouvre la porte à une
multitude de problèmes intrinsèquement simples, non résolus jusqu'ici
faute de moyens techniques appropriés, et posant des problèmes
d'algorithmique ou de modélisation intéressants. Une mine d'or
potentielle, en particulier pour des étudiants aimant à la fois les
mathématiques et l'informatique. Par exemple: Leonid Vainerman m'a
demandé comment déterminer le groupe d'automorphismes d'une certaine
algèbre obtenue par produit croisé de deux algèbres déformées du
groupe symétrique $\sg$. Pour $n=2$, le résultat est trivial, les
groupes étant commutatifs. Résoudre $n=3$ (dimension $36$) a de bonnes
chances de suggérer la solution pour tout $n$. Pour cela, il faudrait
trouver un algorithme de complexité raisonnable pour calculer tous les
isomorphismes entre deux algèbres de Hopf (ou de Kac) de dimension
finie. De même, il serait urgent de tester ce que les outils de
géométrie algébrique réelle peuvent apporter à l'automatisation de la
construction automatique des projecteurs de Jones.

\section{Polynômes harmoniques pour les opérateurs de Steenrod}
\label{section.steenrod}

Assistant en juin 2000 à l'école «Interactions between Algebraic
Topology and Invariant Theory» à Ioannina, mon attention avait été
captivée lorsque Reg
Wood~\cite{Wood.DOSA.1997,Wood.PSA.1998,Wood.HPSA.2001} présenta une
conjecture que l'on peut réécrire sous la forme:
\begin{conjecture}
  \label{conjecture.harmoniques}
  Le sous-espace des polynômes $p$ de $\QQ[x_1,\dots,x_n]$ satisfaisant
  pour tout $k$ l'équation aux dérivées partielles linéaire:
  \begin{displaymath}
    \left(
      \left(1 + x_1\frac{\partial}{\partial x_1}\right) \frac{\partial}{\partial x_1}^k
      + \dots + 
      \left(1 + x_n\frac{\partial}{\partial x_n}\right) \frac{\partial}{\partial x_n}^k
    \right)
    p = 0
  \end{displaymath}
  est isomorphe à la représentation régulière graduée du groupe
  symétrique. En particulier, il est de dimension $n!$\,.
\end{conjecture}
En effet, comme l'avait mentionné Reg Wood, en oubliant les termes
$x_i\frac{\partial}{\partial x_i}$, ce sous-espace devient celui des
\emph{harmoniques des polynômes symétriques} (une réalisation des
coinvariants du groupe symétrique); dans ce cadre, le résultat est
bien connu, avec une multitude de techniques de démonstration (voir
par exemple~\cite{Garsia_Haiman.OHGR}).

Cette conjecture de Reg Wood n'est pas un phénomène isolé, bien au
contraire. Elle est fortement liée à toute une mouvance de problèmes
où interviennent des polynômes harmoniques, problèmes que j'avais déjà
rencontrés lors de mes travaux sur les algèbres d'invariants de
groupes de
permutations~\cite{Aval_Bergeron.CSH.2003,Aval_all.2004,Haiman.2001}.
En gros, d'où vient-elle? Les algèbres d'invariants de groupes finis
en petite caractéristique apparaissent naturellement comme cohomologie
d'espaces topologiques. Les opérateurs de Steenrod (d'origine
topologique, mais de nature algébrique: ils engendrent une algèbre de
Hopf, l'algèbre de Steenrod) agissent sur les invariants, permettant
d'en construire de nouveaux. Un programme important est de décrire
tous les invariants atteignables de la sorte. Rien que pour le groupe
trivial, c'est-à-dire en agissant simplement sur les polynômes, le
problème est difficile; à l'heure actuelle, il n'y en a même pas de
description conjecturale. La conjecture de Reg Wood est l'analogue de
ce problème en caractéristique $0$. Sa résolution donnerait
certainement des pistes d'attaque pour le cas modulaire.

Étudier cette conjecture était pour moi l'occasion rêvée d'utiliser
sérieusement, et donc de mieux maîtriser, un grand nombre
d'outils. D'abord ceux de la combinatoire algébrique phalanstérienne:
représentations du groupe symétrique et des algèbres de Lie,
déformations non commutatives et interpolations entre algèbres de Hopf
(groupes quantiques), polynômes de Schubert et opérateurs sur les
polynômes (différences divisées avec les algèbres de Hecke affine en
arrière plan). Et ensuite des outils du calcul formel: algorithmique
dans les algèbres de Weyl ou de Ore et bases de Gröbner
semi-commutatives. C'était aussi un bon problème pour commencer à
travailler avec Florent Hivert, ce que nous souhaitions de longue
date.

Afin d'interpoler entre les deux problèmes, nous avons défini la
$q$-algèbre de Steenrod par déformation non commutative de l'algèbre
des polynômes symétriques. Cela revient à introduire un facteur $q$
devant les termes $x_i\frac{\partial}{\partial x_i}$ des équations de
la conjecture~\ref{conjecture.harmoniques}. On appelle $q$-harmoniques
les solutions de ces équations déformées; à $q=0$, ce sont les
harmoniques usuels pour les polynômes symétriques, et à $q=1$ les
polynômes de la conjecture~\ref{conjecture.harmoniques}.  Cela permet
d'utiliser des arguments de spécialisation pour réduire cette
conjecture comme suit.
\begin{proposition}[Hivert, T.~\cite{Hivert_Thiery.SA.2002}]
  Pour démontrer la conjecture~\ref{conjecture.harmoniques}, il suffit
  de montrer que:
  \begin{enumerate}
  \item La dimension des $q$-harmoniques est au moins $n!$.
  \item La dimension des $1$-harmoniques est au plus $n!$.
  \end{enumerate}
\end{proposition}

Jusqu'ici, la conjecture avait été testée complètement jusqu'à
$n=3$. La méthode brutale pour tester la conjecture nécessite de faire
de l'algèbre linéaire sur les polynômes de degré $\binom n 2$, ce qui
fait un espace de dimension $\binom{\binom n
  2+n-1}{n-1}$. L'exploitation des symétries pour décomposer l'espace
des polynômes en petits sous-espaces (composantes de Garnir) nous a
permis de mener une exploration informatique complète jusqu'à $n=5$,
et partielle jusqu'à $n=9$. Cela a confirmé la conjecture pour $q=1$,
et en fait pour tout $q$ sauf pour quelques rationnels de la forme
$-\frac ab$ avec $a$ et $b$ petits.  Mais la suite ne s'est pas
déroulée comme prévu: même si nous avons obtenu quelques résultats
dans des cas particuliers (par exemple en petit degré), la conjecture
a résisté à tous nos efforts, toutes les techniques du cas usuel se
révélant inapplicables~\cite{Hivert_Thiery.SA.2002}.

En 2006, Adriano Garsia (grand spécialiste des polynômes
harmoniques) trouva cette conjecture magnifique lorsque nous la lui
avons exposée. Il s'y est depuis intéressé avec François Bergeron et
Nolan Wallach. Ceux-ci ont obtenu des résultats partiels
supplémentaires, comme une borne supérieure sur la dimension du
sous-espace des solutions.  Ils ont aussi découvert que la fameuse
ex-conjecture de $n!$ concernant les harmoniques pour l'action
diagonale de $\sg$ sur deux jeux de paramètres semblait aussi admettre
une généralisation de type Steenrod. Cela renforce l'idée que la
conjecture de Reg Wood n'est pas un accident isolé.

Nous avons nous-mêmes mis de côté cette conjecture en attendant que
l'étude d'un autre problème nous apporte le bon éclairage. Nous
restons en effet convaincus que cette conjecture n'est pas
intrinsèquement difficile et qu'il suffira de poser la bonne question
à l'ordinateur pour dérouler une preuve élémentaire et constructive.

\TODO{Qsym sur Steenrod?}

\chapter{\starcombinatsimple, boîte à outils pour l'exploration informatique }
\label{chapter.combinat}

Le chapitre qui vient est de nature différente des précédents. En
particulier, il ne contient pas de résultats de recherche à proprement
dit. J'y décris différents aspects de mon travail autour du projet
logiciel \starcombinat que je coordonne depuis sa création en 2000.

La mission de ce projet est de fournir une boîte à outils extensible
pour l'exploration informatique en combinatoire algébrique, et
promouvoir la mutualisation de code entre les chercheurs de ce
domaine.

Mon rôle, secondé par Florent Hivert, est multiple: choix de la
plateforme et du modèle de développement, animation de la communauté,
repérage d'intérêts communs et mise en relation, formation permanente,
assistance aux développeurs (conception, modélisation mathématique,
algorithmique), veille sur la qualité (tests, revues de code),
développement d'une vision globale garante de la cohérence interne,
mise en place de l'architecture logicielle (compilation, documentation,
distribution, tests), travail de fond sur le langage (paradigmes et
idiomes), veille technologique (repérage des techniques et composants
les plus intéressants à intégrer), promotion et valorisation. En bref,
libérer les autres contributeurs des contingences informatiques pour
qu'ils puissent se concentrer sur leur tâche: traduire leur expertise
mathématique en composants logiciels réutilisables à même de résoudre
naturellement et efficacement leurs problèmes de calculs. Enfin, je me
dois d'insuffler une dynamique, en étant en permanence en première
ligne, tout en rentabilisant mes propres développements dans des
projets de recherche.

Après un historique (section~\ref{section.combinat.historique}), et
une démonstration courte (section~\ref{section.combinat.demo}), je
présente quelques problématiques spécifiques à nos besoins en calculs
tant du point de vue du modèle de développement
(section~\ref{section.combinat.devel}) que de la conception
(sections~\ref{section.combinat.objet}
et~\ref{section.combinat.representations}).  

Dans ce domaine, nous sommes des praticiens. Autant que possible, nous
avons essayé de réutiliser des techniques et outils préexistants. Mais
dans certains cas, nous avons été amenés à développer des solutions
originales. Leur conception a été guidée puis validée par notre
expérience sur le terrain, tout en restant au plus près des
mathématiques. Il serait maintenant souhaitable d'avoir une analyse
théorique de notre travail, par exemple dans le cadre d'études de cas
par des doctorants en informatique, pour expliquer pourquoi nos
solutions sont effectivement adaptées aux besoins ou, au contraire,
pour en proposer de meilleures.

Ce chapitre se poursuit avec une étude de cas: comment une
collaboration internationale de chercheurs sur le modèle du libre
permet sur le long terme de bouleverser un outil effectif essentiel:
les objets décomposables ou plus généralement les espèces
combinatoires (section~\ref{section.combinat.especes}). Enfin, on
trouvera en section~\ref{section.combinat.publications} une liste de
publications ayant utilisé ou concernant \mupadcombinat.


\section{\starcombinatsimple, Sharing algebraic combinatorics since 2000}
\label{section.combinat.historique}

\subsection{Pourquoi \mupadcombinat}

Toutes mes recherches ont en commun l'utilisation d'outils
informatiques, et notamment du calcul formel, dans des domaines
propices aux explosions combinatoires. L'exploration informatique sert
de guide, suggérant des conjectures, ou au contraire produisant des
contre-exemples. Elle permet d'étudier des exemples suffisamment
conséquents pour être représentatifs; ces exemples sont le plus
souvent intraitables à la main, et souvent hors de portée des
algorithmes classiques. Cela nécessite de trouver les outils
mathématiques appropriés pour développer de nouveaux algorithmes,
sachant que seule l'efficacité de l'implantation finale décide de la
pertinence de ces outils. Réciproquement, il est indispensable que le
logiciel modélise les objets mathématiques au plus près, afin de
pouvoir exprimer des problèmes complexes dans un langage naturel.

Bien entendu, mener à bien de tels calculs sous-entend un important
travail de programmation, et requiert une large panoplie de
techniques. Ainsi, lors de ma thèse, j'avais été amené à utiliser les
logiciels suivants: \magma (invariants, groupes de permutations), \gap
(groupes, représentations), \texttt{Cocoa}, \texttt{Gb}, \texttt{FGb},
\texttt{Macaulay} (bases de Gröbner efficaces), \texttt{Nauty},
\texttt{Graphlet}, \texttt{Lydia} (théorie des graphes), \texttt{ALP},
\texttt{Linbox} (algèbre linéaire creuse exacte), \maple (calcul
formel) et ses modules (\texttt{combstruct}: objets décomposables,
\texttt{gfun}: séries génératrices, \ace, \SF: fonctions symétriques,
etc.), \mupad (calcul formel) et son module \muec (fonctions
symétriques). Pour aller plus loin j'avais besoin non seulement de les
utiliser simultanément, mais aussi de les combiner aisément les uns
avec les autres.
Une tâche d'\emph{intégration} bien trop ambitieuse pour un seul
individu.



À mon sens, il y avait et il y a toujours un besoin, à l'échelle de
la combinatoire algébrique, pour une boîte à outils de référence,
libre, largement diffusée,
développée par une communauté de chercheurs pour leurs propres besoins
et à l'échelle internationale. En bref, jouant le rôle fédérateur de
\gap pour la théorie des groupes.  Elle doit de plus être basée sur
une plateforme généraliste et un langage de programmation répandu et
de qualité (en particulier orienté objet), réutiliser le maximum de
composants existants, et permettre une synergie avec les domaines
avoisinants (calcul formel, graphes, groupes, etc.).  Il est à noter
qu'une telle boîte à outils, en tout cas dans ces fondamentaux, est
utile pour l'enseignement. Elle a aussi des applications potentielles
en physique mathématique (calculs dans les algèbres de Hopf et
opérades venant des problèmes de renormalisation).
De nombreux logiciels existaient, certains importants et bien établis,
mais aucun ne répondait au cahier des charges.

\subsection{Les débuts de \mupadcombinat}

J'en avais besoin pour mes recherches et je n'étais pas le seul. C'est
un projet qui m'intéressait, et j'en avais peut-être les compétences,
à condition de m'associer aux bons experts. J'ai cherché des
partenaires, établi un cahier des charges détaillé, et longtemps
soupesé les plateformes
disponibles~\cite{Hivert_Thiery.MuPAD-Combinat.2004}. Ce fut la
naissance de \mupadcombinat en décembre 2000, d'abord avec Florent
Hivert et le soutien du Phalanstère de Marne-la-Vallée et de l'équipe
du système de calcul formel \mupad à Paderborn. Puis progressivement
se sont rajoutés François Descouens, Teresa Gomez Diaz,
Jean-Christophe Novelli du Phalanstère, ainsi que Christophe Carré,
Éric Laugerotte, Houda Abbad et Janvier Nzeutchap de Rouen, Frédéric
Chapoton de Lyon; Patrick Lemeur de Montpellier. L'équipe s'est enfin
internationalisée avec Xavier Molinéro de Barcelone, Mike Zabrocki de
Toronto, Anne Schilling et Quiang Wang de Davis.  Sans compter
indirectement tous les auteurs des logiciels que nous avions intégrés
(\texttt{ACE}, $\mu$-\texttt{EC}, \texttt{CS}, \texttt{Symmetrica},
\texttt{Nauty}, etc.).

La collaboration avec Anne Schilling, commencée lors de FPSAC 2006 a
été exemplaire. Ayant une bonne expérience de la programmation et de
l'exploration informatique, elle souhaitait implanter dans
\mupadcombinat une bibliothèque sur les graphes cristallins. De
nombreux courriels suivirent, où d'un côté elle m'expliquait la
théorie et les règles de calculs dans les graphes cristallins, et en
retour je lui donnais des conseils de conception et
d'algorithmique. Trois mois plus tard, elle avait le matériel
nécessaire pour tester ses conjectures sur les graphes cristallins de
Kirillov-Reshetikhin en type $D$ affine~\cite{Schilling.2007}, et fit
le commentaire suivant sur une liste de diffusion à propos de son
expérience avec \mupadcombinat: «For me personally, it has been great,
since it already has a lot of features and, having the support from
Nicolas and others, helped me to write code that I needed in such a
way that it can be reused by others, and I myself could build on it
(which has never been the case before with programs I had written in
Maple or Mathematica). Nicolas definitely has the ability to factor
out main features, which is necessary for such a large scale project.»

\subsection{Maturité et migration vers \sage}

À partir de 2006 plusieurs signes m'ont indiqué que \mupadcombinat
était en train de devenir une (la?) bibliothèque généraliste de
référence dans le domaine. J'ai en effet été conférencier invité aux
ateliers Axiom 2006 et 2007, aux 7ème journées \sage, à la session
spéciale «Applications of Computer Algebra in Enumerative and
Algebraic Combinatorics» de l'AMS Joint Mathematics Meeting à San
Diego. D'autre part, \mupadcombinat avait été choisi comme support
logiciel pour le projet NSF trisannuel «Focused Research Group: Affine
Schubert Calculus» emmené par Anne Schilling, Marc Shimozono et
Jennifer Morse (projet qui a financé mon séjour d'un an en tant que
chercheur à Davis). Enfin, le plus important: \mupadcombinat avait
joué un rôle clef dans de nombreuses publications (voir
section~\ref{section.combinat.publications}).

Dans le même temps nous avons commencé à atteindre les limites de la
plateforme \mupad: la communauté était trop petite, à l'échelle du
domaine, comme à l'échelle du langage.  Du coup, nous étions amenés à
implanter, pour nos besoins propres, des fonctionnalités hors de notre
domaine de compétences (théorie des groupes, communications entre
processus, outils de programmation, etc.). Que \mupad ne soit pas
libre y jouait un rôle important.  Cela n'était pas une surprise:
notre choix de \mupad en 2000 avait été pragmatique; c'était
techniquement la meilleure plateforme pour nos besoins, même si
moralement cela n'était pas satisfaisant. Mais entre temps des
alternatives libres commençaient à émerger, en particulier \axiom et
\sage. Réciproquement celles-ci affichaient un grand intérêt pour
\mupadcombinat. J'ai multiplié les contacts avec leurs communautés
respectives pour préparer le terrain, et pour évaluer l'adaptation
technique pour notre projet. Cela a abouti a la création des
bibliothèques soeurs \aldorcombinat (Ralf Hemmecke et Martin Rubey) et
\sagecombinat (Mike Hansen).

Après de longues discussions nous avons décidé en juin 2008 de
basculer vers \sage. Le coût est important: 100k lignes de code à
traduire\footnote{30k ont déjà été portées par Mike Hansen}; cela nous
absorbera au moins jusqu'en 2009, voire 2010. Mais je suis convaincu
que le jeu en vaut la chandelle. L'opération nous a déjà permis de
joindre nos forces avec Mike Hansen, Jason Bandlow, Franco Saliola,
Greg Musiker, Daniel Bump, Justin Walker, Mark Shimozono, Lenny Tevlin
et Kurt Luoto.

Dix ans d'investissement et d'acharnement commencent à payer. Je vais
enfin avoir à ma disposition la plateforme dont je rêvais, sans
concessions. Pour, en fin de compte, pouvoir faire plus de recherche.

\section{Démonstration courte}
\label{section.combinat.demo}

Dans cette section, nous présentons quelques calculs typiques avec
\mupadcombinat. Outre présenter un aperçu rapide de ses
fonctionnalités, l'objectif est de donner un support concret à
quelques spécificités et concepts qui seront développés par la suite.

\subsection{Combinatoire}

Nous commençons par quelques calculs élémentaires de combinatoire.
\mupadcombinat fournit de nombreuses \emph{classes combinatoires}
prédéfinies (une classe combinatoire est un ensemble sur lequel on
souhaite faire des opérations combinatoires comme compter, énumérer,
tirer au sort, etc). Chacune de ces classes est modélisée par un
domaine (une classe) avec une interface standardisée. Commençons par
lister tous les arbres ordonnés non étiquetés à $5$ sommets:

\def\eng{}
\def\Mup{}
\begin{Mexin}
export(combinat):
\end{Mexin}
\begin{Mexin}
trees::list(5)
\end{Mexin}
\begin{Mexout}
     --   o  ,  o ,  o ,  o  ,  o ,  o ,  o ,   o ,  o ,  o ,  o ,  o ,  o , o --
     |  // \\  /|\  /|\  / \   / \  /|\  / \   / \  / \   |    |    |    |   |  |
     |           |   |     /\    |  |    | |  /\    |    /|\  / \  / \   |   |  |
     |                           |                  |           |  |    / \  |  |
     --                                                                      | --
\end{Mexout}
Nous pouvons aussi juste compter ces arbres, une opération beaucoup
plus rapide\footnote{En tout cas pour $n$ grand; mais alors, on
  n'obtient pas d'information fondamentale sur l'univers.}:
\begin{Mexin}
trees::count(6)
\end{Mexin}
\begin{Mexout}
                               42
\end{Mexout}
Voici un arbre aléatoire. L'affichage en ASCII 2D est loin d'être
parfait, mais la structure de données interne, elle, est robuste.
\begin{Mexin}
trees::random(50)
\end{Mexin}
\begin{Mexout}
                           o
                           |
                         // \  \
                            |/ |            \
                              / \
                                |
                              // \     \
                              |   /// //\\\\\
                                     /\|
                                       |
                                       |
                                      /  \
                                        /  \
                                       / \/|\
                                      /|\   |
                                      ||    |
                                            |
\end{Mexout}

Les algorithmes sous-jacents ne sont pas spécifiques aux arbres. Ils
s'appliquent à toute autre famille d'objets pouvant être définis
récursivement par une grammaire. Voici par exemple la relation de
récurrence pour le nombre d'arbres binaires. Elle est calculée
automatiquement à partir de la grammaire pour ces arbres, et est
exploitée pour faire du comptage efficace:

\begin{Mexin}
r := binaryTrees::grammar::recurrenceRelation():
assume(n>0):
u(n) = factor(op(solve(r, u(n)),1))
\end{Mexin}

\begin{Mexout}
                            2 u(n - 1) (2 n - 1)
                     u(n) = --------------------
                                    n + 1
\end{Mexout}
On reconnaît bien évidement la récurrence usuelle des nombres de
Catalan.

\subsection{Graphes cristallins}
\label{section.combinat.demo.cristaux}

Nous allons maintenant manipuler des graphes cristallins, comme ceux
de la section~\ref{section.cristaux}. Ce sont encore des classes
combinatoires, avec des opérations algébriques supplémentaires. Nous
définissons deux graphes cristallins de Kirillov-Reshetikhin de type
$A_2^{(1)}$:
\begin{Mexin}
C1 := crystals::kirillovReshetikhin(2,2,["A",2,1]):
C2 := crystals::kirillovReshetikhin(1,1,["A",2,1]):
\end{Mexin}
Leurs éléments sont des tableaux semi-standard:
\begin{Mexin}
C1::list()
\end{Mexin}
\begin{Mexout}
     -- +---+---+  +---+---+  +---+---+  +---+---+  +---+---+  +---+---+ --
     |  | 2 | 2 |  | 2 | 3 |  | 2 | 3 |  | 3 | 3 |  | 3 | 3 |  | 3 | 3 |  |
     |  +---+---+  +---+---+  +---+---+  +---+---+  +---+---+  +---+---+  |
     |  | 1 | 1 |, | 1 | 1 |, | 1 | 2 |, | 1 | 1 |, | 1 | 2 |, | 2 | 2 |  |
     -- +---+---+  +---+---+  +---+---+  +---+---+  +---+---+  +---+---+ --
\end{Mexout}
sur lesquels agissent des opérateurs montants $e_i$ et descendants
$f_i$:
\begin{Mexin}
x := C1::list()[3]
\end{Mexin}
\begin{Mexout}
                              +---+---+
                              | 2 | 3 |
                              +---+---+
                              | 1 | 2 |
                              +---+---+
\end{Mexout}
\begin{Mexin}
x::e(0), x::e(1), x::e(2),   x::f(0), x::f(1), x::f(2)
\end{Mexin}
\begin{Mexout}
        +---+---+  +---+---+        +---+---+        +---+---+
        | 3 | 3 |  | 2 | 3 |        | 2 | 2 |        | 3 | 3 |
        +---+---+  +---+---+        +---+---+        +---+---+
        | 2 | 2 |, | 1 | 1 |, FAIL, | 1 | 1 |, FAIL, | 1 | 2 |
        +---+---+  +---+---+        +---+---+        +---+---+
\end{Mexout}

Après avoir manipulé leurs éléments, nous faisons des calculs sur les
graphes cristallins eux même. Par exemple, nous vérifions que le
graphe $C_1$ a un unique automorphisme et que les graphes $C_1$ et
$C_2$ ne sont pas isomorphes:
\begin{Mexin}
C1::isomorphisms(C1)
\end{Mexin}
\begin{Mexout}
                         [proc g(x) ... end]
\end{Mexout}
\begin{Mexin}
C1::isomorphisms(C2)
\end{Mexin}
\begin{Mexout}
                                  []
\end{Mexout}
Comme son nom l'indique, la méthode \Mup{C1::isomorphisms} renvoie
la liste de tous les isomorphismes, sous forme de fonctions que l'on
pourrait appliquer aux éléments de $C_1$.

\begin{figure}[h]
  \centering
  \scalebox{0.5}{
    \begin{tikzpicture}[>=latex,join=bevel,]
        \node (N_1) at (178bp,532bp) [draw,draw=none]
        {${\def\lr#1#2#3{\multicolumn{1}{#1@{\hspace{.6ex}}c@{\hspace{.6ex}}#2}{\raisebox{-.3ex}{$#3$}}}\raisebox{-.6ex}{$\begin{array}[b]{cc}\cline{1-1}\cline{2-2}\lr{|}{|}{2}
                & \lr{|}{|}{2}\\\cline{1-1}\cline{2-2}\lr{|}{|}{1} &
                \lr{|}{|}{1}\\\cline{1-1}\cline{2-2}\end{array}$}}
          \otimes
          {\def\lr#1#2#3{\multicolumn{1}{#1@{\hspace{.6ex}}c@{\hspace{.6ex}}#2}{\raisebox{-.3ex}{$#3$}}}\raisebox{-.6ex}{$\begin{array}[b]{c}\cline{1-1}\lr{|}{|}{1}\\\cline{1-1}\end{array}$}}$};
        \node (N_2) at (178bp,446bp) [draw,draw=none]
        {${\def\lr#1#2#3{\multicolumn{1}{#1@{\hspace{.6ex}}c@{\hspace{.6ex}}#2}{\raisebox{-.3ex}{$#3$}}}\raisebox{-.6ex}{$\begin{array}[b]{cc}\cline{1-1}\cline{2-2}\lr{|}{|}{2}
                & \lr{|}{|}{2}\\\cline{1-1}\cline{2-2}\lr{|}{|}{1} &
                \lr{|}{|}{1}\\\cline{1-1}\cline{2-2}\end{array}$}}
          \otimes
          {\def\lr#1#2#3{\multicolumn{1}{#1@{\hspace{.6ex}}c@{\hspace{.6ex}}#2}{\raisebox{-.3ex}{$#3$}}}\raisebox{-.6ex}{$\begin{array}[b]{c}\cline{1-1}\lr{|}{|}{2}\\\cline{1-1}\end{array}$}}$};
        \node (N_3) at (255bp,446bp) [draw,draw=none]
        {${\def\lr#1#2#3{\multicolumn{1}{#1@{\hspace{.6ex}}c@{\hspace{.6ex}}#2}{\raisebox{-.3ex}{$#3$}}}\raisebox{-.6ex}{$\begin{array}[b]{cc}\cline{1-1}\cline{2-2}\lr{|}{|}{2}
                & \lr{|}{|}{3}\\\cline{1-1}\cline{2-2}\lr{|}{|}{1} &
                \lr{|}{|}{1}\\\cline{1-1}\cline{2-2}\end{array}$}}
          \otimes
          {\def\lr#1#2#3{\multicolumn{1}{#1@{\hspace{.6ex}}c@{\hspace{.6ex}}#2}{\raisebox{-.3ex}{$#3$}}}\raisebox{-.6ex}{$\begin{array}[b]{c}\cline{1-1}\lr{|}{|}{1}\\\cline{1-1}\end{array}$}}$};
        \node (N_4) at (152bp,360bp) [draw,draw=none]
        {${\def\lr#1#2#3{\multicolumn{1}{#1@{\hspace{.6ex}}c@{\hspace{.6ex}}#2}{\raisebox{-.3ex}{$#3$}}}\raisebox{-.6ex}{$\begin{array}[b]{cc}\cline{1-1}\cline{2-2}\lr{|}{|}{2}
                & \lr{|}{|}{2}\\\cline{1-1}\cline{2-2}\lr{|}{|}{1} &
                \lr{|}{|}{1}\\\cline{1-1}\cline{2-2}\end{array}$}}
          \otimes
          {\def\lr#1#2#3{\multicolumn{1}{#1@{\hspace{.6ex}}c@{\hspace{.6ex}}#2}{\raisebox{-.3ex}{$#3$}}}\raisebox{-.6ex}{$\begin{array}[b]{c}\cline{1-1}\lr{|}{|}{3}\\\cline{1-1}\end{array}$}}$};
        \node (N_5) at (274bp,360bp) [draw,draw=none]
        {${\def\lr#1#2#3{\multicolumn{1}{#1@{\hspace{.6ex}}c@{\hspace{.6ex}}#2}{\raisebox{-.3ex}{$#3$}}}\raisebox{-.6ex}{$\begin{array}[b]{cc}\cline{1-1}\cline{2-2}\lr{|}{|}{2}
                & \lr{|}{|}{3}\\\cline{1-1}\cline{2-2}\lr{|}{|}{1} &
                \lr{|}{|}{1}\\\cline{1-1}\cline{2-2}\end{array}$}}
          \otimes
          {\def\lr#1#2#3{\multicolumn{1}{#1@{\hspace{.6ex}}c@{\hspace{.6ex}}#2}{\raisebox{-.3ex}{$#3$}}}\raisebox{-.6ex}{$\begin{array}[b]{c}\cline{1-1}\lr{|}{|}{2}\\\cline{1-1}\end{array}$}}$};
        \node (N_6) at (342bp,360bp) [draw,draw=none]
        {${\def\lr#1#2#3{\multicolumn{1}{#1@{\hspace{.6ex}}c@{\hspace{.6ex}}#2}{\raisebox{-.3ex}{$#3$}}}\raisebox{-.6ex}{$\begin{array}[b]{cc}\cline{1-1}\cline{2-2}\lr{|}{|}{3}
                & \lr{|}{|}{3}\\\cline{1-1}\cline{2-2}\lr{|}{|}{1} &
                \lr{|}{|}{1}\\\cline{1-1}\cline{2-2}\end{array}$}}
          \otimes
          {\def\lr#1#2#3{\multicolumn{1}{#1@{\hspace{.6ex}}c@{\hspace{.6ex}}#2}{\raisebox{-.3ex}{$#3$}}}\raisebox{-.6ex}{$\begin{array}[b]{c}\cline{1-1}\lr{|}{|}{1}\\\cline{1-1}\end{array}$}}$};
        \node (N_7) at (185bp,274bp) [draw,draw=none]
        {${\def\lr#1#2#3{\multicolumn{1}{#1@{\hspace{.6ex}}c@{\hspace{.6ex}}#2}{\raisebox{-.3ex}{$#3$}}}\raisebox{-.6ex}{$\begin{array}[b]{cc}\cline{1-1}\cline{2-2}\lr{|}{|}{2}
                & \lr{|}{|}{3}\\\cline{1-1}\cline{2-2}\lr{|}{|}{1} &
                \lr{|}{|}{1}\\\cline{1-1}\cline{2-2}\end{array}$}}
          \otimes
          {\def\lr#1#2#3{\multicolumn{1}{#1@{\hspace{.6ex}}c@{\hspace{.6ex}}#2}{\raisebox{-.3ex}{$#3$}}}\raisebox{-.6ex}{$\begin{array}[b]{c}\cline{1-1}\lr{|}{|}{3}\\\cline{1-1}\end{array}$}}$};
        \node (N_8) at (253bp,274bp) [draw,draw=none]
        {${\def\lr#1#2#3{\multicolumn{1}{#1@{\hspace{.6ex}}c@{\hspace{.6ex}}#2}{\raisebox{-.3ex}{$#3$}}}\raisebox{-.6ex}{$\begin{array}[b]{cc}\cline{1-1}\cline{2-2}\lr{|}{|}{2}
                & \lr{|}{|}{3}\\\cline{1-1}\cline{2-2}\lr{|}{|}{1} &
                \lr{|}{|}{2}\\\cline{1-1}\cline{2-2}\end{array}$}}
          \otimes
          {\def\lr#1#2#3{\multicolumn{1}{#1@{\hspace{.6ex}}c@{\hspace{.6ex}}#2}{\raisebox{-.3ex}{$#3$}}}\raisebox{-.6ex}{$\begin{array}[b]{c}\cline{1-1}\lr{|}{|}{2}\\\cline{1-1}\end{array}$}}$};
        \node (N_9) at (429bp,274bp) [draw,draw=none]
        {${\def\lr#1#2#3{\multicolumn{1}{#1@{\hspace{.6ex}}c@{\hspace{.6ex}}#2}{\raisebox{-.3ex}{$#3$}}}\raisebox{-.6ex}{$\begin{array}[b]{cc}\cline{1-1}\cline{2-2}\lr{|}{|}{3}
                & \lr{|}{|}{3}\\\cline{1-1}\cline{2-2}\lr{|}{|}{1} &
                \lr{|}{|}{1}\\\cline{1-1}\cline{2-2}\end{array}$}}
          \otimes
          {\def\lr#1#2#3{\multicolumn{1}{#1@{\hspace{.6ex}}c@{\hspace{.6ex}}#2}{\raisebox{-.3ex}{$#3$}}}\raisebox{-.6ex}{$\begin{array}[b]{c}\cline{1-1}\lr{|}{|}{2}\\\cline{1-1}\end{array}$}}$};
        \node (N_10) at (78bp,446bp) [draw,draw=none]
        {${\def\lr#1#2#3{\multicolumn{1}{#1@{\hspace{.6ex}}c@{\hspace{.6ex}}#2}{\raisebox{-.3ex}{$#3$}}}\raisebox{-.6ex}{$\begin{array}[b]{cc}\cline{1-1}\cline{2-2}\lr{|}{|}{2}
                & \lr{|}{|}{3}\\\cline{1-1}\cline{2-2}\lr{|}{|}{1} &
                \lr{|}{|}{2}\\\cline{1-1}\cline{2-2}\end{array}$}}
          \otimes
          {\def\lr#1#2#3{\multicolumn{1}{#1@{\hspace{.6ex}}c@{\hspace{.6ex}}#2}{\raisebox{-.3ex}{$#3$}}}\raisebox{-.6ex}{$\begin{array}[b]{c}\cline{1-1}\lr{|}{|}{1}\\\cline{1-1}\end{array}$}}$};
        \node (N_11) at (78bp,360bp) [draw,draw=none]
        {${\def\lr#1#2#3{\multicolumn{1}{#1@{\hspace{.6ex}}c@{\hspace{.6ex}}#2}{\raisebox{-.3ex}{$#3$}}}\raisebox{-.6ex}{$\begin{array}[b]{cc}\cline{1-1}\cline{2-2}\lr{|}{|}{3}
                & \lr{|}{|}{3}\\\cline{1-1}\cline{2-2}\lr{|}{|}{1} &
                \lr{|}{|}{2}\\\cline{1-1}\cline{2-2}\end{array}$}}
          \otimes
          {\def\lr#1#2#3{\multicolumn{1}{#1@{\hspace{.6ex}}c@{\hspace{.6ex}}#2}{\raisebox{-.3ex}{$#3$}}}\raisebox{-.6ex}{$\begin{array}[b]{c}\cline{1-1}\lr{|}{|}{1}\\\cline{1-1}\end{array}$}}$};
        \node (N_12) at (171bp,188bp) [draw,draw=none]
        {${\def\lr#1#2#3{\multicolumn{1}{#1@{\hspace{.6ex}}c@{\hspace{.6ex}}#2}{\raisebox{-.3ex}{$#3$}}}\raisebox{-.6ex}{$\begin{array}[b]{cc}\cline{1-1}\cline{2-2}\lr{|}{|}{2}
                & \lr{|}{|}{3}\\\cline{1-1}\cline{2-2}\lr{|}{|}{1} &
                \lr{|}{|}{2}\\\cline{1-1}\cline{2-2}\end{array}$}}
          \otimes
          {\def\lr#1#2#3{\multicolumn{1}{#1@{\hspace{.6ex}}c@{\hspace{.6ex}}#2}{\raisebox{-.3ex}{$#3$}}}\raisebox{-.6ex}{$\begin{array}[b]{c}\cline{1-1}\lr{|}{|}{3}\\\cline{1-1}\end{array}$}}$};
        \node (N_13) at (253bp,188bp) [draw,draw=none]
        {${\def\lr#1#2#3{\multicolumn{1}{#1@{\hspace{.6ex}}c@{\hspace{.6ex}}#2}{\raisebox{-.3ex}{$#3$}}}\raisebox{-.6ex}{$\begin{array}[b]{cc}\cline{1-1}\cline{2-2}\lr{|}{|}{3}
                & \lr{|}{|}{3}\\\cline{1-1}\cline{2-2}\lr{|}{|}{1} &
                \lr{|}{|}{1}\\\cline{1-1}\cline{2-2}\end{array}$}}
          \otimes
          {\def\lr#1#2#3{\multicolumn{1}{#1@{\hspace{.6ex}}c@{\hspace{.6ex}}#2}{\raisebox{-.3ex}{$#3$}}}\raisebox{-.6ex}{$\begin{array}[b]{c}\cline{1-1}\lr{|}{|}{3}\\\cline{1-1}\end{array}$}}$};
        \node (N_14) at (400bp,188bp) [draw,draw=none]
        {${\def\lr#1#2#3{\multicolumn{1}{#1@{\hspace{.6ex}}c@{\hspace{.6ex}}#2}{\raisebox{-.3ex}{$#3$}}}\raisebox{-.6ex}{$\begin{array}[b]{cc}\cline{1-1}\cline{2-2}\lr{|}{|}{3}
                & \lr{|}{|}{3}\\\cline{1-1}\cline{2-2}\lr{|}{|}{1} &
                \lr{|}{|}{2}\\\cline{1-1}\cline{2-2}\end{array}$}}
          \otimes
          {\def\lr#1#2#3{\multicolumn{1}{#1@{\hspace{.6ex}}c@{\hspace{.6ex}}#2}{\raisebox{-.3ex}{$#3$}}}\raisebox{-.6ex}{$\begin{array}[b]{c}\cline{1-1}\lr{|}{|}{2}\\\cline{1-1}\end{array}$}}$};
        \node (N_15) at (25bp,274bp) [draw,draw=none]
        {${\def\lr#1#2#3{\multicolumn{1}{#1@{\hspace{.6ex}}c@{\hspace{.6ex}}#2}{\raisebox{-.3ex}{$#3$}}}\raisebox{-.6ex}{$\begin{array}[b]{cc}\cline{1-1}\cline{2-2}\lr{|}{|}{3}
                & \lr{|}{|}{3}\\\cline{1-1}\cline{2-2}\lr{|}{|}{2} &
                \lr{|}{|}{2}\\\cline{1-1}\cline{2-2}\end{array}$}}
          \otimes
          {\def\lr#1#2#3{\multicolumn{1}{#1@{\hspace{.6ex}}c@{\hspace{.6ex}}#2}{\raisebox{-.3ex}{$#3$}}}\raisebox{-.6ex}{$\begin{array}[b]{c}\cline{1-1}\lr{|}{|}{1}\\\cline{1-1}\end{array}$}}$};
        \node (N_16) at (171bp,102bp) [draw,draw=none]
        {${\def\lr#1#2#3{\multicolumn{1}{#1@{\hspace{.6ex}}c@{\hspace{.6ex}}#2}{\raisebox{-.3ex}{$#3$}}}\raisebox{-.6ex}{$\begin{array}[b]{cc}\cline{1-1}\cline{2-2}\lr{|}{|}{3}
                & \lr{|}{|}{3}\\\cline{1-1}\cline{2-2}\lr{|}{|}{1} &
                \lr{|}{|}{2}\\\cline{1-1}\cline{2-2}\end{array}$}}
          \otimes
          {\def\lr#1#2#3{\multicolumn{1}{#1@{\hspace{.6ex}}c@{\hspace{.6ex}}#2}{\raisebox{-.3ex}{$#3$}}}\raisebox{-.6ex}{$\begin{array}[b]{c}\cline{1-1}\lr{|}{|}{3}\\\cline{1-1}\end{array}$}}$};
        \node (N_17) at (339bp,102bp) [draw,draw=none]
        {${\def\lr#1#2#3{\multicolumn{1}{#1@{\hspace{.6ex}}c@{\hspace{.6ex}}#2}{\raisebox{-.3ex}{$#3$}}}\raisebox{-.6ex}{$\begin{array}[b]{cc}\cline{1-1}\cline{2-2}\lr{|}{|}{3}
                & \lr{|}{|}{3}\\\cline{1-1}\cline{2-2}\lr{|}{|}{2} &
                \lr{|}{|}{2}\\\cline{1-1}\cline{2-2}\end{array}$}}
          \otimes
          {\def\lr#1#2#3{\multicolumn{1}{#1@{\hspace{.6ex}}c@{\hspace{.6ex}}#2}{\raisebox{-.3ex}{$#3$}}}\raisebox{-.6ex}{$\begin{array}[b]{c}\cline{1-1}\lr{|}{|}{2}\\\cline{1-1}\end{array}$}}$};
        \node (N_18) at (171bp,16bp) [draw,draw=none]
        {${\def\lr#1#2#3{\multicolumn{1}{#1@{\hspace{.6ex}}c@{\hspace{.6ex}}#2}{\raisebox{-.3ex}{$#3$}}}\raisebox{-.6ex}{$\begin{array}[b]{cc}\cline{1-1}\cline{2-2}\lr{|}{|}{3}
                & \lr{|}{|}{3}\\\cline{1-1}\cline{2-2}\lr{|}{|}{2} &
                \lr{|}{|}{2}\\\cline{1-1}\cline{2-2}\end{array}$}}
          \otimes
          {\def\lr#1#2#3{\multicolumn{1}{#1@{\hspace{.6ex}}c@{\hspace{.6ex}}#2}{\raisebox{-.3ex}{$#3$}}}\raisebox{-.6ex}{$\begin{array}[b]{c}\cline{1-1}\lr{|}{|}{3}\\\cline{1-1}\end{array}$}}$};
     \draw [->] (N_1) ..controls (178bp,503bp) and (178bp,486bp) .. (N_2);
     \pgfsetstrokecolor{black}
     \draw (187bp,489bp) node {$1$};
     \draw [->] (N_1) ..controls (204bp,503bp) and
        (221bp,484bp) .. (N_3);
     \draw (232bp,489bp) node {$2$};
     \draw [->] (N_2) ..controls (170bp,417bp) and (165bp,400bp)
        .. (N_4);
     \draw (176bp,403bp) node {$2$};
     \draw [->] (N_3)
        ..controls (261bp,417bp) and (265bp,400bp) .. (N_5);
     \draw (274bp,403bp) node {$1$};
     \draw [->] (N_3) ..controls
        (285bp,416bp) and (303bp,398bp) .. (N_6);
     \draw (314bp,403bp) node {$2$};
     \draw [->] (N_4) ..controls (163bp,331bp) and
        (170bp,314bp) .. (N_7);
     \draw (179bp,317bp) node {$2$};
     \draw [->] (N_5) ..controls (267bp,331bp) and (263bp,314bp) .. (N_8);
     \draw (274bp,317bp) node {$1$};
     \draw [->] (N_5)
        ..controls (302bp,346bp) and (305bp,345bp) .. (308bp,344bp)
        .. controls (338bp,333bp) and (351bp,342bp) .. (379bp,326bp)
        .. controls (388bp,320bp) and (387bp,315bp) .. (395bp,308bp)
        .. controls (398bp,304bp) and (402bp,301bp) .. (N_9);
     \draw (404bp,317bp) node {$2$};
     \draw [<-] (N_1) ..controls
        (135bp,494bp) and (112bp,475bp) .. (N_10);
     \draw (148bp,489bp) node {$0$};
     \draw [->] (N_10) ..controls (78bp,417bp) and (78bp,400bp) .. (N_11);
     \draw (87bp,403bp) node {$2$};
     \draw [->] (N_6) ..controls (384bp,351bp) and (407bp,342bp)
        .. (420bp,326bp) .. controls (426bp,319bp) and (428bp,309bp)
        .. (N_9);
     \draw (435bp,317bp) node {$1$};
     \draw [->] (N_7)
        ..controls (176bp,253bp) and (173bp,246bp) .. (172bp,240bp)
        .. controls (170bp,232bp) and (169bp,222bp) .. (N_12);
     \draw (181bp,231bp) node {$1$};
     \draw [->] (N_7) ..controls
        (218bp,253bp) and (225bp,247bp) .. (231bp,240bp) .. controls
        (237bp,232bp) and (242bp,222bp) .. (N_13);
     \draw (248bp,231bp) node {$2$};
     \draw [<-] (N_2) ..controls (197bp,400bp) and
        (210bp,370bp) .. (222bp,344bp) .. controls (230bp,325bp) and
        (239bp,304bp) .. (N_8);
     \draw (231bp,360bp) node {$0$};
     \draw [->] (N_8) ..controls (220bp,253bp) and (212bp,247bp)
        .. (206bp,240bp) .. controls (198bp,232bp) and (191bp,221bp)
        .. (N_12);
     \draw (215bp,231bp) node {$2$};
     \draw [->] (N_9)
        ..controls (419bp,245bp) and (413bp,228bp) .. (N_14);
     \draw (426bp,231bp) node {$1$};
     \draw [<-] (N_3) ..controls
        (218bp,432bp) and (215bp,431bp) .. (212bp,430bp) .. controls
        (180bp,420bp) and (168bp,427bp) .. (139bp,412bp) .. controls
        (121bp,402bp) and (104bp,387bp) .. (N_11);
     \draw (151bp,403bp) node {$0$};
     \draw [->] (N_11) ..controls (56bp,339bp) and
        (50bp,332bp) .. (46bp,326bp) .. controls (41bp,318bp) and
        (37bp,308bp) .. (N_15);
     \draw (55bp,317bp) node {$1$};
     \draw [<-] (N_4) ..controls (133bp,314bp) and (126bp,284bp)
        .. (133bp,258bp) .. controls (138bp,237bp) and (150bp,217bp)
        .. (N_12);
     \draw (142bp,274bp) node {$0$};
     \draw [->] (N_12)
        ..controls (171bp,159bp) and (171bp,142bp) .. (N_16);
     \draw (180bp,145bp) node {$2$};
     \draw [<-] (N_6) ..controls
        (310bp,299bp) and (277bp,234bp) .. (N_13);
     \draw (314bp,274bp) node {$0$};
     \draw [->] (N_13) ..controls (225bp,158bp) and
        (207bp,140bp) .. (N_16);
     \draw (230bp,145bp) node {$1$};
     \draw [<-] (N_5) ..controls (302bp,323bp) and (316bp,306bp)
        .. (327bp,290bp) .. controls (349bp,260bp) and (375bp,224bp)
        .. (N_14);
     \draw (357bp,274bp) node {$0$};
     \draw [->] (N_14)
        ..controls (380bp,159bp) and (367bp,141bp) .. (N_17);
     \draw (385bp,145bp) node {$1$};
     \draw [<-] (N_10) ..controls
        (44bp,410bp) and (32bp,393bp) .. (26bp,376bp) .. controls
        (16bp,347bp) and (19bp,310bp) .. (N_15);
     \draw (38bp,360bp) node {$0$};
     \draw [<-] (N_11) ..controls (90bp,306bp) and
        (101bp,260bp) .. (115bp,222bp) .. controls (123bp,199bp) and
        (126bp,194bp) .. (137bp,172bp) .. controls (145bp,153bp) and
        (156bp,132bp) .. (N_16);
     \draw (124bp,231bp) node {$0$};
     \draw [->] (N_16) ..controls (171bp,73bp) and (171bp,56bp) .. (N_18);
     \draw (180bp,59bp) node {$1$};
     \draw [<-] (N_8)..controls (283bp,213bp) and (316bp,148bp) .. (N_17);
     \draw (311bp,188bp) node {$0$};
     \draw [->] (N_17) ..controls (285bp,74bp) and (238bp,50bp) .. (N_18);
     \draw (283bp,59bp) node {$2$};
     \draw [<-] (N_15) ..controls (69bp,196bp) and (138bp,74bp) .. (N_18);
     \draw (109bp,145bp) node {$0$};
   \end{tikzpicture}}
  \caption{Le graphe cristallin affine $B^{2,2} \otimes B^{1,1}$ en type $A_2^{(1)}$}
\label{figure.cristal.C1C2}
\end{figure}
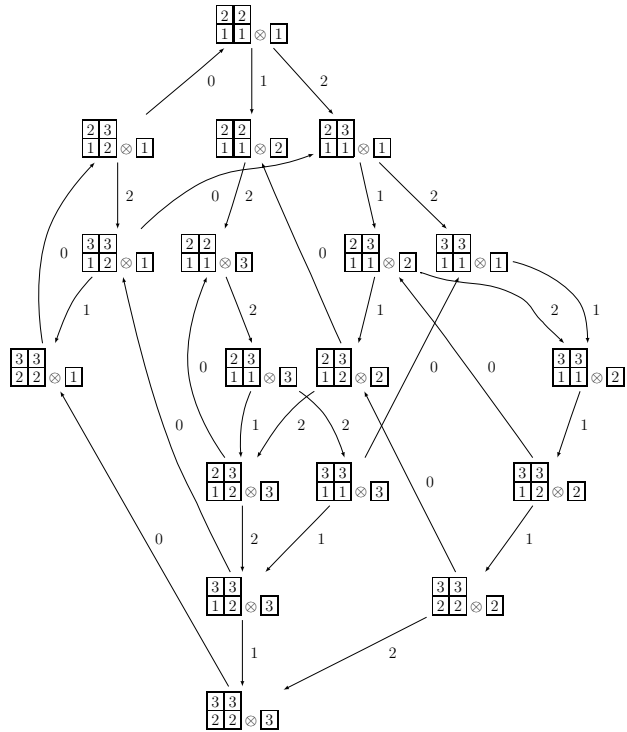

Pour conclure, nous construisons le produit tensoriel des graphes
cristallins $C_1$ et $C_2$ (en fait leur produit cartésien, ce qui
correspond au produit tensoriel des modules dont ils indexent les
bases), et en demandons une représentation graphique (en \LaTeX). Ce
type de représentation est bien entendu un outil important d'exploration.
\begin{Mexin}
operators::setTensorSymbol("#"):
C := C1 # C2:
viewTeX(C::TeXClass())
\end{Mexin}
(voir figure~\ref{figure.cristal.C1C2} pour le résultat)

Cet exemple reste un jouet. De par leur construction, les graphes
cristallins sont de taille explosive. Afin de permettre des calculs
poussés, l'implantation utilise systématiquement l'évaluation
paresseuse. Ainsi, seule la partie du graphe effectivement
étudiée est dépliée en mémoire.

Les calculs précédents illustrent une autre spécificité importante de
nos besoins. Appelons \emph{domaine} un ensemble muni d'opérations
(par exemple l'anneau des matrices $M_n(\RR)$, l'anneau de polynômes
$\CC[x,y,z]$, les fonctions de $\CC$ dans $\CC$). En calcul numérique
ou formel, les manipulations concernent principalement les objets de
ces domaines (pivot de Gauss, base de Gröbner, intégration). En
revanche, en combinatoire algébrique (mais aussi en théorie des
groupes), la richesse des calculs vient le plus souvent de la
manipulation \emph{simultanée} des éléments des domaines (les tableaux
avec l'action des opérateurs cristallins) et des domaines eux-mêmes
(les graphes cristallins). Cela sera discuté dans la
section~\ref{section.combinat.domaines}


\subsection{Algèbre de mélange sur les arbres}

Nous avons vu que les classes combinatoires ont une interface
systématique. Cela est bien entendu souhaitable pour l'utilisateur.
Mais cela permet surtout de les utiliser comme briques logicielles
pour des constructions plus avancées. Nous illustrons cela en donnant
une implantation de l'algèbre graduée dont la base est indexée par les
arbres et le produit est induit par le produit de mélange sur les
arbres (c'est la réalisation de l'algèbre de Loday-Ronco dans la base
$p$):
\begingroup
\microtypesetup{kerning=false}
\begin{listing}{1}
domain ShuffleAlgebraOnTrees(R = Dom::ExpressionField(): Cat::Ring)
    category Cat::GradedAlgebraWithBasis(R);
    inherits Dom::FreeModule(trees, R);

    oneBasis := trees::zero;

    mult2Basis :=
    proc(t1: trees, t2: trees)
        local t;
    begin
        dom::plus(dom::term(t) $ t in trees::shuffle(t1,t2));
    end_proc;
end_domain
\end{listing}
\endgroup
La ligne 1 précise que la constuction de l'algèbre prend un paramètre
$R$, l'anneau de base, avec une valeur par défaut, le corps des
expressions. La ligne 2 fait une promesse au système: l'objet implanté
est une algèbre graduée. La ligne 3 précise la structure de données
pour les éléments: des vecteurs sur $R$ indexés par les arbres.  Il
reste à tenir les promesses, en spécifiant l'unité (ligne 5) et la
règle de calcul du produit, ici exprimée sur la base. La syntaxe
nécessite un apprentissage, mais au final tout le code est signifiant
mathématiquement.

\subsection{Quelques calculs dans les algèbres de dimension finie}
\label{section.combinat.demo.kac}

Nous poursuivons en présentant quelques calculs typiques sur des
algèbres. Ici, nous travaillons sur l'algèbre de Kac $\KD(n)$ obtenue
par déformation de l'algèbre du groupe dihédral $D_{2n}$ étudiée en
section~\ref{section.kac} du
chapitre~\ref{chapter.theorieDesRepresentations}. Il n'est pas
nécessaire de connaître les détails, l'objectif des exemples étant
uniquement d'illustrer le niveau d'abstraction typique des calculs que
nous souhaitons mener.

L'algèbre $\KD(n)$ est pour l'instant implantée dans une feuille de
travail séparée que l'on charge ici:
\begin{Mexin}
read("experimental/2005-09-08-David.mu"):
\end{Mexin}
\Mup{KD(3)} modélise l'algèbre abstraite; après l'avoir construite
\begin{Mexin}
KD3 := KD(3):
\end{Mexin}
\begin{Mexout}
\end{Mexout}
nous demandons quelles sont ses propriétés connues:
\begin{Mexin}
KD3::categories
\end{Mexin}
\begin{Mexout}
     [Cat::HopfAlgebraWithSeveralBases(Q(II, epsilon)),
      TwistedDihedralOrQuaternionGroupAlgebra(3),
      Cat::AlgebraWithSeveralBases(Q(II, epsilon)),
      Cat::Algebra(Q(II, epsilon)), 
      Cat::ModuleWithSeveralBases(Q(II, epsilon)), 
      Cat::Ring, Cat::Module(Q(II, epsilon)),
      Cat::DomainWithSeveralRepresentations, Cat::Rng, Cat::SemiRing,
      Cat::LeftModule(KD(3, Q(II, epsilon))),
      Cat::LeftModule(Q(II, epsilon)), Cat::RightModule(Q(II, epsilon)),
      Cat::UseOverloading, Cat::FacadeDomain, Cat::SemiRng, Cat::Monoid,
      Cat::AbelianGroup, Cat::SemiGroup, Cat::CancellationAbelianMonoid,
      Cat::AbelianMonoid, Cat::AbelianSemiGroup, Cat::Object,
      Cat::BaseCategory]
\end{Mexout}
À quelques exceptions près, ce sont des informations mathématiques.

Pour faire des calculs, nous avons besoin d'une représentation
concrète de cette algèbre. Celle utilisant la base du groupe est
modélisée par \Mup{KD3::group}. Nous commençons par introduire des
notations courtes pour ses générateurs:
\begin{Mexin}
[aa,bb] := KD3::group::algebraGenerators::list()
\end{Mexin}
\begin{Mexout}
                         [B(a), B(b)]
\end{Mexout}
Les quelques calculs suivants, dans la base du groupe, montrent que le
produit n'est pas déformé:
\begin{Mexin}
bb^2
\end{Mexin}
\begin{Mexout}
                             B(1)
\end{Mexout}
\begin{Mexin}
aa^2, aa^6, bb*aa
\end{Mexin}
\begin{Mexout}
                        2            5
                     B(a ), B(1), B(a  b)
\end{Mexout}
\begin{Mexin}
(1 - aa^3)*(bb + aa^3) + 1/2*bb*aa^3
\end{Mexin}
\begin{Mexout}
                                 3
                          1/2 B(a  b)
\end{Mexout}
De part la théorie des représentations, l'algèbre $\KD(3)$ admet une
autre représentation concrète comme algèbre de matrices par blocs,
modélisées par \Mup{KD3::matrix}. Nous utilisons ici l'isomorphisme
entre les deux représentations (transformée de Fourier, d'où la
présence de \Mup{epsilon} qui représente une racine $\epsilon$ de l'unité).
\begin{Mexin}
KD3::M(aa + 2*bb)
\end{Mexin}
\begin{Mexout}
     +-                                                           -+
     |  3,  0,  0, 0,    0,         0,           0,          0     |
     |                                                             |
     |  0, -1,  0, 0,    0,         0,           0,          0     |
     |                                                             |
     |  0,  0, -3, 0,    0,         0,           0,          0     |
     |                                                             |
     |  0,  0,  0, 1,    0,         0,           0,          0     |
     |                                                             |
     |  0,  0,  0, 0, epsilon,      2,           0,          0     |
     |                                                             |
     |  0,  0,  0, 0,    2,    1 - epsilon,      0,          0     |
     |                                                             |
     |  0,  0,  0, 0,    0,         0,      epsilon - 1,     2     |
     |                                                             |
     |  0,  0,  0, 0,    0,         0,           2,      -epsilon  |
     +-                                                           -+
\end{Mexout}
Cet isomorphisme n'a été explicitement implanté que dans un sens (en
donnant l'image de $a$ et $b$); l'isomorphisme inverse est construit
automatiquement par inversion de matrices.

Le coproduit est déformé par un cocycle du sous-groupe $\langle
a^3,b\rangle$. Pour les éléments de ce sous-groupe, le coproduit n'est
donc pas déformé (\# dénote le produit tensoriel $\otimes$):
\begin{Mexin}
coproduct(aa^3), coproduct(bb)
\end{Mexin}
\begin{Mexout}
                     3       3
                  B(a ) # B(a ), B(b) # B(b)
\end{Mexout}
Par contre, le coproduit de $a$ est compliqué:
\begin{Mexin}
coproduct(aa)
\end{Mexin}
\begin{Mexout}
             4         4      /   II        \    4         5    / II        \    4         2     
     1/16 B(a  b) # B(a  b) + | - -- - 1/16 | B(a  b) # B(a ) + | -- - 1/16 | B(a  b) # B(a  b) +
                              \    8        /                   \  8        /                    

                        ... 100 lignes coupées ...

                        2
        -1/16 B(a) # B(a ) + 7/16 B(a) # B(a)
\end{Mexout}

Pour illustrer des calculs tensoriels typiques, nous vérifions que
l'antipode est correcte, c'est-à-dire que la formule $\mu \circ (\id
\otimes S) \circ \Delta$ redonne bien la counité de l'algèbre de Hopf:
\begin{Mexin}
K := KD3::G:   // un simple raccourci
checkAntipode := K::mu @ ( K::id # K::antipode ) @ K::coproduct:
\end{Mexin}
\begin{Mexin}
checkAntipode(x) $ x in K::basis::list()
\end{Mexin}
\begin{Mexout}
     B(1), B(1), B(1), B(1), B(1), B(1), B(1), B(1), B(1), B(1), B(1), B(1)
\end{Mexout}

Notre étude concernait les coidéaux de $\KD(3)$, c'est à dire des
sous-algèbres «stables à droite» pour le coproduit. Notons $e$ les
unités matricielles de l'algèbre.
\begin{Mexin}
e := KD3::e:
\end{Mexin}
Nous calculons maintenant une base du coidéal $K_2:=I(e_1+e_2)$ engendré
par $e_1+e_2$:
\begin{Mexin}
K2basis := coidealAndAlgebraClosure([ e(1)+e(2) ])
\end{Mexin}
\begin{Mexout}
     -- +-                        -+  +-                        -+       +-                        -+ --
     |  |  1, 0, 0, 0, 0, 0, 0, 0  |  |  0, 0, 0, 0, 0, 0, 0, 0  |       |  0, 0, 0, 0, 0, 0, 0, 0  |  |
     |  |                          |  |                          |       |                          |  |
     |  |  0, 1, 0, 0, 0, 0, 0, 0  |  |  0, 0, 0, 0, 0, 0, 0, 0  |       |  0, 0, 0, 0, 0, 0, 0, 0  |  |
     |  |                          |  |                          |       |                          |  |
     |  |  0, 0, 0, 0, 0, 0, 0, 0  |  |  0, 0, 1, 0, 0, 0, 0, 0  |       |  0, 0, 0, 0, 0, 0, 0, 0  |  |
     |  |                          |  |                          |       |                          |  |
     |  |  0, 0, 0, 0, 0, 0, 0, 0  |  |  0, 0, 0, 1, 0, 0, 0, 0  |       |  0, 0, 0, 0, 0, 0, 0, 0  |  |
     |  |                          |, |                          |, ..., |                          |  |
     |  |  0, 0, 0, 0, 0, 0, 0, 0  |  |  0, 0, 0, 0, 0, 0, 0, 0  |       |  0, 0, 0, 0, 0, 0, 0, 0  |  |
     |  |                          |  |                          |       |                          |  |
     |  |  0, 0, 0, 0, 0, 0, 0, 0  |  |  0, 0, 0, 0, 0, 0, 0, 0  |       |  0, 0, 0, 0, 0, 0, 0, 0  |  |
     |  |                          |  |                          |       |                          |  |
     |  |  0, 0, 0, 0, 0, 0, 0, 0  |  |  0, 0, 0, 0, 0, 0, 0, 0  |       |  0, 0, 0, 0, 0, 0, 0, 0  |  |
     |  |                          |  |                          |       |                          |  |
     |  |  0, 0, 0, 0, 0, 0, 0, 0  |  |  0, 0, 0, 0, 0, 0, 0, 0  |       |  0, 0, 0, 0, 0, 0, 0, 1  |  |
     -- +-                        -+  +-                        -+       +-                        -+ --
\end{Mexout}
On reconnaît $e_1+e_2$ dans le premier élément. 

Manipuler directement une telle base n'est pas commode pour étudier
les propriétés algébriques de $K_2$. Nous construisons donc un objet
qui modélise cette sous-algèbre, en promettant à \mupad qu'il s'agit
d'une sous-algèbre de Hopf (ce que nous savions par ailleurs).
\begin{Mexin}
K2 := Dom::SubFreeModule(K2basis,
        [Cat::FiniteDimensionalHopfAlgebraWithBasis(KD3::coeffRing)]):
\end{Mexin}
Nous pouvons maintenant demander si cette sous-algèbre est commutative
ou cocommutative:
\begin{Mexin}
K2::isCommutative(), K2::isCocommutative()
\end{Mexin}
\begin{Mexout}
                          TRUE, FALSE
\end{Mexout}
Nous en déduisons que c'est forcément l'algèbre duale d'une algèbre de
groupe que nous souhaitons retrouver. Pour cela, nous construisons
l'algèbre duale de $K_2$:
\begin{Mexin}
K2dual := K2::Dual():
\end{Mexin}
L'algorithmique sous-jacente est triviale: les opérateurs sont définis
par simple transposition de ceux de l'algèbre originale. Mais il est
pratique de ne pas avoir à s'en soucier.

Calculons les éléments de type groupe de cette algèbre (ici,
l'algorithmique est non triviale!):
\begin{Mexin}
K2dual::groupLikeElements()
\end{Mexin}
\begin{Mexout}
     [B([1, 1]), B([7, 7]), B([3, 3]), B([8, 8]), B([5, 5]) + -II B([6, 5]), B([5, 5]) + II B([6, 5])]
\end{Mexout}
Ces éléments forment un groupe, le groupe intrinsèque. Pour le
déterminer, il reste à reconnaître la règle de produit.
\begin{Mexin}
G := K2dual::intrinsicGroup():
G::list()
\end{Mexin}
\begin{Mexout}
           [[], [1], [1, 1], [2], [1, 2], [1, 1, 2]]
\end{Mexout}
  $G$ modélise le groupe; ses éléments on été exprimés en fonction de
  générateurs choisis au hasard. C'est brutal, \mupad n'ayant que des
  fonctionnalités élémentaires pour les groupes. Cela est cependant
  suffisant ici: on reconnaît à nouveau un groupe dihédral.  Nous
  finissons en calculant à titre de vérification quelques informations
  sur la théorie des représentations:
\begin{Mexin}
K2dual::isSemiSimple()
\end{Mexin}
\begin{Mexout}
                             TRUE
\end{Mexout}
\begin{Mexin}
K2dual::simpleModulesDimensions()
\end{Mexin}
\begin{Mexout}
                           [2, 1, 1]
\end{Mexout}

Résumons ce que nous avons vu. Pour explorer rapidement des structures
algébriques (par exemple une algèbre de dimension finie), nous avons
besoin de poser, si possible en quelques minutes, des questions d'un
relativement haut niveau d'abstraction (l'algèbre est elle
semi-simple?).

Cela nécessite de modéliser simultanément l'algèbre \emph{et} ses
éléments, en s'appuyant sur une représentation concrète; celle-ci
provient usuellement d'une construction basée sur des modèles
combinatoires. Souvent, l'algorithmique sous-jacente est alors une
simple agglomération d'éléments simples (Euclide pour l'arithmétique
dans une extension algébrique, quelques produits s'appuyant sur des
règles combinatoires, un pivot de Gauss). Il serait cependant trop
fastidieux (si ce n'est impossible en pratique) de devoir s'y ramener
explicitement à chaque fois.

Nous avons donc besoin d'un système souple et expressif, permettant de
poser naturellement des questions sur des constructions avancées,
obtenues par composition de briques de base. Vue la multitude de
modèles combinatoires intéressants, les briques et les constructions
doivent être les plus génériques possible, sans cependant entraver la
complexité algorithmique. À ce titre, la conception joue un rôle tout
aussi important que l'algorithmique. Enfin, il devrait être facile de
faire appel à des outils optimisés d'autres domaines mathématiques
(théorie des groupes, algèbre linéaire creuse, etc).

\section{Modèle de développement}
\label{section.combinat.devel}

\subsection{Quelques mots d'ordre de \starcombinatsimple}

Afin d'éclairer la philosophie de \starcombinat, nous donnons, en les
expliquant, quelques-uns de ses mots d'ordre.

\newcommand{\motto}[1]{\textbf{#1.}}

\motto{Toute ligne de code doit être écrite en vue d'une
  application immédiate à un projet de recherche, tout en ayant une
  vision à long terme}

L'objectif est de garantir que le code est utile, qu'il a été testé,
que la conception et la fonctionnalité ont été validées par la
pratique, et que le travail fourni sera valorisé scientifiquement par
des publications. Ce dernier point est essentiel pour maintenir la
motivation, et éviter des frustrations ou des frictions quant à la
paternité du code.

Dans le même temps, il faut viser un code générique et réutilisable,
afin que les investissement soient rentabilisés sur le long terme.

\motto{Écrit par une communauté informelle de chercheurs, pour les
  chercheurs}

Ce point est avant tout un corollaire du point précédent. De plus, des
relations informelles, dans une certaine convivialité, encouragent
l'entraide, les échanges, et favorisent l'émergence d'un corpus
d'expertise qui bénéficie à tous.

\motto{Sous licence libre}

Trop de logiciels sont morts suite à une volonté politique de leurs
institutions de tenter une commercialisation inadéquate avec fermeture
du code. Toute application industrielle de \starcombinat est la
bienvenue! Notre longue et fructueuse collaboration avec SciFace est
la meilleure preuve de notre pragmatisme. Il serait par exemple
parfaitement acceptable de commercialiser des services d'expertise
autour de \starcombinat si le besoin s'en faisait sentir. Mais
l'objectif premier est de mutualiser les efforts de dévelopement des
chercheurs pour la recherche. Le contrôle du logiciel doit rester
entre les mains des chercheurs. Pour cela la meilleure garantie, outre
l'utilisation d'une licence libre type GPL, est de multiplier les
institutions concernées et les sources de financements pour qu'aucune
ne soit dominante.

Le coût de développement de \starcombinat est essentiellement
humain. Que chaque chercheur y participe sur la base du volontariat,
et à la hauteur de ses besoins propres, justifie en soi
l'investissement de son institution. Cela permet à \starcombinat de
fonctionner avec un budget direct très limité.

\motto{Fédérateur à l'échelle internationale}

Trop de logiciels sont restés confidentiels à cause de rivalités entre
communautés. Il est important d'éviter qu'une communauté ne se
l'approprie entraînant un rejet par les autres. Pour cela, il faut
franchir rapidement les frontières (géographiques ou thématiques), par
la base, c'est-à-dire en convaincant un par un les chercheurs.

\motto{Géré par des chercheurs permanents}

Trop de doctorants ont été sacrifiés, la charge de développement ne
leur permettant pas de publier. Réciproquement, trop de logiciels ont
été jetés après le départ de leur responsable ex-doctorant. Le travail
de fond sur une bibliothèque de cette ampleur peut exiger des
investissements massifs rentables uniquement sur le long
terme. Maintenir dans le même temps une productivité scientifique
régulière demande de l'expérience. En revanche, les doctorants sont
bien placés pour intervenir en périphérie, en développant les outils
dont ils ont besoin pour leur propre recherche. Cela doit se faire
sous la supervision amicale d'un permanent dont le rôle est de
garantir une bonne conception et l'intégration à long terme.

\motto{Programmation extrême}

En passant outre la dénomination pompeuse, les démarches de la
programmation extrême s'appliquent fort bien à notre type de
développement: refactorisation permanente, revues de code,
programmation guidée par les tests, intégration constante. Le
développement doit en effet s'adapter aux besoins très variables des
projets de recherche dont les progrès sont par nature
imprévisibles. Seule réserve: de par la dispersion géographique des
chercheurs, la programmation par paires ne peut avoir lieu que
ponctuellement (ateliers, etc.).

\motto{Rendre trivial en pratique ce qui est trivial en théorie}

L'expérience montre que la plupart des calculs repose à 90\% sur une
myriade de petits détails triviaux mathématiquement (calcul de la
longueur d'une permutation; extension par linéarité d'un opérateur,
construction du dual d'une algèbre, etc.). Un rôle essentiel de la
boîte à outils est de prendre en charge la majorité d'entre eux pour
que le programmeur puisse se concentrer sur le cœur de son
problème. C'est à ce niveau là que la mutualisation fonctionne le
mieux, soutenue par une modélisation de haut niveau au plus près du
langage mathématique.

\motto{Algorithmique et conception sont les mamelles de l'efficacité}

La perception de l'efficacité dépend fortement des besoins. Comme nous
l'avons vu, le plus souvent en exploration informatique nous avons à
contrôler une explosion combinatoire. De ce fait, une optimisation
n'est intéressante que si elle permet de calculer vraiment plus loin,
ce qui nécessite un gain en complexité algorithmique.
D'autre part, du point de vue du chercheur, la mesure perçue du temps
inclut non seulement le temps du calcul lui-même, mais aussi le temps
d'exprimer la question.
Enfin, nous devons traiter une grande variété de problèmes, chaque
nouvelle question apportant son lot de spécificités. Notre situation
est donc très différente d'une bibliothèque comme \texttt{GMP} ou
\texttt{Linbox} dont le rôle est de traiter un nombre limité de
problèmes (une dizaine pour \texttt{Linbox}: calcul de rang, de
déterminant, etc), mais ce de manière extrêmement optimisée.

Notre priorité est donc à une bonne complexité, à la souplesse et à
l'expressivité. De ce fait, c'est l'algorithmique et la conception qui
priment.

Notre chance est que, la plupart du temps, notre code a principalement
un rôle de pilotage à partir de briques de base dont la granularité
est relativement élevée (arithmétique sur des entiers longs, algèbre
linéaire, etc.). De ce fait, les surcoûts éventuels liés à
l'utilisation d'un langage souple et de haut niveau (appels de
méthodes, résolution de surcharge, code interprété, etc.)  sont
négligeables par rapport aux calculs sous-jacents. Ainsi, il est
possible de gagner sur les deux tableaux, à condition que la
plateforme permette simultanément de programmer à un haut niveau et
d'utiliser des briques de base externes hautement optimisées.

Les opérations combinatoires élémentaires (par exemple le calcul du
nombre de descentes d'une permutation) sont une exception notable; ici
la granularité est petite (accès aux éléments d'une liste et
manipulations de petits entiers). Cela peut être géré de manière
satisfaisante si la plateforme permet de plus de compiler les sections
critiques en fournissant éventuellement des annotations supplémentaire
(par ex. \sage + \texttt{Cython}).


\subsection{Distribution et outils de développement}

Rendre \mupadcombinat facile à installer et à utiliser par le plus
grand nombre était un point essentiel pour développer notre
communauté. La distribution de \mupadcombinat m'a demandé un travail
considérable:
\begin{itemize}
\item développement d'une architecture de compilation aux standards
  GNU (automake, autoconf, etc.) pour les modules dynamiques (C/C++),
  pour la documentation (\LaTeX puis XML avec adaptation de scripts de
  l'équipe \mupad en \java, \perl et \ruby) et pour les tests;
\item compilation pour toutes les plateformes: Linux, MacOS X, Windows
  (via cygwin sur une machine virtuelle QEMU), voire Zaurus; cela a
  demandé un certain nombre de correctifs ou d'adaptations sur des
  logiciels que nous intégrions (Symmetrica, Nauty, lrcalc);
\item création de paquetages rpm et deb pour Linux, d'un installateur
  NSIS pour Windows, d'un CDROM vif basé sur Morphix;
\item codéveloppement de
  \href{http://mupacs.sf.net/}{\texttt{MuPACS}}, un mode \mupad pour
  \texttt{emacs}, avec dévermineur et aide intégrés;
\item mise en place d'outils collaboratifs sur sourceforge (site Web,
  listes de diffusions, serveur CVS puis subversion, Wiki, etc.).
\end{itemize}

Ce travail d'ingénierie, très gourmand en temps, était difficile à
déléguer faute de compétences. De ce fait, en 2007 comme en 2008, la
sortie stable traditionnelle pour FPSAC n'a pas eu lieu, malgré un
grand nombre de fonctionnalités nouvelles. Fort heureusement,
l'intégration dans le projet \sage nous permet maintenant de
mutualiser ces aspects techniques; ils sont maintenant entièrement
gérés par un ingénieur merveilleusement compétent et motivé (Michael
Abshoff, merci mille fois!).

Le seul point restant à gérer est du ressort du modèle de
développement: le choix des outils collaboratifs et leur bonne
utilisation. L'intégration dans \sage a été l'occasion d'une évolution
pour pallier quelques dérives. J'illustre ici l'une d'entre elles.

La plupart du temps, l'implantation des outils nécessaires à un calcul
donné est transverse, impliquant quelques aspects combinatoires et
algébriques qui sont du ressort de \starcombinat, et d'autres touchant
plus au cœur de la plateforme (correctifs de bogues, algèbre
linéaire, polynômes, etc.). L'organisation de \mupadcombinat comme
\emph{bibliothèque} pour \mupad ne permet pas d'isoler ces
modifications: nous avions dû dupliquer un grand nombre de fichiers de
la bibliothèque de \mupad pour y intégrer un mélange de modifications,
certaines indubitables, d'autres plus expérimentales.  Dans la
pratique la gestion de la réintégration ultérieure de ces
modifications dans la version officielle de \mupad est manuelle,
fastidieuse, et source d'erreurs. De fait, \mupadcombinat et \mupad
étaient progressivement en train de diverger.

Le modèle de développement de \sage résout ce problème en s'organisant
autour des patchs\footnote{Pièces?}: chacun de ceux-ci isole une modification
relativement atomique(ajout d'une fonctionnalité, correctif d'un bogue), mais qui peut concerner plusieurs
fichiers. Dans ce cadre, \sagecombinat devient avant tout une
sous-communauté de \sage partageant une collection de patchs
expérimentaux qui sont développés et utilisés en commun jusqu'à être
suffisamment matures pour être intégrés dans \sage. Par exemple, le
patch jetant les bases des systèmes de racines a bénéficié des
contributions (et donc de l'expertise et du point de vue) de cinq
chercheurs sur une durée de deux mois avant d'être intégré dans \sage;
entre-temps, son utilisation par d'autres membres, en particulier dans
le cadre d'un autre patch sur les graphes cristallins, avait validé sa
conception.

\begin{figure}[h]
  \centering
  \includegraphics[width=\textwidth]{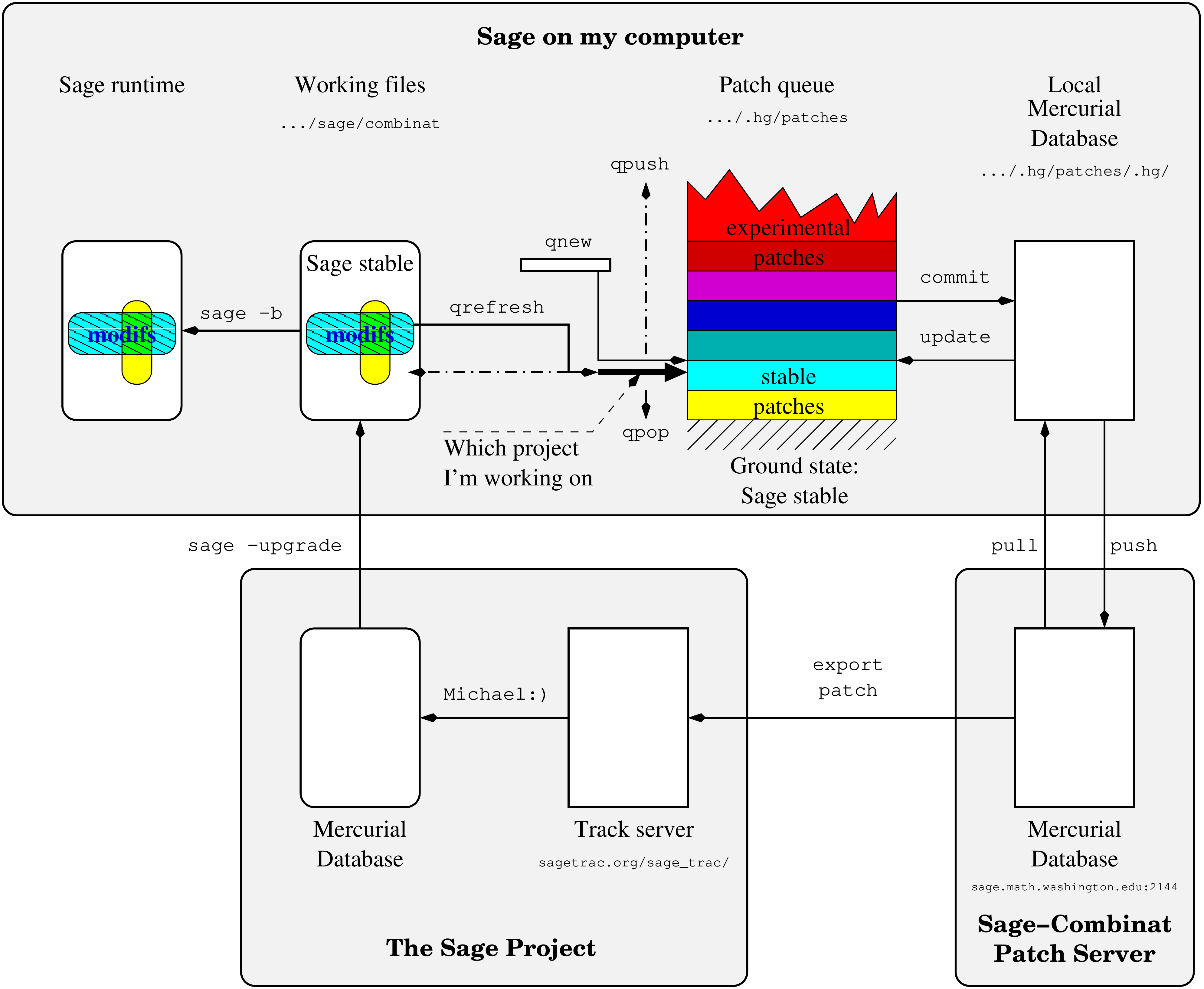}
  \caption[Workflow pour le développement de \sagecombinat \emph{via} un
    serveur de patchs]{Workflow pour le développement de \sagecombinat \emph{via} un
    serveur de patchs; voir
    \url{http://wiki.sagemath.org/combinat/Mercurial} pour les
    explications}
  \label{figure.patchserver}
\end{figure}
\sagecombinat est la première sous-communauté de \sage ayant besoin de
s'organiser de la sorte. Le choix avec Jason Bandlow, Mike Hansen et
Anne Schilling du bon outil collaboratif (les
piles\footnote{littéralement des files (\emph{patch queue}), mais dans
la pratique, nous les utilisons plutôt comme piles} de patchs du système de gestion de
version décentralisé \texttt{Mercurial}) et d'un bon
workflow\footnote{Quel nom affreux! Flot de développements?}
(voir figure~\ref{figure.patchserver}) a nécessité plusieurs semaines
d'expérimentation et surtout la grande patience de nombreux cobayes,
afin de concilier des contraintes inédites:

\begin{itemize}
\item Encourager une intégration rapide dans \sage pour s'adapter à
  son rythme de développement rapide (une nouvelle version par mois).
\item S'adapter à la dispersion géographique et thématique des
  chercheurs et à des rythmes de travail très irréguliers; un patch
  peut rester dormant quelques mois (enseignement, habilitation);
  d'autres peuvent être développés en trois jours par quatre personnes
  à l'occasion d'une conférence.
\item Encourager la collaboration autour de chaque patch.
\item Contrôler la complexité et garantir la robustesse, pour que des
  utilisateurs ou contributeurs occasionnels puissent en faire une
  utilisation simple.
\end{itemize}

\section{Conception: modélisation objet}
\label{section.combinat.objet}
\subsection{Importance de la programmation orientée objet}

\begin{quotation}
  «La programmation orientée objet (POO) ou programmation par objet,
  est un paradigme de programmation informatique qui consiste en la
  définition et l'assemblage de briques logicielles appelées objets;
  un objet représente un concept, une idée ou toute entité du monde
  physique, comme une voiture, une personne ou encore une page d'un
  livre.»\\
  \href{http://fr.wikipedia.org/wiki/Orienté_objet}{\cite[Orienté objet]{Wikipedia_FR}}
\end{quotation}

Je rajoute, dans notre contexte, tout objet mathématique sur lequel on
souhaite faire un calcul: un nombre, un vecteur, un arbre, une
partition; mais aussi, j'y reviendrai en
section~\ref{section.combinat.domaines}, un groupe, une classe
combinatoire, une algèbre.

Pour un informaticien, l'utilisation de la programmation orientée
objet (apparue dès les années 1970!) semble aller de soi pour un
projet logiciel de l'ampleur d'un système de calcul formel. Certains
de ces systèmes, comme \axiom, l'ont d'ailleurs intégrée très tôt.
Singulièrement, son utilité reste encore contestée, certains mettant
en valeur que des systèmes comme \maple (ou \mathematica), au succès
d'ailleurs incontestable jusqu'ici, en particulier en combinatoire,
semblent fonctionner très bien sans. À cela, deux facteurs me
semble-t-il. Le premier est qu'ils ont été conçus avec très peu de
structures de données: nombres (entiers, flottants, etc.), listes,
tables, matrices, et surtout l'universelle expression.  Ces structures
sont implantées dans le noyau, sans possibilité d'extension. À
quelques inconvénients près (opérateur \&* pour la multiplication de
matrices), cela est effectivement suffisant pour l'essentiel du calcul
formel universitaire de licence: algèbre linéaire, calculus, etc. Au
prix de contorsions, beaucoup de problèmes peuvent s'y ramener. On
peut aussi songer à la gamme d'applications de \matlab, autour d'une
seule structure de données (en caricaturant à peine): la matrice à
coefficients flottants. L'autre facteur est que l'essentiel du code
mathématique de \maple, incluant les multiples contributions externes
est composé de petites bibliothèques indépendantes.

Mais la programmation orientée objet devient fondamentale à notre
échelle:
\begin{itemize}
\item \mupadcombinat contient plus de 600 briques logicielles (600
  domaines, ou classes concrètes, pour 5000 méthodes et 130k lignes de
  code), que l'on peut composer entre elles. Cela permet, par exemple,
  de construire une algèbre dont la base est indexée par des arbres et
  les coefficients sont dans une extension algébrique $\QQ[\epsilon]$ avec
  $\epsilon$ une racine de l'unité (voir la démonstration en
  section~\ref{section.combinat.demo}).  En \maple, calculer avec des
  matrices à coefficients dans $\ZZ/2\ZZ$ est déjà malaisé.
\item \mupadcombinat étend considérablement la hiérarchie de
  catégories (classes abstraites) de \mupad (voir
  figure~\ref{figure.categories}). Cela a permis de factoriser une
  grande quantité de code, en le rendant générique, alors que, par
  exemple, 90\% du code de \ace était dupliqué entre les bibliothèques
  de fonctions symétriques commutatives et non commutatives. 
\end{itemize}

\begin{sidewaysfigure}
  \begin{bigcenter}
    \fbox{\resizebox{15cm}{!}{\input{Fig/CategoriesMuPAD}}}
  \end{bigcenter}
  \begin{bigcenter}
    \fbox{\resizebox{25cm}{!}{\input{Fig/Categories}}}
  \end{bigcenter}
  \caption[Hiérarchie des catégories (classes abstraites) dans
    \mupad et \mupadcombinat]{Hiérarchie des catégories (classes abstraites) dans
    \mupad (en haut) et \mupadcombinat (en bas)}
  \label{figure.categories}
\end{sidewaysfigure}

Les mathématiques sont au cœur de la discussion. D'une part, la
hiérarchie de catégories \emph{modélise} naturellement la hiérarchie
usuelle des catégories mathématiques (d'où le nom). D'autre part le
choix de la catégorie où est implantée une fonctionnalité donnée est
conditionné par le niveau d'abstraction auquel l'algorithme utilisé
s'applique. Par exemple, l'algorithme de calcul du radical, que nous
avons utilisé à l'origine pour des algèbres de type Hecke, a été
implanté au niveau des algèbres de dimension finie. Plus tard, il a
été remonté au niveau des algèbres non nécessairement associatives (un
simple déplacement du code), ce qui permet de l'appliquer aux algèbres
de Lie.

\subsection{Améliorations de l'héritage multiple en \mupad}

La modélisation de la hiérarchie usuelle des catégories mathématiques
requiert un mécanisme d'héritage multiple. Par exemple, une algèbre
est à la fois un anneau et un module sur le corps de base. En \mupad,
ce mécanisme se situe naturellement au niveau des catégories (classes
abstraites). L'ordre dans lequel on recherche une méthode dans la
hiérarchie est notoirement sensible. D'un point de vue mathématique,
on souhaite que cet ordre soit au moins une extension linéaire de la
hiérarchie; on veut en effet garantir qu'une implantation spécialisée
d'une méthode dans une sous-catégorie A surcharge systématiquement
celle générique définie dans une sur-catégorie B.  La construction de
notre hiérarchie de catégories a révélé que ce n'était pas le cas dans
\mupad. L'ordre de recherche suivait en effet un parcours en largeur
(comme jusqu'à très récemment en \texttt{Python}), mais leur hiérarchie
originale était suffisamment simple et proche d'être graduée pour que
cela n'ait jamais posé de problème.

C'est un exemple typique de modification, plutôt simple mais profonde
et donc sensible, que j'ai apportée à \mupad, avec la bénédiction de ses
développeurs. De même, j'ai rajouté un système d'initialisation des
domaines, extensible par greffons au niveau de chaque catégorie. Il
permet typiquement de faire des déclarations de conversions implicites
ou de surcharge pour les opérateurs adéquats de la catégorie. J'ai
aussi introduit une notation objet plus usuelle, \Mup{objet::methode()}.

\subsection{Modélisation des domaines}
\label{section.combinat.domaines}

Dans cette section, j'explique pourquoi nous avons atteint les limites
du modèle objet de type Domaine/Catégorie de \axiom et \mupad.

Nous avons vu dans la démonstration
(section~\ref{section.combinat.demo} et en
particulier~\ref{section.combinat.demo.cristaux}) que nous souhaitions
mener simultanément des calculs sur les éléments d'un domaine et sur
les domaines eux-mêmes (un \emph{domaine} étant un ensemble muni
d'opérations, comme par exemple l'anneau des matrices $M_n(\RR)$,
l'anneau de polynômes $\CC[x,y,z]$, un groupe ou un graphe
cristallin).

Traditionnellement, ce sont les éléments des domaines que l'on
manipule, mettons les matrices dans $M_n(\RR)$. Ces éléments sont donc
modélisés par des objets. Le domaine, $M_n(\CC)$ lui n'a pas
d'existence. À la rigueur, on peut considérer que la classe des
matrices le modélise, mais typiquement dans un langage typé
statiquement, cette classe peut ne pas avoir d'existence à
l'exécution.

Dans le modèle objet d'\axiom, repris par la suite par \mupad, les
domaines prennent une existence à part entière. Ils sont modélisés par
des classes. Cela requiert plusieurs points techniques qui excluent
l'implantation de ce modèle dans la plupart des langages de
programmation classiques. Tout d'abord, les classes sont des entités
de premier niveau dans le langage. On peut les stocker dans des
variables, les passer en paramètre à des fonctions, etc. D'autre part,
ces classes sont paramétrées: pour modéliser $M_{10}(\RR)$, on
construit, à l'aide d'un foncteur générique pour les matrices carrées,
la classe des matrices $10\otimes 10$ sur le corps des réels (glissons
sur le problème de représentation des réels par des flottants). On
peut de plus créer au vol de nouvelles classes, éventuellement en
grand nombre. Enfin, les classes sont réflexives: on peut, à
l'exécution, demander leurs propriétés (ainsi, l'algorithmique pour
les matrices pourra changer selon si l'anneau de base est un corps, un
anneau principal, etc. À cet effet, le modèle introduit un deuxième
niveau de classes (les catégories), jouant en gros le rôle de classes
abstraites. Point fort du système, ces catégories associent aux
domaines des informations mathématiques.

Ce modèle fonctionne très bien en calcul formel traditionnel. Les
domaines servent à décrire la représentation concrète des éléments
(structure de données et opérations de bas niveau), tandis que les
catégories décrivent les propriétés mathématiques et les algorithmes
génériques associés. Mais nous en avons atteint les limites.

Prenons un exemple: les ordres partiels. Il y a d'un côté les petits
ordres partiels, que l'on veut décrire en stockant explicitement le
graphe sous-jacent. On les modélise donc naturellement par des objets
dans une classe «PetitOrdrePartiel». Cette classe implante une
opération «prédécesseurs immédiats(O,x)» qui, étant donné un ordre
partiel $O$ (une instance de PetitOrdrePartiel) et un de ses sommets
$x$ renvoie les successeurs immédiats de $x$ dans $O$.  La classe
PetitOrdrePartiel hérite d'une classe OrdrePartiel qui implante des
opérations génériques à partir des opérations de base; par exemple
«section initiale(O,x)» renvoie tous les éléments plus petits que $x$
dans $O$.  Maintenant se présente une autre situation: l'ensemble $P$
des partitions est lui aussi muni d'une structure d'ordre partiel (le
treillis de Young). Mais cette fois, dans notre modèle, $P$ est une
classe; elle ne peut donc pas être elle-même une instance de la classe
OrdrePartiel. On ne peut donc pas réutiliser les opérations génériques
qui y sont implantées.

Autre exemple: la résolution de surcharge, par exemple dans une
expression $a*b$ se fait généralement en considérant les domaines,
donc les classes, de $a$ et $b$. Comment doit-on procéder pour une
opération $A\otimes B$ qui prend deux domaines? Quelle est la classe
d'une classe?

Quel est le cœur du problème? Dans \axiom ou \mupad, on modélise la
relation «$x$ est un élément du domaine $X$» par «\texttt{x} est une
instance de la classe \texttt{x}». Mais alors les domaines et les
éléments ne jouent plus des rôles semblables. Introduire une double
hiérarchie n'est qu'un palliatif. Comment modéliser, ce qui est
naturel et utile mathématiquement, des ordres partiels dont les
éléments sont eux même des ordres partiels? Ou des monoides de
monoides de monoides? Comment implanter du code générique fonctionnant
quel que soit le niveau du domaine considéré?

Dans \sage (tout comme dans \magma dont il est inspiré), le modèle
choisi est différent et, je pense, plus adapté à notre besoin.  Au
fond, on en revient à la base de la programmation objet: tout concept
mathématique que l'on souhaite manipuler doit être représenté par un
objet; et cela vaut aussi pour un groupe, un ordre partiel, ou
$M_n(\CC)$. La relation «$x$ est un élément du domaine $X$» est une
relation entre deux objets. 

Cela dit, la conception d'une hiérarchie de catégories systématique
n'en est qu'à ses débuts dans \sage, et il reste de nombreuses
questions pratiques à régler. Pour commencer, le fait qu'un domaine
soit dans une certaine catégorie (mettons une algèbre de dimension
finie) donne des algorithmes génériques à la fois pour les éléments du
domaine (test d'inversibilité) et pour le domaine lui-même (calcul de
la théorie des représentations). Peut-on éviter l'établissement d'une
double hiérarchie de classes en parallèle, l'une donnant les
opérations sur les éléments et l'autre sur les domaines?  Comment
gérer, lors de la construction d'un produit cartésien $C=A\times B$
toutes les propriétés et opérations de $C$ induites par celles de $A$
et $B$ (des groupes, des classes combinatoires, etc.)?

Cette réflexion doit bien entendu se faire en lien avec des
experts. Nous renvoyons par exemple à~\cite{Hardin_Rioboo.2004} pour
des réflexions complémentaires sur la modélisation objet pour le
calcul formel, en particulier dans le cas plus contraignant d'un
langage fonctionnel à typage statique.

\section{Conception: représentations multiples}
\label{section.combinat.representations}

\subsection{Représentations multiples}

Une stratégie classique du calcul formel est d'avoir plusieurs
représentations pour les objets, et de convertir de l'une à l'autre
selon le calcul à réaliser. Typiquement, on peut représenter un
polynôme univarié comme combinaison linéaire de monômes ou par
évaluation sur des points bien choisis (transformée de Fourier), la
seconde représentation donnant un calcul de produit en temps
linéaire. 

\begin{sidewaysfigure}
  \centering
  \begin{bigcenter}
    \fbox{\includegraphics[height=0.4\textheight]{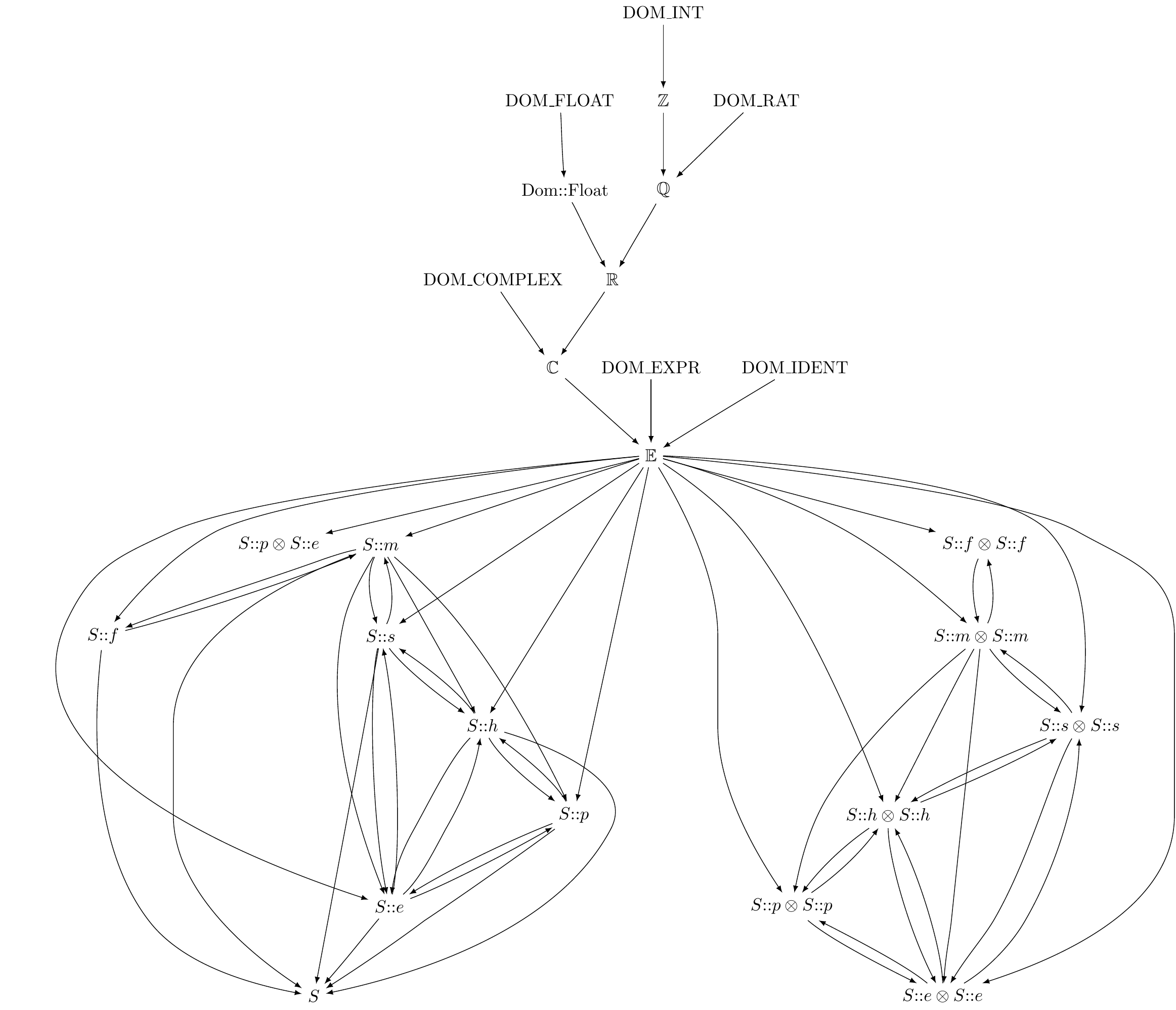}}

    \fbox{\includegraphics[width=\textwidth]{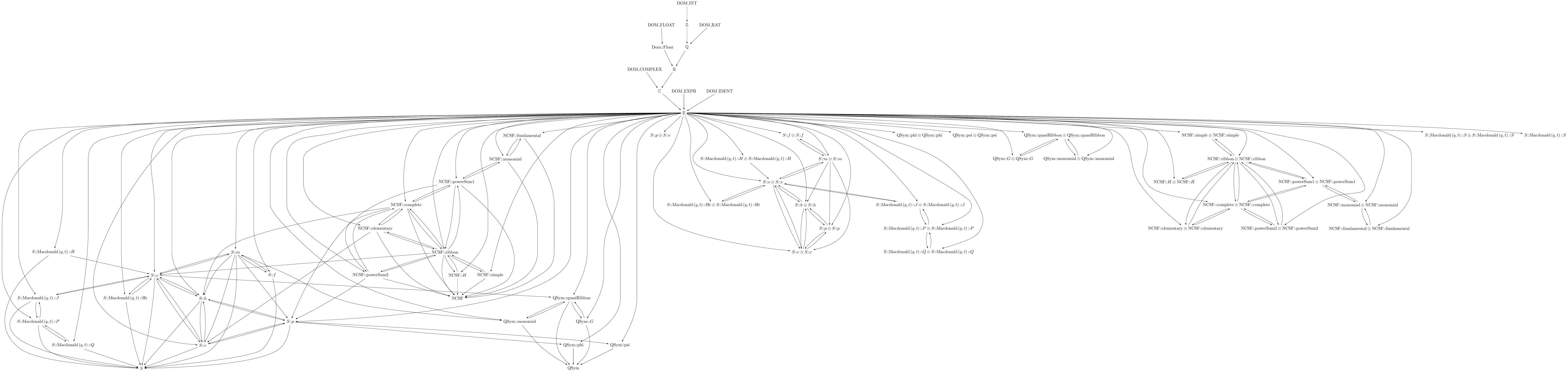}}
  \end{bigcenter}
  \caption[Graphes des conversions dans deux sessions \mupadcombinat
    typiques]{Graphes des conversions dans deux sessions \mupadcombinat
    typiques; en haut: avec juste les fonctions symétriques usuelles;
    en bas: avec en sus les polynômes de Maconald, les fonctions
    symétriques non commutatives et les fonctions quasisymétriques}
\label{figure.conversions}
\end{sidewaysfigure}

En combinatoire algébrique, cette idée est démultipliée; il y a par
exemple plus d'une vingtaine de bases implantées pour les fonctions
symétriques dans \starcombinat (voir figure~\ref{figure.conversions}),
chacune ayant son utilité propre. On veut typiquement faire des
calculs comme:
\begin{Mexin}
S := examples::SymmetricFunctions():
x := S::p( S::m[1] * ( S::e[3]*S::s[2] + 1 )) # S::e[2,1]
\end{Mexin}
\begin{Mexout}
     p[1] # e[2, 1] + 1/6 p[3, 2, 1] # e[2, 1] + 1/6 p[3, 1, 1, 1] # e[2, 1] -

        1/4 p[2, 2, 1, 1] # e[2, 1] - 1/6 p[2, 1, 1, 1, 1] # e[2, 1] + 1/12 p[1, 1, 1, 1, 1, 1] # e[2, 1]
\end{Mexout}
où \Mup{\#} est le symbole utilisé pour les produits tensoriels et,
par exemple, \Mup{S::s} modélise l'algèbre des fonctions
symétriques exprimées sur la base des fonctions de Schur.

Il est alors hors de question d'implanter toutes les conversions
possibles.  De même, il n'est pas envisageable d'implanter toutes les
opérations possibles (produit, coproduit, antipode, omega, etc) dans
toutes les bases; souvent on ne connaît pas (encore) de meilleur
algorithme que de changer de base et d'y faire le calcul! Cela devient
vital lorsque l'on souhaite explorer informatiquement une nouvelle
algèbre de Hopf ou une nouvelle base.

Il est aussi néfaste de devoir préciser, pour chaque calcul, les
conversions à utiliser. Cela requiert de la part de l'utilisateur
(novice ou expérimenté) une connaissance approfondie de ce qui est
implanté. Et surtout, cela fige le code dans un état donné, ne lui
permettant pas de bénéficier automatiquement d'implantations
ultérieures de conversions ou d'opérations.

Cette réflexion, guidée par la pratique, nous a amenés Florent Hivert
et moi à dégager l'idée clef suivante:

\emph{L'implantation d'une algèbre (de Hopf) doit être appréhendée
  comme une base de données évolutive de représentations (différentes
  bases) et de formules, règles de calculs ou autres algorithmes
  (Jacobi-Trudi, Pieri Littlewood Richardson, etc.). À charge pour le
  système de choisir, pour chaque calcul, les changements de
  représentations et les algorithmes appropriés.}

En somme, l'esprit est de modéliser l'analogue du «\emph{formulaire
  raisonné des fonctions
  symétriques}»~\cite{Lascoux_Schutzenberger.1985}, tout en
encapsulant l'ensemble des opérations et conversions effectivement
codées, pour en faire un détail d'implantation.

La réalisation de cette idée passe par trois points: un mécanisme de
conversions implicites, un mécanisme de surcharge multiparamètres, et
une infrastructure basée sur ces mécanismes pour construire aisément
des domaines avec plusieurs représentations. La difficulté était la
mise au point de modèles souples, robustes et généralistes, faisant
des choix raisonnables, prédictibles et avant tout corrects
mathématiquement, le tout avec de bonnes performances. La simplicité
était de mise. Par exemple, nous avons choisi de ne pas essayer
d'inclure dans la base de données des informations de complexité
d'algorithmes pour guider le choix du système. Il est préférable de
laisser \emph{in fine}, le contrôle à l'utilisateur avancé, lui
permettant, lorsque cela s'avère nécessaire, de rajouter quelques
règles supplémentaires (par ex: pour calculer le produit d'une
fonction de Schur par une fonction puissance, passer par telle base et
utiliser tel algorithme), plutôt que de l'inciter à ajuster des
paramètres de complexité pour corriger des travers éventuels. Cela
laisse typiquement la porte ouverte à l'utilisation de stratégies
différentes selon la taille des données, avec ajustement d'un
paramètre de seuil.

Je me suis chargé des deux premiers points, que je décris ci-dessous,
et j'ai secondé Florent Hivert sur le troisième. Nous aurions bien
préféré pouvoir transposer directement un des modèles existant dans
des langages comme \texttt{C++}, \texttt{Java}, \texttt{Python},
\texttt{Lisp}, \gap, \axiom; cependant l'ampleur de nos besoins
spécifiques nous a amené à des solutions originales, quoique largement
inspirées de ces modèles.

\subsection{Résolution de surcharge multiparamètres et conversions}

Tous les langages orientés objet, y compris \mupad, implantent la
surcharge simple des méthodes: lors d'un appel \texttt{f(x,y)}, le
choix de la fonction appelée sera déterminé par la classe de
x. Usuellement cela est mis en valeur par la syntaxe, l'appel
s'écrivant plutôt sous une forme comme \texttt{x.f(y)}.

Ce système a ses limites. Prenons un exemple dans le contexte du
calcul formel (voir~\cite{Sinha.2005.MultipleDispatch} pour des
exemples dans d'autres contextes). Supposons que \texttt{x} soit un
nombre rationnel, et \texttt{y} un vecteur de $\QQ^3$. Quelle méthode
est effectivement appelée par $x*y$? Ce devrait être celle de
multiplication des vecteurs de $\QQ^3$ par les scalaires dans $\QQ$,
dont l'implantation est rattachée à $\QQ^3$. Cependant, dans \mupad,
mais aussi dans \sage jusqu'à très récemment, c'était la méthode
\texttt{mult} des rationnels qui était appelée, à charge pour celle-ci
de renvoyer l'appel au bon endroit (l'exemple est simplifié pour
l'exposé, en fait il y avait un traitement spécifique pour les
rationnels, mais le principe reste). Cela est contraire au principe de
localisation: la méthode de multiplication des rationnels ne devrait
avoir à traiter que de la multiplication des rationnels (ou à la
rigueur des rationnels par des entiers), sans avoir à gérer et donc
dans une certaine mesure à connaître, tout ce qui peut et pourra être
construit au dessus des rationnels. En particulier, il est difficile
de s'assurer que tous les objets implanteront le renvoi de manière
consistante (ce n'était pas le cas en \mupad).

Ce n'est certainement pas un problème original. Tous les langages de
programmation implantent au moins un système minimal de surcharge
permettant d'écrire $x*y$ pour x et y entiers, flottants, etc. Le plus
souvent, cela se fait avec un traitement au cas par cas. Ce qui est
original dans notre cas, c'est l'\emph{ampleur} du problème. Il n'est
pas inusuel dans nos calculs de créer temporairement plusieurs
milliers de domaines, modélisant autant de structures algébriques dans
lesquelles ont lieu les calculs. J'avais réalisé un premier prototype,
inspiré par celui de \texttt{C++} (quoique adapté à un langage où le
typage est dynamique) ou de \gap.  Ce prototype a été
essentiel pour le développement de l'infrastructure pour les algèbres
(de Hopf) et en particulier les outils de théorie des représentations
(Florent Hivert considère par exemple que, sans ce prototype, la
publication~\cite{HNTAriki} n'aurait simplement pas existée). Il a
permis de simplifier considérablement et donc de rendre robuste et
cohérentes les conversions et les surcharges dans \mupad.

Les spécifications suivantes ont été conditionnées par le retour
d'expérience lors d'utilisations intensives. En particulier, les
considérations de complexité en temps et en mémoire proviennent toutes
de bogues rencontrés lors de gros calculs (jusqu'à 20\% du temps perdu
dans la résolution des surcharges) ou lors de l'intégration dans
\mupad (fuites de mémoire dans les tests de factorisation, ces
derniers nécessitant la construction de multiples tours d'extensions
algébriques).
\begin{itemize}
\item Les résolutions de surcharge ou de conversions sont déterminées
  par le domaine du ou des arguments.
\item Le système maintient un graphe des conversions implicites entre
  domaines existants. Chacune de ces conversions correspond à
  l'implantation d'un morphisme naturel canonique (le plongement du
  corps de base dans une algèbre, les isomorphismes entre fonctions
  symétriques dans différentes bases, etc.).

  L'aspect canonique garantit que, s'il existe deux chemins de
  conversions implicites entre deux domaines, alors les conversions
  doivent donner le même résultat (diagramme commutatif).

  L'aspect naturel garantit que la conversion est un morphisme pour
  toutes les opérations surchargées. Ce dernier point est flou: quelle
  catégorie mathématique ce graphe modélise-t-il?

  Deux domaines dans la même composante fortement connexe sont
  considérés isomorphes; c'est-à-dire qu'ils implantent deux
  représentations pour la même structure mathématique. Un résultat
  peut alors être exprimé dans l'une quelconque des représentations.

  Noter que ce graphe (voir figure~\ref{figure.conversions}) a une
  structure très spécifique: l'ensemble des prédécesseurs d'un domaine
  est relativement petit (<20), tandis que l'ensemble des successeurs
  peut représenter presque tout le graphe (pour $\NN$).

\item La résolution de surcharge, lorsqu'il n'y a pas coïncidence
  exacte des signatures, recherche la solution la plus approchée en
  autorisant l'utilisation de conversions implicites. Il n'est pas
  souhaitable de remonter dans la hiérarchie de classe, comme dans le
  modèle des multiméthodes: deux $\QQ$ algèbres $A$ et $B$ de
  dimension $4$ hériterons typiquement de la même structure
  sous-jacente d'espace vectoriel $\QQ^4$; mais cela ne donne pas pour
  autant un sens à l'addition d'un élément de $A$ et d'un de $B$.

  Lorsqu'il y a plusieurs solutions possibles, les hypothèses
  garantissent que les résultats sont identiques (éventuellement
  \emph{via} un isomorphisme). Le rôle du choix est donc uniquement de
  minimiser le nombre de conversions pour gagner en efficacité. Une
  erreur dans ce choix n'aura donc pas de conséquence sémantique.
\item Dans une résolution approchée, tous les arguments peuvent
  simultanément subir une conversion. Mais au moins l'une de ces
  conversions doit être entre deux domaines isomorphes.
\item Le mécanisme doit supporter la création au vol, éventuellement
  temporaire, de plusieurs milliers de domaines. En particulier:
  \begin{itemize}
  \item La résolution de conversions ou de surcharge ne doit avoir
    lieu que si elle est nécessaire (évaluation paresseuse). La
    première résolution ne doit considérer que la portion du graphe
    strictement nécessaire (typiquement de l'ordre de la dizaine de
    domaines). Les résolutions suivantes doivent être en temps
    constant (petit), \emph{via} l'utilisation d'un cache.
  \item Le cache mémoire ne doit pas dépasser un espace
    essentiellement linéaire vis-à-vis du nombre de domaines et
    d'opérations effectivement utilisés.
  \item La représentation en mémoire du graphe et du cache ne doivent
    pas empêcher la libération de la mémoire occupée par un domaine
    lorsque celui-ci n'est plus utilisé (que ce soit par comptage de
    référence ou par glanage de cellules lorsqu'il fait partie d'une
    composante connexe non triviale).
  \item Il n'est pas essentiel que le rajout (ou la suppression) a
    posteriori d'une conversion invalide le cache pour les conversions
    et surcharges indirectes déjà résolues.
  \end{itemize}
\end{itemize}
L'implantation actuelle dans \mupad vérifie ces spécifications grâce à
une structure de données répartie, pour le graphe de conversion comme
pour le cache, et une algorithmique appropriée. Dans le même temps le
modèle suit suffisamment les mathématiques pour que son utilisation
reste intuitive, voire transparente.

\section[Objets décomposables et espèces combinatoires]{Étude de cas: objets décomposables et espèces combinatoires}
\label{section.combinat.especes}

Les classes combinatoires décomposables (ou espèces combinatoires)
sont un des moteurs de calcul essentiels pour une bibliothèque de
combinatoire (algébrique). Elles permettent de traiter génériquement
toutes les classes combinatoires que l'on peut définir en composant
des classes plus simples, typiquement récursivement. Comme nous
l'avons vu dans la section~\ref{section.combinat.demo} cela inclut
bien entendu les arbres sous toutes leur formes (plus d'une vingtaine
dans \mupadcombinat), mais aussi les partitions, les langages définis
par automate ou grammaire, etc. Pour mentionner une application plus
originale, cet outil a fait de \mupadcombinat un maillon important
d'une chaîne logicielle pour faire du test statistique de
programmes~\cite{Gouraud.2004,Denise_Gaudel_Gouraud.2004,DeGaGoLaPe06,DBLP:journals/corr/abs-cs-0606086,Oudinet.2007,rasta08cbrTR}.

La première implantation, par Paul Zimmermann, remonte à
Gaia~\cite{Zimmermann.1994,Flajolet_Zimmermann.1994}, devenue par la
suite la bibliothèque \texttt{combstruct} de \maple. Cette
bibliothèque permet le comptage et le tirage aléatoire pour les
classes décomposables non étiquetées. Par la suite, Paul Zimmermann,
Alain Denise et Isabelle Dutour ont porté cette bibliothèque vers
\mupad 1.4.2 sous le nom de \texttt{CS}~\cite{Denise_all.1998.CS}, et
surtout l'ont étendue: tirage aléatoire de structures de grande taille
grâce au comptage approché en flottant, calcul automatique de relations
de récurrence, utilisation du produit de Karatsuba sur les séries
génératrices, génération de code C autonome. Sébastien Cellier s'est
chargé sous ma direction de son adaptation à \mupad 2.0.0.

Ma première contribution a été d'utiliser l'architecture objet de
\mupadcombinat pour réduire considérablement la taille du code et
surtout le rendre réentrant: il est maintenant possible de construire
plusieurs classes décomposables simultanément. En utilisant les
fermetures et des générateurs, j'ai rajouté les fonctionnalités de
listage et d'itération (qui, pour nous, sont les plus importantes).
Enfin, j'ai intégré cette bibliothèque dans \mupadcombinat, lui
permettant d'une part de prendre une classe combinatoire quelconque de
\mupadcombinat comme brique de base, et d'autre part d'être utilisé
comme moteur de calcul interne. M'inspirant des travaux de Xavier
Molinero et Conrado
Martinez~\cite{Martinez_Molinero.2003,Martinez_Molinero.2005}, j'ai
commencé à étendre la bibliothèque pour traiter les objets
étiquetés. Florent Hivert s'est chargé de l'implantation des
constructeurs \texttt{Cycle} et \texttt{Set} en étiqueté; suite à un
\emph{quiproquo}, il a en fait introduit pour cela une technique
nouvelle (\texttt{Div}/\texttt{MultX}) qui permet en fait de traiter
des exemples nouveaux par rapport au «boxed product» de Xavier
Molinero et Conrado Martínez~\cite{Martinez_Molinero.2003}. Lors de
deux séjours croisés, Xavier Molinero a rajouté les opérations de
\emph{unrank}, tandis que j'implantais l'itération pour les
constructeurs de cycles, de mots de Lyndon et d'ensembles pour les
objets non étiquetés~\cite{Martinez_Molinero_Thiery.2006}, en
m'appuyant sur un algorithme de Joe
Sawada~\cite{Martinez_Molinero.2004,Sawada.2003}. À l'automne 2005,
j'étais dans le jury de thèse de Xavier Molinero.

Il est clairement apparu à ce moment là que la conception monolithique
de la bibliothèque atteignait ses limites, rendant fastidieuse sa
maintenance et son extension.  Aussi, lorsque invité à l'atelier
\axiom 2006 j'ai décrit les grandes étapes pour établir une
bibliothèque de combinatoire, j'ai insisté sur l'importance des objets
décomposables et de leur implantation de façon très modulaire. Pour
illustrer mon propos, j'avais écrit un petit prototype en \aldor,
proche d'un autre prototype en \texttt{C++} que j'avais écrit quelques
mois auparavant après une discussion avec Conrado Martínez.

Martin Rubey et Ralf Hemmecke ont tout de suite pris le projet en
main, marquant la naissance
d'\aldorcombinat~\cite{Hemmecke_Rubey.2006.Aldor-Combinat}. Martin
Rubey a introduit un leitmotiv qui s'est avéré brillant: suivre au
plus près la théorie des
espèces~\cite{Joyal.1981.Especes,Bergeron_Labelle_Leroux.1994.Especes}
pour unifier le cas étiqueté et non étiqueté.
Cela l'a amené à traiter aussi les cas intermédiaires (étiquetage
semi-standard à contenu fixé) afin de pouvoir générer les structures
non-étiquetées pour la substitution.

Ralf Hemmecke a, de son côté, apporté une solution logicielle
originale: modéliser la définition récursive de la classe combinatoire
par une définition récursive des classes par foncteurs
paramétrés. Cela permet de construire simultanément la structure de
données récursive des objets, et l'algorithmique récursive pour les
opérations combinatoires sur la classe. Je ne connais pas d'autre
langage de programmation permettant cela; en fait, lorsque nous en
avions discuté, il n'était même pas clair pour nous que le compilateur
d'\aldor le permette. Ce deuxième aspect, quoique fort esthétique, est
à mon sens loin d'être essentiel, et ne doit pas être vu comme un
frein pour une adaptation à d'autres langages.

Ensemble, ils ont implanté une bibliothèque extrêmement bien
documentée~\cite{Hemmecke_Rubey.2006.Aldor-Combinat}, que l'on peut
voir comme une version effective du
livre~\cite{Bergeron_Labelle_Leroux.1994.Especes}.

Cela a servi de base inestimable de travail pour Mike Hansen qui, sur
ma suggestion, s'est chargé d'une bonne partie de l'adaptation à \sage
lors de l'été 2008, sur un financement de
\texttt{Google}\footnote{Voir \url{http://blog.phasing.org/}.}. Il est
encore trop tôt pour juger de son apport personnel pour
l'algorithmique et la conception. Mais son travail est prometteur, en
particulier du point de vue de l'intégration avec la génération
incrémentale d'objets combinatoires à un isomorphisme près --- telle
qu'utilisée par Brendan McKay pour les graphes dans \texttt{Nauty} ---
implantée dans \sage par Robert Miller.

En résumé, l'implantation actuelle dans \sagecombinat des objets
décomposables (ou plus généralement des espèces combinatoires) est le
fruit de l'intervention de très nombreux contributeurs, sur le long
terme, chacun apportant sa pierre à l'édifice. C'est une illustration
exemplaire de l'importance pour ce type d'outil d'être développés sur
un modèle libre. Dans ce cadre, et grâce à la pollinisation croisée,
même le développement sur plusieurs plateformes en parallèle a été
bénéfique, chaque portage ayant été utilisé à bon escient pour essayer
de nouvelles techniques logicielles.

\clearpage
\section{Liste de publications utilisant ou concernant \mupadcombinat}
\label{section.combinat.publications}

\begin{bibunit}
\makeatletter
\renewcommand{\@bibunitname}{MuPAD-Combinat}
\makeatother
\nocite{*}   
\renewcommand{\section}[2]{}
\renewcommand{\chapter}[2]{}
\putbib[MuPAD-Combinat]
\end{bibunit}

\bibliographystyle{alpha}
\bibliography{main}

\cleardoublepage
\thispagestyle{empty}

\bigskip
\hrule
\bigskip

\textbf{Abstract:}
\smallskip

This manuscript synthesizes almost fifteen years of research in
algebraic combinatorics, in order to highlight, theme by theme, its
perspectives.

In part one, building on my thesis work, I use tools from commutative
algebra, and in particular from invariant theory, to study isomorphism
problems in combinatorics. I first consider algebras of graph
invariants in relation with Ulam's reconstruction conjecture, and
then, more generally, the age algebras of relational structures. This
raises in return structural and algorithmic problems in the invariant
theory of permutation groups.

In part two, the leitmotiv is the quest for simple yet rich
combinatorial models to describe algebraic structures and their
representations. This includes the Hecke group algebras of Coxeter
groups which I introduced and which relate to the affine Hecke
algebras, but also some finite dimensional Kac algebras in relation
with inclusions of factors, and the rational Steenrod algebras.
Beside being concrete and constructive, such combinatorial models shed
light on certain algebraic phenomena and can lead to elegant and
elementary proofs.

My favorite tool is computer exploration, and the algorithmic and
effective aspects play a major role in this manuscript. In particular,
I describe the international open source project \starcombinat which I
founded back in 2000, and whose mission is to provide an extensible
toolbox for computer exploration in algebraic combinatorics and to
foster code sharing among researchers in this area. I present specific
challenges that the development of this project raised, and the
original algorithmic, design, and development model solutions I was
led to develop.

\bigskip
\hrule
\bigskip

\textbf{Keywords:}
\smallskip

Algebraic Combinatorics~-- Computer Algebra~-- Computer exploration

Graphs and Isomorphisms~-- Groups~-- Invariants ~--
Representations

Commutative Algebras~-- Hopf and Kac Algebras~-- Iwahori-Hecke algebras

\bigskip
\hrule
\bigskip
\clearpage
\thispagestyle{empty}

\bigskip
\hrule
\bigskip

\textbf{Résumé:}
\smallskip

Ce mémoire fait la synthèse de presque quinze années de recherche en
combinatoire algébrique afin d'en dégager, thème par thème, les
perspectives.

Dans un premier volet, issu de ma thèse, j'utilise des outils
d'algèbre commutative, et notamment de théorie des invariants, pour
étudier des problèmes d'isomorphisme en combinatoire. Je m'intéresse
tout d'abord aux algèbres d'invariants de graphes en lien avec la
conjecture de reconstruction de Ulam puis, plus généralement, aux
algèbres d'âges des structures relationnelles. Cela pose en retour des
problèmes algorithmiques et structurels en théorie des invariants des
groupes de permutations.

Dans un deuxième volet, le leitmotiv est la recherche de modèles
combinatoires simples, mais riches, pour décrire des structures
algébriques et leurs représentations. Cela inclut notamment les
algèbres de Hecke groupe que j'ai associées aux groupes de Coxeter (en
lien avec les algèbres de Hecke affine), mais aussi les algèbres de
Kac de dimension finie (en lien avec les inclusions de facteurs) et
les algèbres de Steenrod. Outre un aspect concret et effectif, de tels
modèles apportent un éclairage sur certains phénomènes algébriques, et
en particulier des démonstrations élégantes et élémentaires.

Mon outil principal est l'exploration informatique. Aussi, les aspects
algorithmiques et effectifs tiennent une place particulière dans ce
mémoire. En effet, je coordonne le développement du projet logiciel
international \starcombinat depuis sa création en 2000. Sa mission est
de fournir une boîte à outils extensible pour l'exploration
informatique en combinatoire algébrique et de promouvoir la
mutualisation de code entre les chercheurs de ce domaine. Je détaille
notamment les défis particuliers rencontrés lors de son développement,
et les solutions originales que ceux-ci m'ont amené à mettre au point,
tant du point de vue de l'algorithmique que de la conception ou du
choix du modèle de développement.

\bigskip
\hrule
\bigskip

\textbf{Mots clefs:}
\smallskip

Combinatoire Algébrique~-- Calcul Formel~-- Exploration
Informatique

Graphes et Isomorphisme~-- Groupes~-- Invariants ~--
Représentations

Algèbres commutatives~-- Algèbres de Hopf et de Kac~-- Algèbres de
Iwahori-Hecke

\bigskip
\hrule
\bigskip

\end{document}